\documentclass{amsart}

\usepackage{amsfonts}
\usepackage{amsmath}
\usepackage{amsthm}
\usepackage{amssymb}
\usepackage[english]{babel}
\usepackage{graphics}
\usepackage{latexsym}
\usepackage{longtable} 
\usepackage{mathrsfs} 
\usepackage{graphicx}
\usepackage{wrapfig} 
\usepackage[pdftex,colorlinks=true,linkcolor=blue,citecolor=blue]{hyperref}
\usepackage{enumitem}
\usepackage{slashed}
\usepackage{esint}
\usepackage{mathcomp}
\usepackage{stmaryrd}
\usepackage{wasysym}
\usepackage{fancyhdr}

\newtheorem{theorem}{Theorem}
\newtheorem{corollary}[theorem]{Corollary}
\newtheorem{lemma}[theorem]{Lemma}
\newtheorem{proposition}[theorem]{Proposition}

\newtheorem{claim}[theorem]{Claim}
\newtheorem{example}[theorem]{Example}

\theoremstyle{definition}
\newtheorem{definition}[theorem]{Definition}
\newtheorem{remark}[theorem]{Remark}

\renewcommand{\d}{\mathrm{d}}

\newcommand{\mL}{\mathcal{L}}
\newcommand{\mH}{\mathcal{H}}
\newcommand{\mF}{\mathcal{F}}

\newcommand{\K}{\mathcal{K}}

\newcommand{\mM}{\mathcal{M}}
\newcommand{\mO}{\mathcal{O}}

\newcommand{\mD}{\mathcal{D}}

\newcommand{\E}{\mathcal{E}}

\newcommand{\F}{\mathrm{F}}
\newcommand{\D}{\mathrm{D}}

\renewcommand{\:}{\!:\!} 

\newcommand{\W}{\mathcal{W}}

\newcommand{\A}{\mathbb{A}}

\newcommand{\R}{\mathbb{R}}

\newcommand{\N}{\mathbb{N}}

\newcommand{\mS}{\mathbb{S}}
\newcommand{\mB}{\mathbb{B}}

\newcommand{\noi}{\noindent}
\newcommand{\ms}{\medskip}

\newcommand{\al}{\alpha}

\newcommand{\ga}{\gamma}
\newcommand{\Ga}{\Gamma}
\newcommand{\de}{\delta}
\newcommand{\De}{\Delta}
\newcommand{\e}{\varepsilon}
\newcommand{\si}{\sigma}

\newcommand{\la}{\lambda}
\newcommand{\La}{\Lambda}

\newcommand{\Om}{\Omega}
\newcommand{\om}{\omega}

\newcommand{\ze}{\zeta}

\newcommand{\av}{-\hspace{-10.5pt} \int}

\newcommand{\weak }{\, -\!\!\!\!\!-\!\!\!\!\rightharpoonup}
\newcommand{\weakstar }{ \overset{\, *_{\phantom{|}}}{{\smash{\, -\!\!\!\!-\!\!\!\!\rightharpoonup}}\, } }

\newcommand{\larrow}{\longrightarrow}
\newcommand{\ot}{\otimes}

\newcommand{\LL}{\text{\LARGE$\llcorner$}}
\newcommand{\p}{\partial}
\newcommand{\sub}{\subseteq}
\newcommand{\set}{\setminus}
\newcommand{\by}{\times}

\newcommand{\tr}{\mathrm{tr}}
\DeclareMathOperator{\sgn}{sgn} 
\newcommand{\ess}{\mathrm{ess}}
\newcommand{\dist}{\mathrm{dist}}

\renewcommand{\div}{\mathrm{div}}

\newcommand{\osc}{\mathrm{osc}}

\newcommand{\bt}{\begin{theorem}}\newcommand{\et}{\end{theorem}}
\newcommand{\bd}{\begin{definition}}\newcommand{\ed}{\end{definition}}
\newcommand{\bl}{\begin{lemma}}\newcommand{\el}{\end{lemma}}
\newcommand{\beq}{\begin{equation}}\newcommand{\eeq}{\end{equation}}
\newcommand{\bc}{\begin{claim}}\newcommand{\ec}{\end{claim}}
\newcommand{\bex}{\begin{example}}\newcommand{\eex}{\end{example}}
\newcommand{\bcor}{\begin{corollary}}\newcommand{\ecor}{\end{corollary}}
\newcommand{\bp}{\begin{proof}}\newcommand{\ep}{\end{proof}}

\newcommand{\BPL}{\medskip \noindent \textbf{Proof of Lemma} }
\newcommand{\BPC}{\medskip \noindent \textbf{Proof of Claim} }
\newcommand{\BPCOR}{\medskip \noindent \textbf{Proof of Corollary} }
\newcommand{\BPP}{\medskip \noindent \textbf{Proof of Proposition} }
\newcommand{\BPT}{\medskip \noindent \textbf{Proof of Theorem} }

\numberwithin{equation}{section}

\makeatletter
\@namedef{subjclassname@2020}{\textup{2020} Mathematics Subject Classification}
\makeatother

\begin{document}

\title[Local minimisers in Calculus of Variations in $\mathrm L^{\infty}$]{Existence, uniqueness and characterisation of local minimisers in higher order Calculus\\ of Variations in $\mathrm L^{\infty}$}

\author{Nikos Katzourakis}

\address{N.K., Department of Mathematics and Statistics, University of Reading, Whiteknights Campus, Pepper Lane, Reading RG6 6AX, UNITED KINGDOM}

\email{n.katzourakis@reading.ac.uk}

\author{Roger Moser}

\address{R.M., Department of Mathematical Sciences,
University of Bath, Bath BA2 7AY, UNITED KINGDOM}

\email{r.moser@bath.ac.uk}


\date{}


\thanks{N.K.\ has been partially financially supported through the EPSRC grant EP/X017206/1.\ R.M.\ has been partially financially supported through the EPSRC grant EP/X017109/1.}

\begin{abstract} We study variational problems for second order supremal functionals $\mathrm F_\infty(u)= \|F(\cdot,u,\D u,\mathrm{A}\!:\!\mathrm D^2u)\|_{\mathrm L^{\infty}(\Omega)}$, where $F$ satisfies certain natural assumptions, $\mathrm A$ is a positive matrix, and $\Omega \Subset \R^n$. Higher order problems are very novel in the Calculus of Variations in $\mathrm L^{\infty}$, and exhibit a strikingly different behaviour compared to first order problems, for which there exists an established theory, pioneered by Aronsson in 1960s.  The aim of this paper is to develop a complete theory for $\mathrm F_\infty$. We prove that, under appropriate conditions, ``localised" minimisers can be characterised as solutions to a nonlinear system of PDEs, which is different from the corresponding Aronsson equation for $\mathrm F_\infty$; the latter is only a necessary, but not a sufficient condition for minimality. We also establish the existence and uniqueness of localised minimisers subject to Dirichlet conditions on $\partial \Omega$, and also their partial regularity outside a singular set of codimension one, which may be non-empty even if $n=1$.
\end{abstract}

\maketitle

\tableofcontents

\section{Introduction and main results}   \label{Section1}

Let $\Om \Subset \R^n$ be a bounded open set with Lipschitz boundary $\p\Om$, $n\in\N$. Suppose also that $F : \Om \by (\R \by \R^n \by \R) \larrow \R$ is a given Carath\'eodory function, whose arguments will be denoted by $(x,\eta,\mathrm p,\xi)$. In this paper we are interested in studying $\mathrm L^{\infty}$ variational problems for supremal functionals involving higher order derivatives of admissible functions. Specifically, we will consider the second order functional
\beq
\label{1.1}
\ \ \ \ \F_\infty (u) := \big\| F(\cdot,u, \D u,\mathrm{A}\!:\!\mathrm D^2u)\big\|_{\mathrm L^{\infty}(\Omega)},\ \ \ \ u \in \W^{2,\infty}(\Om),
\eeq
where the class $\W^{2,\infty}(\Om)$ of weakly twice differentiable functions is defined as
\beq
\label{1.2}
\W^{2,\infty}(\Om) \,:=\, \bigcap_{1<p<\infty}\Big\{ u \in\mathrm W^{2,p}(\Om)\ \big|\  \mathrm{A}\!:\! \D^2u \in \mathrm L^{\infty}(\Om) \Big\} .
\eeq
In \eqref{1.1}-\eqref{1.2}, $\mathrm A \in \R^{n\by n}$ denotes a positive symmetric matrix, and ``$:$" denotes the Euclidean inner product in  $\R^{n\by n}$ (for details on all the symbolisations we employ in this paper, we refer to Subsection \ref{subsection2.1}). Variational problems for second (or higher) order supremal functionals are very novel in the Calculus of Variations in $\mathrm L^{\infty}$, and no general theory is yet available. {\it The purpose of the present paper is to develop a complete theory for the general second order functional \eqref{1.1}}. 

Since $\F_\infty$ is not Gateaux differentiable, general critical points cannot be defined in the classical sense. The lack of smoothness is a standard intricacy of all supremal functionals, and persists even if e.g.\ $F$ is $\mathrm{C}^\infty$ and we consider only the restriction $\smash{\F_\infty|_{\mathrm{C}^2(\overline{\Om})}}$: by Danskin's theorem (see \cite{Da}), the Gateaux derivative exists only when the argmax  set is a singleton. We therefore turn our attention to the {\it existence, uniqueness and regularity of appropriately defined minimisers, as well in their variational characterisation through PDEs.}

Before introducing our assumptions and results, and for the sake of comparison, let us discuss briefly the more classical case of first order $\mathrm L^{\infty}$ functionals, which presents some striking differences. In this case, at least when restricted to scalar-valued admissible functions, there exists a rather complete theory. Let $\mathrm H : \Om \by \R \by \R^n \larrow \R$ be given, and consider
\beq
\label{1.1A}
\ \ \ \mathrm E_\infty (u,\mO) :=  \big\| \mathrm H(\cdot,u, \D u)\big\|_{\mathrm L^{\infty}(\Omega)},\ \ \ \ u \in\mathrm W^{1,\infty}(\Om), \ \mO \Subset \Om.
\eeq
Variational problems for \eqref{1.1A} were pioneered by Aronsson in the 1960s, who is considered the founder of the theory (see e.g.\ \cite{A1}-\cite{A4}). The trademark feature of \eqref{1.1A} is that minimisers of $\mathrm E_\infty (\cdot ,\Om)$ are in general highly non-unique and may not satisfy any PDE, in the sense of some analogue of the Euler-Lagrange equations of integral functionals. If however one builds locality in the minimality notion, namely require \emph{by definition} that $\mathrm E_\infty (\cdot ,\mO)$ is minimised for all $ \mO \Subset \Om$ with respect to the function's own boundary values on $\p\mO$, then any such  function $u$, satisfies the  \emph{Aronsson equation} in the viscosity sense and is called an \emph{absolute minimiser} of \eqref{1.1A} (for a pedagogical introduction, we refer e.g.\ to \cite{K0}). The Aronsson equation arising from \eqref{1.1A} is the following quasilinear second order degenerate elliptic PDE
\beq
\label{1.1B}
\mathrm H_{\mathrm p} (\cdot,u, \D u) \D \big(\mathrm H(\cdot,u, \D u)\big) =0, \ \ \text{ in }\Om.
\eeq
It turns out that the opposite is true as well, and the PDE \eqref{1.1B} is sufficient for absolute minimality, whilst the corresponding Dirichlet problem has a unique solution. All these deep results require appropriate assumptions to be imposed, typically involving convexity of $\mathrm H$ in $\mathrm P$, and no explicit dependence in $u$. Without any attempt to be exhaustive, for some important contributions in the area, we refer e.g.\ to \cite{AP, BJW2, BDP, BK, CDP, JWY, KZ, MWZ, PWZ, PZ, RZ}.

The phenomena arising in the study of higher order supremal functionals are radically different from those governing the first order case of \eqref{1.1A}-\eqref{1.1B}. Firstly, at least for the case of \eqref{1.1} where the dependence on second derivatives is through the elliptic operator $u\mapsto \mathrm A \: \D^2u$, absolute minimisers are not a fundamental notion, as \emph{usual (global) minimisers of \eqref{1.1} satisfy astounding uniqueness properties}, even with $u$-dependence being allowed and under very weak convexity requirements. In fact, one can take this one step further, and show that even \emph{localised minimisers}, namely those which minimise only in an appropriate neighbourhood of $\W^{2,\infty}(\Om)$, are also locally unique. On the other hand, the corresponding higher order Aronsson equation for \eqref{1.1}, although it can be easily derived, it plays only an ancillary role in this setting, as it \emph{is a necessary but not a sufficient condition for minimality}, even if we restrict our attention to classical solutions. Additionally, it is fully nonlinear and of third order, failing to be elliptic in any reasonable sense. Furthermore, the very definition of generalised solutions for the Aronsson equation is a non-trivial matter, as neither weak nor viscosity solutions apply. 

In this paper, following an entirely different approach, \emph{we identify a different set of equations which completely characterises minimisers of \eqref{1.1}, even localised ones, and it allows to ascertain their uniqueness for given prescribed Dirichlet boundary conditions}, under natural assumptions. We may also include explicit $(u,\D u)$-dependence, a fact which is interesting in itself, as $u$-dependence typically precludes uniqueness assertions in the first order case (see e.g.\ \cite{Y}). Our investigations are in part inspired by curvature problems in Riemannian geometry, specifically the Yamabe and Nirenberg problems, see \cite{MS,S} and subsequently. However, our results herein have a different focus and as they stand they do not directly apply to differential geometry.

We now present the main contributions of this work, beginning with the minima- lity notions. To state them concisely, let us define the following subspace of \eqref{1.2}:
\beq
\W^{2,\infty}_0(\Om)  \, := \, \W^{2,\infty}(\Om) \cap\mathrm W^{2,2}_0(\Om).
\label{1.3}
\eeq
Then, $\W^{2,\infty}_0(\Om)$ becomes a dual Banach space with separable predual when endowed with the norm (we postpone a thorough discussion until Subsection \ref{subsection2.4}):
\[
\| u \|_{\W^{2,\infty}_0(\Om)}  :=\, \|   \mathrm{A}\!:\!\D^2u \|_{\mathrm L^{\infty}(\Om)}.
\]
\begin{definition}[Global, narrow, local minimisers]
\label{definition1} Consider a function $u \in \W^{2,\infty}(\Om)$.

\smallskip 
\noi (i) We say that $u$ is a \emph{(global) minimiser} of \eqref{1.1} when
\[
 \ \F_\infty (u) \leq \F_\infty (u+\psi), \ \ \ \ \text{for all } \psi \in \W^{2,\infty}_0(\Om).
\]

\noi (ii) We say that $u$ is a \emph{narrow minimiser} of \eqref{1.1} when there exists $\rho>0$ such that
\[
 \ \ \ \ \ \ \ \ \ \ \ \ \ \ \ \ \F_\infty (u) \leq \F_\infty (u+\psi), \ \ \ \ \text{for all } \psi \in \mB^{\mathrm W^{1,\infty}(\Om)}_\rho \cap \W^{2,\infty}_0(\Om).
\]
(iii) We say that $u$ is a \emph{local minimiser} of \eqref{1.1} when there exists  $\rho>0$ such that
\[
  \ \F_\infty (u) \leq \F_\infty (u+\psi), \ \ \ \ \text{for all }\psi \in \mB^{\W^{2,\infty}_0(\Om)}_\rho.
\]
(iv) If additionally the inequalities in either of (i)-(ii)-(iii) are strict when $\psi \neq 0$, then $u$ will be called a \emph{strict global/narrow/local minimiser} respectively.
\end{definition}

One can define global/narrow/local minimisers of integral functionals in an analogous fashion. Evidently, every global minimiser is narrow, and every narrow minimiser is local, but in general neither of the converse implications are true. The significance of the concept of narrow minimisers is justified by the fact that they satisfy stronger existence, uniqueness and stability properties than local minimisers. Narrow minimisers relate to the concept of strong local minimisers (see e.g.\ \cite{KT}), but they are different. The space $\mathrm W^{1,\infty}(\Om)$ appearing in their definition is not essential, and can be replaced by any Banach space $(\mathfrak{X},|\hspace{-1pt}\| \cdot \| \hspace{-1pt}|)$ into which $ \W^{2,\infty}_0(\Om)$ has a compact embedding (e.g.\ $\mathrm L^{\infty}(\Om), \mathrm{C}^{0,\ga}(\overline{\Om}),\mathrm{C}^{1,\ga}(\overline{\Om}),\mathrm W^{2-\e,p}(\Om)$, etc.), and analogous results are obtainable in this case too. Let us now state our assumptions:
\beq
\label{1.6}
\left\{\ 
\begin{split} 
&\text{$F : \Om\by (\R\by \R^n \by \R) \larrow \R$ is a Carath\'eodory function, satisfying }
\\
&\text{that $F(x,\cdot) \in \mathrm{C}^2(\R\by\R^n \! \by \R)$, a.e.\,\,$x\in\Om$, and $F(\cdot,\eta,\mathrm p, \xi) \in \mathrm{L}^\infty(\Om)$}
\\
&\text{uniformly when $(\eta,\mathrm p, \xi)\in \mathcal M$, for any bounded $\mathcal M \sub \R\by\R^n \! \by \R$.}
\end{split}
\right.
\eeq
\beq
\label{1.7}
\left\{ \
\begin{split}
\text{There exist $C \in \mathrm{C}^\nnearrow [0,\infty)$ and  $c>0$} & \text{ such that}
\\
\p_\xi F(x,\eta,\mathrm p,\xi) \geq \ & c,  \phantom{^\big|}
\\
F(x,\eta,\mathrm p, \cdot)^{-1}(\{0\}) \sub \ & \big[\!-\!C\big(|\eta|+|\mathrm p|\big) , C\big(|\eta|+|\mathrm p|\big)\big],
\\
\big| \p_* F(x,\eta,\mathrm p,\xi) \big| \, \leq \ & C\big(|\eta|+|\mathrm p|+|\xi|\big), \ \ \ \ \ \ \ \ \ \ \ \ \ \ \ \
\\
\big| \p^2_{**} F(x,\eta,\mathrm p,\xi) \big| \, \leq \ &  C\big(|\eta|+|\mathrm p|+|\xi|\big) , \phantom{_\big|}
\\
\text{for a.e.\ $x\in\Om$ and all  $(\eta,\mathrm p, \xi)\in$ } \R\by\, & \R^n \by \R.
\end{split}
\right.
\eeq
In \eqref{1.7} we used the notation $\mathrm{C}^\nnearrow [0,\infty)$ for the set of continuous increasing functions on $[0,\infty)$, and also the following compact symbolisations for the \emph{partial} gradient $\p_* F$ and the \emph{partial} Hessian $\p^2_{**} F$  with respect to the variables $(\eta,\mathrm p, \xi) \in \R\by\R^n\by \R$:
\[
\p_* F  := \big(0,\p_\eta F, \p_{\mathrm p} F, \p_\xi F\big), \ \ \  \p^2_{**} F := \left[ 
\begin{array}{c|ccc}
0 & 0 & 0 & 0
\\
\hline
0 &  \p^2_{\eta \eta} F&  \p^2_{\eta \mathrm p}F & \p^2_{\eta \xi} F \phantom{\Big]}
\\
0 &  \p^2_{\mathrm p \eta} F&  \p^2_{\mathrm p \mathrm p}F &  \p^2_{\mathrm p \xi} F \phantom{\Big]}
\\
0 &  \p^2_{\xi \eta} F&  \p^2_{\xi \mathrm p}F &  \p^2_{\xi \xi} F \phantom{\Big]}
\end{array}
\right] .
\]
For notational convenience, let us also introduce the following symbolisations:
\[
\left\{ \ \ 
\begin{split}
\mathrm J^2u &:= \big(\cdot,u,\D u, \mathrm{A}\!:\!\D^2u\big)  :\  \Om \larrow \Om \by \R \by \R^n \by \R,
\\
\mathrm J^2_*u &:= \big(u,\D u, \mathrm{A}\!:\!\D^2u\big)\ \  : \, \, \Om \larrow \R \by \R^n \by \R.
\end{split}
\right.
\]
We will call $\mathrm J^2u$ the (second order) \emph{Jet} of $u$ and $\mathrm J^2_*u$ the (second order) \emph{reduced} \emph{Jet} of $u$. In this notation, the functional \eqref{1.1} can be compactly rewritten as
\[
\F_\infty (u) = \| F(\mathrm J^2u) \|_{\mathrm L^{\infty}(\Om)}.
\]
Finally, let $\sgn : \R \larrow [-1,1]$, be the sign function, defined by  $\sgn(t):={t}/{|t|}$ for $t\neq0$ and $\sgn(0):=0$. We may now state our first main result, which identifies a system of equations as a sufficient condition for minimality. Given $u_\infty \in \W^{2,\infty}(\Om)$ and $f_\infty \in \mathrm L^1(\Om)$, we consider the differential equations
\begin{align}
\ \ \ \ & F(\mathrm J^2u_\infty) =  \F_\infty ( u_\infty ) \sgn(f_\infty), & \text{ a.e.\ on } \Om,
\label{1.8}\\
& \div \bigg( \! \mathrm A \D f_\infty - \frac{\p_{\mathrm p} F(\mathrm J^2u_\infty)}{\p_\xi F(\mathrm J^2u_\infty)} f_\infty \! \bigg) + \frac{\p_\eta F(\mathrm J^2u_\infty)}{\p_\xi F(\mathrm J^2u_\infty)} f_\infty =\, 0 , &\ \   \text{ in } \big(\mathrm C^\infty_c(\Om)\big)^{\!*}.\!\!
\label{1.9}
\end{align}
In the above strongly coupled system, \eqref{1.8} is a second order fully nonlinear PDE satisfied in the strong a.e.\ sense by $u_\infty$, and \eqref{1.9} is a second order linear divergence elliptic  PDE for $f_\infty$, satisfied in the sense of distributions. The function $f_\infty$ should be regarded as a ``parameter" corresponding to $u_\infty$ (this will be the minimiser).

\begin{theorem}[Sufficiency of the equations] 
\label{theorem1}
Consider the functional \eqref{1.1}, where the supremand $F$ is assumed to satisfy \eqref{1.6}-\eqref{1.7}. Let $u_\infty \in \W^{2,\infty}(\Om)$ and $f_\infty \in \mathrm L^1(\Om)\set\{0\}$ be given functions which satisfy the system of equations \eqref{1.8}-\eqref{1.9}. Then, for any compact set $\K \sub \mS^{\W^{2,\infty}_0(\Om)}$, there exists $\rho >0$ such that 
\beq
\label{1.10}
\F_\infty ( u_\infty ) < \F_\infty ( u_\infty + t \psi), \ \ \ \ \ \psi \in \K,\ 0<t<\rho.
\eeq
\end{theorem}

We note that the conclusion \eqref{1.10} is weaker than local minimisers (cf.\ Definition \ref{definition1}(iii)), as it requires a compact set of directions of variations. However, since $-\K$ is also compact if $\K$ is compact, inequality \eqref{1.10} is in fact true for all $t\in (-\rho,\rho)\set\{0\}$ (perhaps for a smaller $\rho>0$). In order to obtain that $u_\infty$ is a global/narrow/local minimiser, we need some stronger hypotheses on $F$, see Theorem \ref{theorem2} that follows. Nevertheless, in Subsection \ref{subsection2.5}, we provide a necessary and sufficient condition for strong compactness in $\W^{2,\infty}_0(\Om)$, along the lines of the classical Ascoli-Arzela and Frech\'et-Kolmogorov theorems. In Corollary \ref{corollary3} (Subsection \ref{subsection2.4}) we demonstrate particular interesting choices of compact sets of variations.

We now state our second main result which gives sufficient conditions for $u_\infty$ to be a local/narrow/global minimiser for $\F_\infty$, under additional assumptions on $F$. 

\begin{theorem}[Sufficiency of the equations and strictness of miminisers] 
\label{theorem2}
In the setting of Theorem \ref{theorem1} and under the same assumptions, let $L,R \sub \R$ denote arbitrary open intervals and let $M\sub \R^n$ denote an open bounded convex set, and consider the following hypotheses:
\beq
\left\{ \ \ \ \ \ 
\begin{split}
\label{1.13}
& \text{$\exists\ L,M,R$ such that $\ess \big( \mathrm J^2_* u_\infty(\Om)\big) \! \sub L \!\by\! M \!\by\! R$, and } \exists 
\\
&\text{$\bar p \geq 2$: $|F(x,\cdot)|^{\bar p}$ is convex on $L \!\by\! M \!\by\! R$, for a.e.\ $x\in\Om$.}  
\end{split}
\right.
\eeq
\beq
\left\{ \ \ \ \ \
\begin{split}
\label{1.14}
& \text{$\exists\ L,M$ with $\ess ( (u_\infty,\D u_\infty)(\Om))\sub L\!\by\!M$, and $\forall\, R\sub \R$,} 
\\
&\text{$\exists\ \bar p \geq 2$: $|F(x,\cdot)|^{\bar p}$ is convex on $L \!\by\! M \!\by\! R$, for a.e.\,$x\in\Om$.} 
\end{split}
\right.
\eeq
\beq
\label{1.15}
 \text{$\forall\ L, M, R $, $\exists\ \bar p \geq 2$: $|F(x,\cdot)|^{\bar p}$ is convex on $L \!\by\! M \!\by\! R$, for a.e.\,$x\in\Om$.}
\eeq

Then, the following are true:

\smallskip

\noi \emph{(i)} If \eqref{1.13} is also satisfied, then $u_\infty$ is a strict local minimiser of $\F_\infty$, and hence it is unique within the neighbourhood of $\W^{2,\infty}_{u_\infty}(\Om)$ that minimises.

\smallskip

\noi \emph{(ii)} If \eqref{1.14} is also satisfied, then $u_\infty$ is a strict narrow minimiser of $\F_\infty$, and hence it is unique within the neighbourhood of $\W^{2,\infty}_{u_\infty}(\Om)$ that minimises.

\smallskip

\noi \emph{(iii)} If \eqref{1.15} is also satisfied, then $u_\infty$ is a strict global minimiser of $\F_\infty$, and hence it is globally unique in the space $\W^{2,\infty}_{u_\infty}(\Om)$. 
\end{theorem}

The meaning of the inclusions in \eqref{1.13} and \eqref{1.14} are understood in the sense of the essential range of measurable functions. In Subsection \ref{subsubsection2.2.1} we provide equivalent and sufficient conditions for the convexity of $|F(x,\cdot)|^{\bar p}$ for any $\bar p \geq 2$ (see \eqref{2.6} and  \eqref{2.7}). In particular, this requirement does \emph{not} impose convexity on $|F(x,\cdot)|$, but rather a semi-convexity type condition. In Subsection \ref{subsubsection2.2.3} we discuss explicit examples of functions which satisfy our assumptions. For convenience, we define the affine space counterpart of $\W_0^{2,\infty}(\Om)$, namely for a given $u_0 \in \W^{2,\infty}(\Om)$, we set
\beq
\W^{2,\infty}_{u_0}(\Om)  \,:= \,   u_0 + \W^{2,\infty}_0(\Om). \label{1.4}
\eeq 

Now we consider the opposite problem regarding the necessity of equations \eqref{1.8}-\eqref{1.9} for minimality, as well as the related problem of existence of minimisers, given appropriate Dirichlet boundary conditions on $\p\Om$. To this aim, we need to strengthen our assumptions on $F$. In addition to \eqref{1.6}-\eqref{1.7}, we will suppose that:
\beq
\label{1.16}
\left\{\ \
\begin{split}
& \text{Exist $c>0$ and $\al\in(0,1)$ such that,}
\\
& \hspace{100pt} \p^2_{\xi\xi} F(x,\eta,\mathrm p, \xi) \geq -\frac{1}{c},  
\\
& \hspace{22pt}  F(x,\eta,\mathrm p,\cdot)^{-1}(\{0\}) \sub \bigg[ \!-\frac{(|\eta|+|\mathrm p|)^\al+1}{c} , +\frac{(|\eta|+|\mathrm p|)^\al+1}{c} \bigg],
\\
&\text{for a.e.\ $x\in\Om$ and all $(\eta,\mathrm p, \xi)\in\R\by \R^n\by \R$.}
\end{split}
\right.
\eeq
Also, that:
\beq
\label{1.17}
\left\{\ \
\begin{split}
& \text{Exists $C \in \mathrm{C}^\nnearrow [0,\infty)$ such that,}
\\
& \hspace{50pt} \big| \p_\xi F(x,\eta,\mathrm p, \xi) \big| \, \leq\, C\big(|\eta|+|\mathrm p|\big), \ \ \ \ \ \ \ \ \ \ \ \ \ \ \ \
\\
& \hspace{50pt} \big| \p_{\mathrm p} F(x,\eta,\mathrm p, \xi) \big| \, \leq\, C\big(|\eta|+|\mathrm p|\big)(1+|\xi|), \ \ \ \ \ \ \ \ \ \ \ \ \ \ \ \
\\
& \hspace{50pt} \big| \p_\eta F(x,\eta,\mathrm p, \xi) \big| \, \leq\, C\big(|\eta|+|\mathrm p|\big)(1+|\xi|) ,
\\
&\text{for a.e.\ $x\in\Om$ and all $(\eta,\mathrm p, \xi)\in\R\by \R^n\by \R$.}
\end{split}
\right.
\eeq
Finally, we we will need to assume that the restriction $u_0|_{\p\Om}$ of the boundary data determined by $u_0 \in \W^{2,\infty}_0(\Om)$, as well as the boundary $\p\Om$ itself, are more regular:
\beq
\label{1.17A}\text{The boundary $\p\Om$ is $\mathrm{C}^2$, and the restriction $u_0|_{\p\Om}$ is $\mathrm{C}^2$ on $\p\Om$.}
\eeq

We will discuss the exact meaning of the regularity assumption \eqref{1.17A} in Subsection \ref{subsubsection2.2.2}. Below is our third main result, which concerns the necessity of equations \eqref{1.8}-\eqref{1.9} for minimality.

\begin{theorem}[Necessity of the equations]
\label{theorem6}
Suppose that $\{F,\p\Om,u_0\}$ satisfy \eqref{1.6}-\eqref{1.7}, and also \eqref{1.16}-\eqref{1.17A}. Let $u_\infty \in \W^{2,\infty}_{u_0}(\Om)$ be a narrow minimiser of \eqref{1.1}. Then, there exists $f_\infty \in \mathrm L^1(\Om) \set\{0\}$ such that the pair of functions $(u_\infty,f_\infty)$ satisfies equations \eqref{1.8}-\eqref{1.9}, and $\{f_\infty=0\}$ is countably $(n-1)$-rectifiable (and hence a Lebesgue nullset). In addition, the equation \eqref{1.9} is satisfied weakly in $(\W_0^{2,\infty}(\Om))^*$, whilst $f_\infty \in \mathrm W^{1,p}_{\mathrm{loc}}(\Om)$ for all $p\in(1,\infty)$, and the set $\{f_\infty=0\}$ is closed. 
\end{theorem}

By combining Theorems \ref{theorem2} and \ref{theorem6}, we have the following important uniqueness assertion for narrow and global minimisers (Definition \ref{definition1}).

\smallskip

\begin{theorem}[Uniqueness of narrow and global minimisers] 
\label{UniquenessTheorem} \!\!\! Suppose that $\{F,\p\Om,u_0\}$ satisfy \eqref{1.6}-\eqref{1.7}, and also \eqref{1.16}-\eqref{1.17A}. 

\smallskip

\noi \emph{(I)} If \eqref{1.14} holds true, then any narrow minimiser $u_\infty \in \W^{2,\infty}_{u_0}(\Om)$ of  \eqref{1.1} is the unique narrow minimiser within  $\smash{\mB_\rho^{\mathrm W^{1,\infty}(\Om)}(u_\infty)} \cap \W_{u_0}^{2,\infty}(\Om)$ (perhaps for a smaller $\rho>0$).

\smallskip

\noi \emph{(II)} If \eqref{1.15} holds true, then global minimisers of \eqref{1.1} in the space $\W_{u_0}^{2,\infty}(\Om)$ are unique. 
\end{theorem}

Now we consider the question of regularity of local/narrow/global minimisers, and of solutions to the system of equations \eqref{1.8}-\eqref{1.9}. In \cite[Section 8]{KP}, {\it it was proved that even when $n=1$ and no lower-order terms are present, the minimisers for $u \mapsto \|u''\|_{L^\infty(a,b)}$ in $\mathrm W^{2,\infty}_{u_0}(a,b)$ are absolute, unique, but non-$\mathrm C^2$}. More precisely, {\it minimisers are always piecewise quadratic with exactly one singular point $\bar x \in (a,b)$, and they are smooth if and only if the boundary data $\{u_0(a),u_0'(a), u_0(b),u_0'(b)\}$ can be interpolated by a  quadratic polynomial, in which case this is the minimiser}. In the generality of this work, we have the following (partial) regularity result.

\begin{corollary}[Partial regularity] \label{corollary7} Suppose that $F$ satisfies \eqref{1.6}-\eqref{1.7}, and also that $F \in \mathrm C^k\big(\Om\by \R \by \R^n \by \R\big)$, for some $k\in\N$. 

\smallskip

\noi \emph{(a)} Let $(u_\infty,f_\infty)$ $\in \W^{2,\infty}(\Om) \by (\mathrm L^1(\Om)\set\{0\})$ satisfy the equations \eqref{1.8}-\eqref{1.9}. Then,
\[
\ \ \ u_\infty \in \mathrm C^{k+2,\ga}\big(\Om \set \{f_\infty=0\} \big), \ \ \ \ \text{for all } \ga\in(0,1), 
\]
and the set $\{f_\infty=0\}$ is closed and countably $(n-1)$-rectifiable, hence in particular a Lebesgue nullset in $\Om$. 

\smallskip

\noi \emph{(b)} Under the additional assumptions of Theorem \ref{UniquenessTheorem}, any narrow/global minimiser $u_\infty$ satisfies the above $C^{k+2,\ga}$ partial regularity on $\Om \set \{f_\infty=0\}$.
\end{corollary}

We now connect the system of equations \eqref{1.8}-\eqref{1.9} to the Aronsson equation which is associated with \eqref{1.1}. By a standard formal computational argument based on $\mathrm L^p$ approximations and an appropriate rescaling of the Euler-Lagrange equations (as e.g.\ in \cite{KPr, AB}), we can derive the Aronsson equation for \eqref{1.1}, which is a fully nonlinear third order PDE:
\beq
\label{1.29}
\p_\xi F(\mathrm J^2u_\infty)\hspace{1pt} \mathrm A  \: \D\big(|F( \mathrm J^2u_\infty )| \big) \ot  \D\big(|F(\mathrm J^2u_\infty )| \big) =0, \ \ \text{ in }\Om.
\eeq
Since by assumption \eqref{1.7} we have $\p_\xi F >0$, and also $\mathrm A \in \R^{n\by n}$  is a (strictly) positive symmetric matrix, \eqref{1.29} is (at least formally) equivalent to
\[
 \D\big(|F(\mathrm J^2u_\infty )|\big) =0, \ \ \text{ in }\Om.
\]
Therefore, \emph{the Aronsson equation \eqref{1.29} corresponding to \eqref{1.1} is equivalent to the statement that $|F(\mathrm J^2u_\infty )|$ is constant on connected components of $\Om$. This is a genuinely weaker statement that the statement encoded in the system of equations \eqref{1.8}-\eqref{1.9}}. Indeed, on the one hand, \eqref{1.8} guarantees that $|F(\mathrm J^2u_\infty )|$ is not only constant throughout $\Om$, but in fact it is everywhere equal to the exact value $\F_\infty(u_\infty)$ (even if $\Om$ has more than one connected components). On the other hand, \eqref{1.9}, provides sharper information about the sign of $F(\mathrm J^2u_\infty )$, which is determined by the nodal set of the weak solution $f_\infty$ to this divergence PDE. By the results in \cite{HS}, the nodal set is actually a lower dimensional set and in particular has zero Lebesgue measure. \emph{Therefore, the Aronsson equation is only a necessary but not a sufficient condition for minimality.} More precisely, directly from Theorem \ref{theorem6}, we have that any narrow minimiser of \eqref{1.1} is a strong solution to \eqref{1.29} in \emph{contracted} form (since $|F(\mathrm J^2u_\infty )|$ is constant throughout $\Om$ and hence it is differentiable). Further, by using \cite[Th.\ 30, p.\ 665]{K1} and by arguing as in \cite[Sec.\ 7, p.\ 18-20]{KPr}, one can additionally show also that $u_\infty$ is in fact a (generalised) $\mD$-solution to the \emph{expanded} counterpart of the third order fully nonlinear Aronsson equation \eqref{1.29}.

Now we supplement Theorem \ref{theorem6} with certain approximability properties, originating from the method of construction through penalised $\mathrm L^p$-approximating functionals. Our method is based on a novel adaptation of the trademark method of $\mathrm L^p$ approximations, which consists of the introduction of an extra penalisation term. To the best of our knowledge, this adaptation first appeared in the literature in \cite{M1}. \emph{The advantage of this modification is that it forces the approximating  $\mathrm L^p$ minimisers to converge to a pre-selected $\mathrm L^\infty$ minimiser as $p\to\infty$}. For any $u_0 \in \W^{2,\infty}(\Om) $, $\bar u \in \W^{2,\infty}_{u_0}(\Om) $, $p>n$, and $\epsilon>0$, we set
\beq
\label{1.19}
\F_{p,\epsilon,\bar u}(u) \,:= \, \left(\, \av_{\Om} |F(\mathrm J^2u)|^p \, \mathrm d \mL^n \! \right)^{\!\!{1}/{p}} \,+\,  \frac{\epsilon}{2} \, {\av_{\Om}} |u-\bar u|^2 \, \mathrm d \mL^n,
\eeq
for all $u \in\mathrm W^{2,p}_{u_0}(\Om) \ (=\mathrm W^{2,p}_{\bar u}(\Om) )$.

\begin{theorem}[Necessity of the equations and $\mathrm L^p$ approximation with penalisation]
\label{theorem7}
In the setting of Theorem \ref{theorem6} and under the same hypotheses, we have that the pair $(u_\infty,f_\infty)$ satisfies the following approximation properties: 

\smallskip

\noi \emph{(i)} There exists a family of functions $\{u_{p}\}_{n<p<\infty} \sub\mathrm W^{2,n}_{u_0}(\Om)$, where for each $p$, $u_p$ is a narrow minimiser of the penalised $\mathrm L^p$ functional $\F_{p,1,u_\infty}$ (given by \eqref{1.19}) in the space $\mathrm W^{2,p}_{u_0}(\Om) $, such that, 
\beq
\label{1.20}
\left\{ \ \ 
\begin{split}
&\F_{p,1,u_\infty}(u_p) \larrow  \F_\infty(u_\infty) ,
\\
&u_p \weak  u_\infty, \ \ \text{ in }\mathrm W^{2,q}(\Om), \ \forall\ q\in (1,\infty),
\\
&u_p \larrow   u_\infty, \ \ \, \text{ in }\mathrm{C}^{1,\ga}(\overline{\Om}), \ \forall\ \ga \in (0,1),
\end{split}
\right.
\eeq
along a sequence of indices $(p_j)_1^\infty$ satisfying $p_j \to \infty$ as $j\to\infty$.

\smallskip

\noi  \emph{(ii)} For any $p \in (n,\infty)$, the narrow minimiser $u_p$ is a weak solution to the Euler-Lagrange equations of $\F_{p,1,u_\infty}$, namely
\beq
\label{1.21}
e_p^{1-p}\int_\Om \Big[ \big(|F|^{p-2}F \p_* F\big)(\mathrm J^2u_p) \cdot \mathrm J^2_* \phi \Big]  \mathrm d \mL^n  + \int_\Om  (u_p-u_\infty) \phi  \, \mathrm d \mL^n =0,
\eeq
for all $\phi \in\mathrm W^{2,p}_0(\Om)$, where
\beq
\label{1.22}
e_p \, := \left(\, \av_{\Om} |F(\mathrm J^2u_p)|^p \, \mathrm d \mL^n \! \right)^{\!\!\frac{1}{p}} .
\eeq
Further, by setting
\beq
\label{1.23}
f_p \, :=\, e_p^{1-p} \big(|F|^{p-2}F \p_\xi F\big)(\mathrm J^2u_p),
\eeq
we have
\beq
\label{1.24}
\left\{ \ \ 
\begin{split}
& f_p \mL^n \LL_\Om \weakstar f_\infty \mL^n\LL_\Om,  \ \text{ in }\mM(\overline{\Om}),
\\
& f_p   \larrow f_\infty,  \hspace{53pt}  \text{ in }\mathrm L^{q}_{\mathrm{loc}}(\Om), \ \forall\ q\in (1,n'),
\end{split}
\right.
\eeq
and also
\beq
\label{1.25}
\ \ \ \mathrm{A}\:\D^2 u_p \larrow \mathrm{A}\:\D^2u_\infty \ \  \text{ in }\mathrm L^{q}\big(\Om\set\{f_\infty=0\}\big), \ \forall\ q\in (1,\infty),
\eeq
along a sequence of indices $(p_j)_1^\infty$ which satisfies $p_j \to \infty$ as $j\to\infty$.
\end{theorem}

Our final result concerns an immediate generalisation of all our previous results to more general functionals, which may involve supremands which are perhaps non-differentiable and with arbitrary growth at infinity.

\begin{corollary}[Extension to non-differentiable supremands with arbitrary growth] \label{corollary10} Let us consider the supremal functional
\[
\ \ \ \ \ \ \ \bar \F_\infty (u) := \|\mathcal F(\mathrm J^2u)\|_{\mathrm L^{\infty}(\Omega)},\ \ \ \ u \in \W^{2,\infty}(\Om),
\]
and suppose the supremand $\mF : \Om\by (\R\by \R^n \by \R) \larrow \R$ can be factorised as follows:
\[
\left\{\ 
\begin{split} 
&\text{There exist a strictly increasing odd function $\Phi : \R \larrow \R$ and a Cara-}
\\
&\text{th\'eodory function $F : \Om\by (\R\by \R^n \by \R) \larrow \R$, such that $\mathcal F = \Phi \circ F$},
\end{split}
\right.
\]
where $F$ satisfies \eqref{1.6}-\eqref{1.7}. Then, all results hold true for $\bar \F_\infty$ as well, if $F$ satisfies the additional corresponding assumptions \eqref{1.13}, \eqref{1.14}, \eqref{1.15}, \eqref{1.16} and \eqref{1.17}.
\end{corollary}

Let us now discuss the state of the art for higher order (scalar-valued) $\mathrm L^{\infty}$ variational problems. The relevant higher order Aronsson equation was first derived formally by Aronsson himself and Barron in \cite{AB}, but noting its peculiarities, it was not explored further. In the paper \cite{KPr} by the first author and Pryer, higher order $\mathrm L^{\infty}$ problems were considered systematically for the first time, and several initial results where established for the functional $u\mapsto \| H(\D^2u)\|_{\mathrm L^{\infty}(\Om)}$, including existence of generalised solutions to the Dirichlet problem for the corresponding fully nonlinear Aronsson equation. For this task, the general framework of $\mD$-solutions was utilised, a systematic approach for fully nonlinear PDE systems, introduced by the first author in \cite{K1}. In \cite{PP}, Papamikos and Pryer studied special solutions to the same higher order Aronsson equation, classifying them by using Lie symmetries. In \cite{KP}, the first author with Parini studied a higher order $\infty$-eigenvalue problem, generalised subsequently to the fully nonlinear setting by the first author with Clark in \cite{CK2}. In \cite{CKM}, Clark and Muha, jointly with the first author, inspired by questions arising in data assimilation, studied a second order PDE constrained $\mathrm L^{\infty}$ problem involving the Navier-Stokes equations. Higher order $\mathrm L^{\infty}$ variational problems arising in Riemannian geometry were previously studied by the second author with Schwetlick and also by Sakellaris in \cite{MS, S}, and more recently again by the second author in \cite{M1}. In \cite{KM}, the authors considered the special case of \eqref{1.1} given by $u\mapsto \|\mathcal F(\cdot,\De u)\|_{\mathrm L^{\infty}(\Om)}$ which, unsurprisingly, is immensely simpler than the case of $\F_\infty$ which involves explicit $(u,\D u)$-dependence. 

This paper is organised as follows. This rather long introduction is followed by Section \ref{section2}, in which we discuss in some detail the preliminaries, our functional space setup, our assumptions, and also provide explicit classes of integrands (Subsection \ref{subsection2.2}). Therein we also establish some auxiliary results of independent interest. In particular, we identify $\W^{2,\infty}_0(\Om)$ as a dual Banach space with separable predual, and we also prove a general lemma regarding the existence of predual spaces for weakly* closed subspaces in dual Banach spaces  (Subsection \ref{subsection2.3}). We also establish certain auxiliary results of independent interest, including characterisations and examples of compact sets in $\W^{2,\infty}_0(\Om)$ (Subsection \ref{subsection2.4}), approximability of constrained minimisers in $\mathrm L^\infty$ via $\mathrm L^p$ constrained minimisers as $p\to\infty$ (Subsection \ref{subsection2.5}), and a lemma about strong compactness in $\mathrm L^1_{\mathrm{loc}}(\Om)$ of solutions to divergence PDE with $\mathrm L^1$ bounded coefficients  (Subsection \ref{subsection2.6}), and some results regarding the size of nodal sets and other properties of solutions to equations of the form \eqref{1.9} (Subsection \ref{subsection2.7}). Finally, in Section \ref{section3} we establish all our main results, whose statements are given earlier in this introduction.

We close this introduction by noting that, in our recent paper \cite{KM2}, we consider the case of the functional $u\mapsto \|  \mathrm A(\cdot) \:\D^2u +b(\cdot,u,\D u)\|_{\mathrm L^{\infty}(\Om)}$, which involves a semilinear differential operator. Even though there is some overlap of \cite{KM2} with the results obtained herein, the phenomena are significantly different. In particular, in \cite{KM2} the natural minimality notion involves \emph{almost minimisers} with a first order error term, and not global/narrow/local minimisers as herein. Additionally, unlike the present work, uniqueness in \cite{KM2} is an open question for almost minimisers. Therefore, for endogenous mathematical reasons, as well as for practicality purposes of space and size, we treat these related problems separately, as an attempt for a unified approach would considerably complicate and impair the exposition, rather than enhance it.


\section{Preliminaries and auxiliary results} 
\label{section2}

We begin by discussing in some detail our notation and our functional space setup, and also establish some auxiliary results of independent interest.


\subsection{Preliminaries and notational conventions}  \label{subsection2.1}  Let ``$:$" denotes the Euclidean (Frobenius) inner product in  $\R^{n\by n}$, namely $A\!:\!B = \tr( A^\top B)$. For $a,b\in \R^n \equiv \R^{n \by 1}$, we symbolise $a\ot b := a b^\top$. Further, for any open $\Om\sub \R^n$, for a twice (weakly) differentiable function $u : \Om \larrow \R$, the notation $\D u =(\D_1u,...,\D_nu): \Om \larrow \R^n$ will symbolise its gradient, and $\D^2 u =(\D^2_{ij}u)_{i,j=1}^n: \Om \larrow \R^{n\by n}$ will denote its Hessian, where obviously $\D_i \equiv \p/\p x_i$. For any function $u :\R^n \larrow \R^n$ and $z\in \R^n$, we define the translation $\tau_z u$ by setting $\tau_z u(x):= u(x+z)$. The closed ball, open ball and sphere of radius $r$ centred at a point $x\in \mathfrak X$ of a normed space $(\mathfrak{X},\| \cdot\|)$ will be denoted by $\bar\mB^{\mathfrak X}_\rho(x),  \mB^{\mathfrak X}_\rho(x), \mS^{\mathfrak X}_\rho(x)$ respectively. If $\rho=1$ and/or $x=0$, we will omit explicitly displaying them, and if $X$ is clear from the context, it may be omitted as well. The functional space notation we employ is either self-explanatory, or otherwise a convex combination of standard symbolisations (as e.g.\ in \cite{Ad, D, FL, KV, R, T}). For example, $[ \, \cdot \, ]_{\mathrm{C}^\al(\overline{\Om})} $ symbolises the H\"older seminorm over $\Om$, $\smash{\mathrm{C}^\nnearrow [0,\infty)}$ symbolises the set of continuous increasing functions on $[0,\infty)$, and $\smash{\mathrm{C}_{0^+}^\nnearrow [0,\infty)}$ symbolises the subset of those that vanish at the origin (i.e.\ moduli of continuity). The $n$-dimensional Lebesgue measure will be denoted by $\mL^n$, and the $k$-dimensional Hausdorff measure by $\mH^k$. Further, all $\mathrm L^p$ norms appearing in this paper with respect to the Lebesgue measure $\mL^n$ over $\Om$ will always be the normalised ones, namely with the average over $\Om$ instead of the integral, as e.g.\ in \eqref{1.19}. Let now $(\mathfrak X,\|\cdot\|)$ be a normed space, with dual Banach space $(\mathfrak X^*,\|\cdot\|_*)$. Given a subset $\mathcal A \sub \mathfrak X$ and a subset $\mathcal B\sub \mathfrak X^*$, the annihilator $\mathcal A^\bot \sub \mathfrak X^*$ and the pre-annihilator $^\bot  \mathcal B \sub \mathfrak X$ are defined as 
\[
\begin{split}
\mathcal A^\bot  &:= \big\{f\in \mathfrak X^* : f(x)=0, \ \forall \, x \in \mathcal A \big\} = \big\{f\in \mathfrak X^* : f^{-1}(\{0\}) \supseteq \mathcal A \big\}, 
\\
^\bot  \mathcal B & := \big\{x\in \mathfrak X : f(x)=0,\ \forall \, f \in \mathcal B \big\} =\bigcap \big\{ f^{-1}(\{0\}) : {f \in \mathcal  B} \big\}.
\end{split}
\]
(Details can be found e.g.\ in \cite[Ch.\ 1 and Ch.\ 4]{R}.) Then, $\mathcal A^\bot$ is a weakly*-closed subspace in $\mathfrak X^*$, and $^\bot \mathcal B $ is a norm-closed subspace in $\mathfrak X$. If $ \mathcal N\sub \mathfrak X$ is any norm-closed subspace, the equivalence relation $x' \sim x'' \Leftrightarrow x''-x'\in \mathcal N$ partitions $\mathfrak X$ to cosets, the collection of which comprises the quotient vector space $\mathfrak X/_{\mathcal N} :=\{x+{\mathcal N} : x\in \mathfrak X\}$. Let $\pi : \mathfrak X \larrow \mathfrak X/_{\mathcal N}$ denote the quotient map, given by $\pi(x):=x + {\mathcal N}$. Then, $\mathfrak X/_{\mathcal N}$ becomes a normed space, when endowed with the quotient norm 
\[
\ \ \ \|\pi(x)\|_{\mathfrak X/_{\mathcal N}}:= \inf \big\{ \|x-n\| \, :\, n\in {\mathcal N} \big\} = \dist(x,{\mathcal N}), \ \ \ x\in \mathfrak  X.
\]


\subsection{Discussion of assumptions and examples} \label{subsection2.2}

We now focus on assumptions \eqref{1.13}-\eqref{1.15}, which involve convexity of $|F(x,\cdot)|^{\bar p}$ on convex subsets of $\R \by \R^n\by \R$ for large $\bar p\geq 2$, and also on assumption \eqref{1.17A}, which requires additional regularity of the restriction of the boundary condition $u_0|_{\p\Om}$ to $\p\Om$.

\subsubsection{The restricted convexity assumptions \eqref{1.13}-\eqref{1.15}} \label{subsubsection2.2.1} A necessary and sufficient condition for the convexity of $|F(x,\cdot)|^{\bar p}$ on an open convex subset $L\by M\by R \sub \R \by \R^n \by \R$, is the following (semi-convexity type) matrix inequality
\beq
\label{2.6}
F\p^2_{**}F + (\bar p-1)\p_* F \ot  \p_* F \geq 0,
\eeq
which must hold a.e.\ on $\Om \by L \by M \by R$. The equivalence follows from the identity
\beq
\label{2.9A}
\p^2_{**}\big(|F|^p\big) = \,  p|F|^{p-2}\Big(F\p^2_{**}F  + (p-1)\p_* F \ot  \p_* F\Big).
\eeq
As a consequence, if $|F(x,\cdot)|^p$ is convex on $L\by M\by R$, so is $|F(x,\cdot)|^q$  for any $q\geq p$. In the case that there is no explicit gradient dependence (and hence no $\mathrm p$-variable), a sufficient condition for \eqref{2.6} is the following: let $\star \in \R^{2\by 2}$ be the Hodge matrix, satisfying $\star^2=-\mathrm I$ (rotation by $\pi/2$ in $\R^{2}$). Then, \eqref{2.6} follows from: 
\beq
\label{2.7}
\left\{ 
\begin{split}
& \text{$(\p_* F, \star \,\p_* F)$ is an orthogonal eigenframe  of $\p^2_{**}(|F|^2)$ with eigenvalues}
\\
& \text{$(\si_1,\si_2)$, satisfying $\si_2 \geq 0$, $\si_1 \geq - m $ for $m\in \N$, a.e.\ on $\Om\by L\by R$.} 
\end{split}
\right.
\eeq
Condition \eqref{2.7} is a stronger requirement than \eqref{2.6}, but might be easier to verify. The implication ``\eqref{2.7}$\Rightarrow$\eqref{2.6}" follows by the spectral theorem applied to $\p^2_{**}(|F|^2)$:
\[
\begin{split}
\p^2_{**}(|F|^2) &= \frac{\si_1}{|\p_* F|^2} \p_* F \ot  \p_* F \, +\, \frac{\si_2}{|\p_* F|^2} (\star \, \p_* F) \ot  (\star \, \p_* F)
\\
&\geq \frac{\si_1}{|\p_* F|^2} \p_* F \ot  \p_* F 
\\
&\geq -\frac{m}{c} \p_* F \ot  \p_* F,
\end{split}
\]
which holds a.e.\ on $\Om \by L\by R$. This implies \eqref{2.6} when $\bar p\geq (m/c)+1$. In the above inequalities we have used that $|\p_* F| \geq c >0$, which follows from \eqref{1.7}.


\subsubsection{The Whitney differentiability assumption \eqref{1.17A}} \label{subsubsection2.2.2} By the (Lipschitz) regularity of $\p\Om$, we have the Sobolev imbedding $\W^{2,\infty}(\Om) \sub \mathrm{C}^{1,\ga}(\overline{\Om})$, for all $\ga \in (0,1)$ in the usual sense, namely the precise representative from each a.e.-equivalence class has an extension to a function in $\mathrm{C}^{1,\ga}(\R^n)$. In particular, $u_0|_{\p\Om}$ is at the outset $\mathrm{C}^{1,\ga}$ on $\p\Om$, for all $\al\in(0,1)$. The assumption that $u_0|_{\p\Om}$ is $\mathrm{C}^2$ on the compact set $\p\Om$ is understood in the standard sense of Whitney \cite[Sec.\ 3,\ p.\ 64]{W}, meaning that it possesses (uniform) second order Taylor expansions. More precisely, we assume that there exists a symmetric matrix-valued map $X\in \mathrm C(\p\Om;\R^{n\by n})$ and a modulus of continuity $\om \in \mathrm{C}^\nnearrow_0[0,\infty)$, such that
\[
\left\{ \ \ \ 
\begin{split}
\left|u_0(x+z) -u_0(x) -\D u_0(x)z -\frac{1}{2}X(x) : z\ot z\right| &\leq \om(|z|)|z|^2, \ \ \ \
\\
\Big|\D u_0(x+z) -\D u_0(x) - X(x) z\Big| &\leq \om(|z|)|z|,
\end{split}
\right.
\]
for all $x,z\in\R^n$ such that $x,x+z \in \p\Om$. It can be easily seen that, due to the smoothness of $\p\Om$,  $\mathrm{C}^1$ differentiability for a scalar-valued function in the above sense is equivalent to tangential differentiability on $\p\Om$ (regarded as a $\mathrm{C}^2$ compact manifold with the induced Riemannian metric). However, requiring a function to be $\mathrm{C}^2$ on $\p\Om$, or requiring its Euclidean gradient vector field be tangentially $\mathrm{C}^1$ on $\p\Om$, are in general slightly stronger requirements that being tangentially $\mathrm{C}^2$ on $\p\Om$.


\subsubsection{Explicit classes of supremands} \label{subsubsection2.2.3} We demonstrate some classes of functions $F$ which satisfy our assumptions, noting that the results proved herein apply to supremands of the form $\mathcal F = \Phi \circ F$, for \emph{any} increasing odd function $\Phi : \R \larrow \R$ (Corollary \ref{corollary10}). For the sake of clarity, we assume no explicit gradient dependence in these examples (and hence no $\mathrm p$-variable). A similar (but more complicated) analysis can be performed by assuming explicit dependence on the gradient as well.

\smallskip

\noi {\bf Additive ansatz.} Let $A : \R \larrow \R$ and $a : \Om \by \R \larrow \R$ be given. We set 
\[
F(x,\eta,\xi) : = a(x,\eta)+A(\xi).
\]
If we have $A \in \mathrm{C}^2(\R)$ with $A_\xi \geq c>0$, and for all $x\in\Om$, we have $a(x,\cdot) \in \mathrm{C}^2(\R)$ with $a,a_\eta , a_{\eta\eta} \in C(\overline{\Om} \by \R)$, then \eqref{1.6}-\eqref{1.7} are readily satisfied. In particular, since $A$ is strictly increasing, $F(x,\eta,\xi)=0$ iff $\xi=A^{-1}(-a(x,\eta))$, and therefore
\[
F(x,\eta,\cdot)^{-1}(\{0\}) \sub \Big[A^{-1}\Big(-\max_{x\in\overline{\Om}}a(x,\eta)\Big) , A^{-1}\Big(-\min_{x\in \overline{\Om}}a(x,\eta)\Big) \Big].
\]
Regarding the satisfaction of \eqref{1.13}-\eqref{1.14}-\eqref{1.15}, a direct computation gives that, if we set $\mathfrak{F}:=F\p_{**}^2F+M \p_*F \ot \p_*F$ where $M\geq 1$, then we have $\mathfrak{F}\geq 0$ if and only if (suppressing the arguments $(x,\eta,\xi)$ for brevity)
\begin{align}
\tr (\mathfrak{F}) = (a+A)(A_{\xi\xi} +a_{\eta\eta}) + M\big[ (A_\xi)^2 + (a_\eta)^2 \big] &\geq 0, \label{2.4B}
\\
\det (\mathfrak{F}) = (a+A)^2A_{\xi\xi}a_{\eta\eta} + M(a+A) \big[ a_{\eta\eta}(A_\xi)^2 + A_{\xi\xi}(a_\eta)^2 \big] & \geq 0.  \label{2.5B}
\end{align}
Let us also consider the \emph{hypergraph} of the function $A^{-1}\circ (-a)$, namely
\[
\mH := \Big\{(x,\eta,\xi) \in \Om \by \R \by \R :\ \xi \geq A^{-1}\big(\! -a(x,\eta)\big) \Big\},
\]
which has the property that $F(x,\eta,\xi) \geq 0$ iff $(x,\eta,\xi) \in \mH$. Then:

\smallskip

\noi (a) Let $L,R \sub \R$ be given open intervals such that $\Om \by L\by R \sub \mH$ (the hypergraph). If for all $x\in \Om$, $a(x,\cdot)$ is convex on $L$, and $A$  is convex on $R$, then \eqref{2.4B}-\eqref{2.5B} are satisfied on $\Om \by L\by R$, for all $M>0$. Hence, \eqref{1.13} is satisfied. 

\smallskip

\noi (b) If $a$ is bounded below by $-m$ for some $m\in \R$, then $\mM:=\Om \by \R \by (A^{-1}(m),\infty )$  $\sub \mH$. If for all $x\in \Om$, $a(x,\cdot)$ is convex on $\R$, and $A$  is convex on $(A^{-1}(m),\infty )$, then \eqref{2.4B}-\eqref{2.5B} are satisfied on  $\mM$, for all $M>0$. Hence, \eqref{1.13} is satisfied. 

\smallskip

\noi (c) If for all $x\in \Om$, $a(x,\cdot)$ is convex on $\R$, and $A$  is also convex on $\R$, then \eqref{2.4B}-\eqref{2.5B} are satisfied on the entire hypergraph $\mH$, for all $M>0$. Hence, \eqref{1.13} is satisfied.

\smallskip

\noi (d) If for all $x\in \Om$, $a(x,\cdot)$ is affine on $\R$, and $A$  is also affine on $\R$, then \eqref{2.4B}-\eqref{2.5B} are satisfied on $\Om \by \R \by \R$ for all $M>0$, hence \eqref{1.15} is satisfied.

\smallskip

\noi Regarding the satisfaction of \eqref{1.16}-\eqref{1.17}, they are satisfied when $A$ is semiconvex on $\R$ ($A_{\xi\xi}$ bounded below) and $A_\xi$ is bounded above, with the growth of $a(x,\cdot)$ smaller than the growth of $A$ at infinity, namely when we have $|A(\xi)|\leq C(1+|\xi|^p)$ and $|a(x,\eta)|\leq C(1+|\eta|^q)$, for some $C>0$ and $0<q<p$.

\ms


\noi {\bf Multiplicative ansatz.} Let $A : \R \larrow \R$ and $a : \Om \by \R \larrow \R$ be given. We set 
\[
F(x,\eta,\xi) : = a(x,\eta)A(\xi).
\]
If $A \in \mathrm{C}^2(\R)$ satisfies $A_\xi \geq c>0$, and also for all $x\in\Om$, we have $a(x,\cdot) \in C^2(\R)$, with $a,a_\eta , a_{\eta\eta} \in C(\overline{\Om} \by \R)$ and $a\geq c>0$, then \eqref{1.6}-\eqref{1.7} are satisfied with $F(x,\eta,\cdot)^{-1}(\{0\})=\{0\}$. Regarding the satisfaction of \eqref{1.13}-\eqref{1.14}-\eqref{1.15}, let again $M\geq 1$. Suppressing the arguments $(x,\eta,\xi)$ for brevity, a computation gives that, if $\mathfrak{F}:=F\p_{**}^2F+M \p_*F \ot \p_*F$, then $\mathfrak{F}\geq 0$ if and only if
\begin{align}
\tr (\mathfrak{F}) &= A^2 \big[ a a_{\eta\eta} + M(a_\eta)^2 \big] + a^2 \big[ A A_{\xi\xi} + M(A_\xi)^2 \big] \geq 0, \label{2.6B}
\\
\label{2.7B} 
\begin{split}
\det (\mathfrak{F}) &= (aA)^2\bigg\{ 
 (AA_{\xi\xi})(a a_{\eta\eta}) -(A_\xi)^2(a_\eta)^2 
\\
&\ \ \ \, + M(A_\xi)^2 \big[ a a_{\eta\eta} - (a_\eta)^2 \big] + M(a_\eta)^2 \big[ A A_{\xi\xi} - (A_\xi)^2 \big]
\bigg\}  \geq 0.   
\end{split}
\end{align}
\noi (a) Let $L,R\sub \R$ be open intervals, and let us note the following elementary identities, which are valid for any $h\in \mathrm{C}^2(\R)$: $(|h|^p)''=p|h|^{p-2}[hh'' +(p-1)(h')^2]$ on $\R$, and $(\ln |h|)''=h^{-2}[hh''-(h')^2]$ on $\{h\neq 0\}$. It follows that, \eqref{2.6B}-\eqref{2.7B} are satisfied on $\Om \by L\by R$, if $\ln a(x,\cdot)$ is convex on $L$ for all $x\in\Om$, and $\ln |A|$ is convex on $R$. Therefore, we may select $\mathrm{C}^2$ functions $A:=e^B$, $a:=e^b$, where $B : R \larrow \R$ is convex and strictly increasing, and $b : \overline{\Om} \by L \larrow \R$ is convex in the second variable. We extend $A$ as a strictly increasing function in $\mathrm{C}^2(\R)$ which is linear outside an interval $[-m,m]$ for some $m\in\N$ large, and extend also $a$ as an arbitrary $\mathrm{C}^2$ function on $\Om\by \R$ bounded below. It follows that \eqref{1.13} is satisfied.

\smallskip

\noi (b) Let $A$ be affine. It can be easily confirmed that \eqref{2.6B}-\eqref{2.7B} can be deduced from
\[
(\ln a)_{\eta \eta} \geq \frac{M+1}{M}\big((\ln a)_\eta\big)^2.
\]
Then, for any interval $L\sub \R$, and for any $b \in \mathrm{C}^2(\overline{\Om} \by \R)$ which is bounded below and is strictly convex in $\eta$, we can find $\e>0$ small such that $a:=e^{\e b}$ satisfies the above inequality on $\Om \by L$. Hence, \eqref{2.6B}-\eqref{2.7B} hold true on $\Om \by L\by \R$, and therefore \eqref{1.14} is satisfied. Finally, assumptions \eqref{1.16}-\eqref{1.17} are also trivially satisfied.


\subsection{The functional space setup and duality} \label{subsection2.3} Now we expound on the functional spaces we utilise herein, and study some of their properties. We set
\begin{align}
\mathbb W_{\mathrm{loc}}^{2,\infty}(\Om) \,:=\, \bigcap_{1<p<\infty} \mathrm W_{\mathrm{loc}}^{2,p}(\Om) ,
\label{2.1}
\\
\mathbb W^{2,\infty}(\Om) \,:=\, \bigcap_{1<p<\infty}\mathrm W^{2,p}(\Om),
\label{2.2}\\
\mathbb W_0^{2,\infty}(\Om) \,:=\, \bigcap_{1<p<\infty}\mathrm W_0^{2,p}(\Om).
\label{2.3}
\end{align}
Then, we have $\mathbb W_{\mathrm{loc}}^{2,\infty}(\Om) \supset \mathbb W^{2,\infty}(\Om) \supset \mathbb W_0^{2,\infty}(\Om)$, and all three spaces are Fr\'echet spaces, topologised by the family of seminorms $\big\{ \| \cdot \|_{\mathrm W^{2,p}(\Om')} : \Om' \! \Subset \Om,\, p\in(1,\infty) \big\}$. If now $\A :\mathrm W^{2,1}(\Om)\larrow \mathrm L^1(\Om)$ is the (strong) linear differential operator defined by $\A u := \mathrm{A}\!:\!\D^2u$, the spaces defined in \eqref{1.2} and \eqref{1.3} can be written as
\begin{align}
\W^{2,\infty}(\Om) \, &=\, \mathbb W^{2,\infty}(\Om) \cap \A^{-1}\big(\mathrm L^{\infty}(\Om)\big) ,
\label{2.4}
\\
\W^{2,\infty}_0(\Om) \, &=\,  \mathbb W^{2,\infty}_0(\Om) \cap \A^{-1}\big(\mathrm L^{\infty}(\Om)\big).
\label{2.5}
\end{align}
As we will show in Subsections \ref{subsubsection2.4.1}-\ref{subsubsection2.4.2}, $\W^{2,\infty}_0(\Om)$ becomes a dual Banach space with (separable) predual $\A(\mathrm L^1(\Om))$ and convenient weak* density and compactness properties, when endowed with either of the equivalent norms 
\beq \label{norms}
\left\{\ \ \
\begin{split}
\| u \|_{\W^{2,\infty}_0(\Om)}  &:=\, \|   \mathrm{A}\!:\!\D^2u \|_{\mathrm L^{\infty}(\Om)},
\\
\| u \|'_{\W^{2,\infty}_0(\Om)}  &:=\, \|   \mathrm J_*^2 u \|_{\mathrm L^{\infty}(\Om)} .
\end{split}
\right.
\eeq


\subsubsection{Completeness and a priori estimates} \label{subsubsection2.4.1} We now study some important properties of the space $\W_0^{2,\infty}(\Om)$ which have independent interest, but are also utilised in the proofs of our main results in Section \ref{section3}.

\begin{lemma}
\label{lemma1}
The functional space \eqref{2.5} is a Banach space when endowed with either of the Lipschitz-equivalent norms defined in \eqref{norms}. Moreover, for any $p\in(1,\infty)$, there exists a constant $K=K(p)>0$ (depending also on $\Om,\mathrm A$), such that
\[
\| u \|_{\mathrm W^{2,p}(\Om)} \leq  K \| \mathrm{A}\!:\!\D^2 u\|_{\mathrm L^p(\Om)} \leq  K \| \mathrm{A}\!:\!\D^2 u\|_{\mathrm L^{\infty}(\Om)},
\]
for all $u\in\mathrm W^{2,p}_0(\Om)$. 
\end{lemma}

\BPL \ref{lemma1}. We begin by noting that the Lipschitz equivalence between the two norms in \eqref{norms} is an immediate consequence of Morrey's imbedding theorem, since $\|u\|_{\mathrm L^{\infty}(\Om)}\leq K \|u\|_{\mathrm W^{2,n}(\Om)} $ for some $K=K(\Om)>0$. Next, by the $\mathrm L^p$ estimates for the Laplacian on a bounded open set $U \Subset \R^n$ (see e.g.\ \cite[Cor.\ 9.10]{GT}), for any $p\in (1,\infty)$, there exists a $K=K(p,U)>0$ such that
\[
\| v \|_{\mathrm W^{2,p}(U)} \, \leq \, K \| \De v\|_{\mathrm L^p(U)},
\]
for all $ v \in \mathrm W^{2,p}_0(U)$. Let $\mathrm A = O^\top \! \La^2O$ be the spectral representation of $\mathrm A$, namely $O\in \mathrm O(n,\R)$ and $\La=\mathrm{diag}(\si(\sqrt{\mathrm A}))$ is a diagonal matrix with entries the eigenvalues of $\sqrt{\mathrm A}$. Then, by setting for brevity $\Ga:= \La O$ and selecting $v(x) := u(\Ga^{\top} x)$ and also $U:= \Ga^{-\!\top} \Om$, we have $\mathrm A=\Ga^\top\Ga$ and the identities $\De v( x) = \mathrm{A}\: \D^2 u(\Ga^{\top} x) $ and $\D v(x) = \Ga \hspace{1pt} \D u(\Ga^{\top}x)$, for all $x\in \Om$.  By the change of variables' formula and H\"older's inequality, there exists a new constant $K=K(p,\Om,\mathrm{A})>0$, such that 
\[
\| u \|_{\mathrm W^{2,p}(\Om)} \, \leq \, K \| \mathrm{A}\!:\!\D^2 u\|_{\mathrm L^p(\Om)}\, \leq \, K \|  u\|_{\W^{2,\infty}_0(\Om)},
\]
for all $u\in\mathrm W^{2,p}_0(\Om)$. 
If $(u_m)_1^\infty \sub \W^{2,\infty}_0(\Om) $ is a Cauchy sequence, the above implies
\[
\| u_m -u_n \|_{\mathrm W^{2,p}(\Om)} \leq \, K \big\| \mathrm{A}\: \D^2 u_m-\mathrm{A}\: \D^2 u_n \big\|_{\mathrm L^\infty(\Om)} \larrow 0,
\]
as $m,n \to \infty$. By a diagonal argument and by the completeness of $\mathrm L^{\infty}(\Om)$ and of $\mathbb{W}^{2,\infty}_0(\Om)$, there exist functions $u \in \mathbb{W}^{2,\infty}_0(\Om)$ and $f\in \mathrm L^{\infty}(\Om)$ such that $u_m \larrow u $ in $\mathrm W^{2,p}(\Om)$ for any fixed $p\in(1,\infty)$, and also $\mathrm{A}\!:\!\D^2 u_m \larrow f$ in $\mathrm L^{\infty} (\Om)$, both as $m\to\infty$. Further, for any fixed $\phi \in \mathrm{C}^\infty_c(\Om)$, we have
\[
\begin{split}
\int_\Om \phi (\mathrm{A}\!:\!\D^2 u_m) \,\d\mL^n & =  \int_\Om u_m (\mathrm{A}\!:\!\D^2 \phi)  \,\d\mL^n 
\\
&\!\!\! \larrow \int_\Om u (\mathrm{A}\!:\!\D^2  \phi )\, \d\mL^n  
\\
&= \int_\Om \phi (\mathrm{A}\!:\!\D^2 u) \,\d\mL^n,
\end{split}
\]
as $m\to\infty$, which establishes that $\mathrm{A}\!:\! \D^2 u = f \in \mathrm L^{\infty}(\Om)$. This in turn yields that $u\in \W^{2,\infty}_0(\Om)$, as claimed. \qed \ms

\begin{remark} The space $\W^{2,\infty}_0(\Om)$ is not closed in the relative Fr\'echet topology of $\mathbb{W}^{2,\infty}_0(\Om)$. Indeed, fix any function $u \in \mathbb{W}^{2,\infty}_c(\Om) \set \W^{2,\infty}(\Om)$, and let $u^\e := u * \eta^\e \in \mathrm{C}^\infty_c(\Om) \sub \W_0^{2,\infty}(\Om)$ be the standard mollifier, for $\e>0$ small (see e.g.\ \cite[Ch.\ 9.6]{KV}).  Then, $u^\e \larrow u$ in $\mathrm W^{2,p}_0(\Om)$ as $\e\to0$ for any fixed $p\in(1,\infty)$, but $u \not\in \W^{2,\infty}_0(\Om)$.
\end{remark}


\subsubsection{Identification of the predual of $\W^{2,\infty}_0(\Om)$ and weak* approximation} \label{subsubsection2.4.2} The next important functional properties have their own independent interest, but will also be utilised in Section \ref{section3}.

\begin{proposition}
\label{lemma2}
\emph{(i)} The space $\W^{2,\infty}_0(\Om)$ is a dual Banach space, and in particular
\[
\W^{2,\infty}_0(\Om) \cong \big( \A(\mathrm L^1(\Om))\big)^*,
\]
with the duality pairing given by
\[
\ \ \ \ \ \ \ \ \langle \psi , \A u\rangle := \int_\Om u (\mathrm{A}\!:\! \D^2 \psi)\,\d \mL^n, \ \ \ \psi \in \W^{2,\infty}_0(\Om),\, u\in \mathrm L^1(\Om).
\]
In the above, $\A : \mathrm L^1(\Om) \larrow (\W^{2,\infty}_0(\Om))^*$ denotes the bounded linear operator given by $u\mapsto \langle \cdot, \A u\rangle$, and we endow $\A(\mathrm L^1(\Om))$ with the natural norm induced by inclusion into $(\W^{2,\infty}_0(\Om))^*$, for which it becomes a separable Banach space, namely
\[
\| \A u \| :=  \sup\bigg\{ \int_\Om u (\mathrm{A}\!:\! \D^2 \phi)\,\d \mL^n \ \Big| \ \phi \in \W^{2,\infty}_0(\Om),\ \|\phi\|_{\W^{2,\infty}_0(\Om)}\leq 1\bigg\}.
\]

\noi \emph{(ii)} The space $\mathrm{C}^\infty_c(\Om)$ is sequentially weakly*-dense in $\W^{2,\infty}_0(\Om)$. Namely, for any $\psi \in \W^{2,\infty}_0(\Om)$, there exists $(\psi_m)_{m=1}^\infty \sub \mathrm{C}^\infty_c(\Om)$ such that, as $m\to\infty$ we have
\begin{align}
 & \phantom{aaaaaa} \mathrm{A}\!:\! \D^2 \psi_m \weakstar \mathrm{A}\!:\! \D^2 \psi,  &  \text{ in $\mathrm L^{\infty}(\Om)$}, \phantom{aaaaaaaaaaaaaa}
\\
 & \phantom{aaaaaa} \psi_m \! \weak \psi, &   \text{ in $\mathrm W_0^{2,p}(\Om)$, $\forall\ p\in (1,\infty)$}. \phantom{aaaaaaaaaaaaaa}
\end{align}
\end{proposition}

\begin{remark} Using similar techniques to those used to establish Proposition \ref{lemma2}, one can also show that
\[
\big(\W^{2,\infty}_0(\Om)\big)^* \cong \A(\mathrm{ba}(\Om)),
\]
where $\mathrm{ba}(\Om) \cong (\mathrm L^{\infty}(\Om))^*\cong (\mathrm L^1(\Om))^{**}$ is the space of real-valued bounded finitely additive measures, which are absolutely continuous with respect to the Lebesgue measure on $\Om$ (see e.g.\ \cite{FL,T}). However, the identification of the dual space of $\W^{2,\infty}_0(\Om)$ is not required for the main results, therefore we refrain from this task.
\end{remark}

The proof of Proposition \ref{lemma2} is based on the following general result of independent interest about subspaces of dual Banach spaces (Lemma \ref{lemma12}). This is perhaps well-known to the experts of duality theory, but since we could not locate a specific reference of this fact, we do provide a proof for the sake of completeness.

\begin{lemma} \label{lemma12} Let $(\mathfrak X,\| \cdot\|)$ be a Banach space, and let $(\mathfrak X^*,\| \cdot\|_*)$ be its dual Banach space. If $\mathfrak M\sub \mathfrak X^*$ is any weakly*-closed subspace, then $\mathfrak M$ is a dual Banach space itself, and in particular
\[
 (\mathfrak X/_{ ^\bot  \mathfrak M})^* \cong \, \mathfrak M.
\]
The isometric isomorphism above is given by $\phi \mapsto \phi \circ \pi$, where $\pi : \mathfrak X \larrow \mathfrak X/_{ ^\bot \mathfrak M}$ is the quotient map, and $^\bot \mathfrak M$ denotes the pre-annihilator of $\mathfrak M$ in $\mathfrak X$ (Subsection \ref{subsection2.1}).
\end{lemma}

\BPL \ref{lemma12}. We endow $\mathfrak X/_{ ^\bot \mathfrak M}$ with the quotient norm, and the spaces $\mathfrak X^*$ and $(\mathfrak X/_{ ^\bot \mathfrak M})^*$ with their natural dual norms, and define the mapping
\[
\ \ \ \pi^* \ : \ \  (\mathfrak X/_{ ^\bot \mathfrak M})^* \larrow \mathfrak X^*, \ \ \ \pi^* (\phi) := \phi \circ \pi.
\]
Then, $\pi^* $ is linear and well-defined, as for any $x\in \mathfrak X$, we have
\[
\big| \big(\pi^* (\phi)\big)(x)\big| = |\phi(\pi(x))| \leq \| \phi \|_{(\mathfrak X/_{  ^\bot \mathfrak M})^*}  \|\pi(x)\|_{\mathfrak X/_{  ^\bot \mathfrak M}} \leq \| \phi \|_{(\mathfrak X/_{  ^\bot \mathfrak M})^*} \|x\|.
\]
Hence, $\pi^* (\phi) \in \mathfrak X^*$. We will now show that
\beq \label{quotientballs}
\pi \big( \mB^{\mathfrak X}\big) = \mB^{\mathfrak X/_{  ^\bot \mathfrak M}}. 
\eeq
(This fact is claimed without proof in \cite[p. 97]{R}.) On the one hand, the inclusion $\pi ( \mB^{\mathfrak X}) \sub \mB^{\mathfrak X/_{  ^\bot \mathfrak M}}$ is an immediate consequence of the inequality $ \|\pi(x)\|_{\mathfrak X/_{  ^\bot \mathfrak M}} \leq  \|x\|$, which is valid for all $x\in \mathfrak X$. Conversely, if we have any $\pi(x)\in \mathfrak X/_{^\bot  \mathfrak M}$ which satisfies $\|\pi(x)\|_{\mathfrak X/_{  ^\bot \mathfrak M}}<1$, for any $\e>0$ small enough there exists ${m_\e \in} \,^\bot  \mathfrak M$ such that $\|x-m_\e \| < \|\pi(x)\|_{\mathfrak X/_{  ^\bot \mathfrak M}} +\e <1$. Since $\pi(x-m_\e )=\pi(x)$, it follows that there exists a representative $x':=x-m_\e$ such that $\pi(x)=\pi(x') \in \pi \big( \mB^{\mathfrak X}\big)$ and $ \|x'\|<1$. This proves the reverse inclusion $\smash{\mB^{\mathfrak X/_{  ^\bot \mathfrak M}} \sub \pi ( \mB^{\mathfrak X})}$, therefore establishing \eqref{quotientballs}. We will now show that $\pi^* $ is an isometry. In view of \eqref{quotientballs}, we have
\[
\begin{split}
\| \pi^* (\phi) \|_* & = \sup_{x\in \mathfrak X, \|x\|<1} |\phi(\pi(x))| 
=  \sup_{\pi(x)\in \mathfrak X/_{^\bot\! M}, \|\pi(x)\|_{\mathfrak X/_{  ^\bot \mathfrak M}}<1} |\phi(\pi(x))| 
=  \| \phi \|_{(\mathfrak X/_{  ^\bot \mathfrak M})^*} .
\end{split}
\]
We will now show that the image of $\pi^* $ lies in $ \mathfrak M$, namely
\[
\pi^* \big( \big(\mathfrak X/_{   ^\bot \mathfrak M}\big)^*\big) \sub  \mathfrak M.
\]
For any fixed $\phi \in (\mathfrak X/_{ ^\bot \mathfrak M})^*$, we have that $(\pi^* (\phi))(x) = \phi(\pi(x))=0$, for all ${x \in} ^\bot  \mathfrak  M$, because $\pi(x)=0$ when ${x \in} ^\bot  \mathfrak M$. This implies that $\pi^* (\phi) \in (^\bot  \mathfrak  M)^\bot$. Since $ \mathfrak M$ is by assumption a weakly* closed subspace of $\mathfrak X^*$, by \cite[Th.\ 4.7, p. 96]{R}, we have that $(^\bot \mathfrak  M)^\bot = \mathfrak M$. This implies the claimed inclusion. We finally show that $\pi^* $ is onto $ \mathfrak M$. To this aim, note that any functional $f \in  \mathfrak M$ induces a functional $\bar{f} \in ( \mathfrak X/_{ ^\bot \mathfrak M})^*$ in a natural way, by setting $\bar{f}(\pi(x)):= f(x)$, $x\in  \mathfrak X$. This is well-defined regardless of the choice of representative, because if we have $\pi(x'')=\pi(x')$, then $x''-x' \in { ^\bot \mathfrak M}$ and hence $f(x'')=f(x')$, because $f|_{  ^\bot \mathfrak M} \equiv 0$ when $f\in  \mathfrak M$. Further, by the identity
\[
|\bar{f}(\pi(x))| = |f(x)| = |f(x-z)|, \ \ \ \forall \ {z \in} { ^\bot \mathfrak M},
\]
for any $x\in \mathfrak X$, we estimate
\[
|\bar{f}(\pi(x))|  \leq \|f\|_* \inf_{{z \in} { ^\bot \mathfrak M}} \|x-z\| = \|f\|_* \|\pi(x)\|_{\mathfrak X/_{  ^\bot \mathfrak M}}.
\]
This shows that $\bar{f} \in (\mathfrak X/_{ ^\bot \mathfrak M})^*$ and $\pi^* (\bar{f})=f$, completing the proof.
\qed
\ms

\begin{remark} The assumption that the subspace $ \mathfrak M\sub  \mathfrak  X^*$ be weakly*-closed is essential, and in general it cannot be replaced by norm-closedness. For example, the space of sequences $c_0(\N)$ vanishing at infinity is norm-closed (and weakly*-dense) in $\ell^\infty(\N)$, but as it is well known, it possesses no predual space.
\end{remark}

Having procured the necessary tools, we main now establish Proposition \ref{lemma2}.

\BPP \ref{lemma2}. (i) We begin by noting that $\mathrm L^1(\Om)$ is a separable Banach space and $\A$ is a bounded linear operator, because for any $u \in \mathrm L^1(\Om)$, we have
\[
\begin{split}
\| \A u\|  =   \sup_{\psi \in \W^{2,\infty}_0(\Om), \|\psi\|_{\W^{2,\infty}_0(\Om)}\leq 1} \int_\Om u (\mathrm{A}\!:\! \D^2 \psi)\,\d \mL^n \leq \| u\|_{\mathrm L^1(\Om)}.
\end{split}
\]
Hence, the image $\A(\mathrm L^1(\Om))$ is a separable subset of $(\W^{2,\infty}_0(\Om))^*$. Thus,  $\A(\mathrm L^1(\Om))$ becomes a separable normed space when endowed with the dual norm induced by $(\W^{2,\infty}_0(\Om))^*$. Now we show that it is complete, and hence a Banach space, by showing that it is isometrically isomorphic to  a quotient of $\mathrm L^1(\Om)$ by a closed subspace. We define
\[
\mathcal{W} := \Big\{f\in \mathrm L^{\infty}(\Om) \ \Big| \ \exists \ \phi \in \W^{2,\infty}_0(\Om)\, :\, f=\mathrm{A} \: \D^2 \phi \text{ a.e.\ on }\Om\Big\} = \A\big( \W^{2,\infty}_0(\Om)\big).
\]
Then, by arguing as in the proof of Lemma \ref{lemma1}, we have that $\mathcal{W}$ is sequentially weakly* closed in $\mathrm L^{\infty}(\Om)$, and therefore it is also weakly* closed, because its predual $\mathrm L^1(\Om)$ is separable and complete, as a consequence of the Krein-Smulian theorem (we refer e.g.\ \cite[Sec.\ 12, p.\ 159]{C} or \cite{O}). By applying Lemma \ref{lemma12}, it follows that $\pi^* : \mathcal{W} \larrow \left( \mathrm L^1(\Om)/_{\! ^{\bot}{\mathcal{W}}}\right)^*$ is an isometric isomorphism, and we also have that
\[
\! ^{\bot}{\mathcal{W}} = \bigg\{ u \in \mathrm L^1(\Om) : \int_\Om u(\mathrm{A} \: \D^2 \phi)\, \mathrm d \mL^n=0, \ \forall \ \phi \in \W^{2,\infty}_0(\Om) \bigg\} = \mathrm{N}(\A).
\]
Hence, $\mathcal{W} \cong  ( \mathrm L^1(\Om)/_{\mathrm{N}(\A)})^{*}$. Next, note that the  induced map on the quotient space
\[
\A \ :\ \  \mathrm L^1(\Om)/_{\mathrm{N}(\A)} \larrow \A(\mathrm L^1(\Om)),\ \ \ \ \A(u+ \mathrm{N}(\A)):=  \A u,
\]
is a well-defined linear bijection, independent of the representative used in the coset. Further, for any fixed $u\in \mathrm L^1(\Om)$, by standard duality results (see e.g.\ \cite[Corollary 11.46]{KV}) and by the expression of the quotient norm on $\mathrm L^1(\Om)/_{\mathrm{N}(\A)}$, it follows that
\[
\begin{split}
\| \A u\| &=  \sup_{\psi \in \W^{2,\infty}_0(\Om), \, \| \mathrm{A} : \D^2 \psi \|_{\mathrm L^\infty(\Om)}\leq 1} \int_\Om u (\mathrm{A}\!:\! \D^2 \psi)\,\d \mL^n
\\
&= \sup_{f\in \mathcal{W}, \, \| f\|_{\mathrm L^{\infty}(\Om)}\leq 1} \int_\Om u f\, \mathrm d \mL^n.
\end{split}
\]
Hence,
\[
\begin{split}
\| \A u\| &=   \sup_{f\in ( \mathrm L^1(\Om)/_{\mathrm{N}(\A)})^*, \, \| f\|_{( \mathrm L^1(\Om)/_{\mathrm{N}(\A)})^*}\leq 1} \int_\Om u f\, \mathrm d \mL^n
\\
& = \big\| u+ \mathrm{N}(\A)\big\|_{\mathrm L^1(\Om)/_{\mathrm{N}(\A)}}.
\end{split}
\]
Thus, the mapping $\A : \mathrm L^1(\Om)/_{\mathrm{N}(\A)} \larrow \A(\mathrm L^1(\Om))$ is an isometric isomorphism, and therefore $\A(\mathrm L^1(\Om)) $ is a Banach  space. It follows that the induced dual map
\[
\A^*  : \ \ \big(\A(\mathrm L^1(\Om))\big)^* \larrow \big(\mathrm L^1(\Om)/_{\mathrm{N}(\A)}\big)^*,
\]
given by
\[
\ \ \ \ \langle \A^*(\psi),u \rangle :=  \langle \psi, \A u\rangle, \ \ \ \ u\in \mathrm L^1(\Om),
\]
is also an isometric isomorphism between the corresponding dual spaces. Further, the linear map
\[
\A \ : \ \ \W^{2,\infty}_0(\Om) \larrow \mathcal W, \ \ \ \phi \mapsto \A \phi,
\]
is an isometric isomorphism between Banach spaces, since $\|\A \phi\|_{\mathrm L^\infty(\Om)}=\|\phi\|_{\W^{2,\infty}_0(\Om)}$. In conclusion, we have the following isometric isomorphism
\[
(\A^*)^{-1}\circ \pi^*\circ \A \ :\ \ \W^{2,\infty}_0(\Om) \larrow \big(\A(\mathrm L^1(\Om))\big)^*,
\]
which, in view of the expression that the induced duality pairing takes, establishes the desired conclusion.

\smallskip

\noi (ii) We will use the regularisation scheme introduced and studied in \cite[Sec.\ 5]{K2}, whose definition we recall for the convenience of the reader. Since $\p\Om$ is assumed to be Lipschitz, there exists a smooth vector field $\xi \in \mathrm{C}^\infty_c(\R^n;\R^n)$ which is transversal to the boundary $\p\Om$, namely $\xi \cdot \mathrm{n} \geq \de_0$ which holds $\mH^{n-1}$-a.e.\ on $\p\Om$, for some $\de_0>0$, where $\mathrm{n}: \p\Om \larrow \mS$ is the unit outer normal vector field, and additionally we have $|\xi|\equiv1$ in an open collar around the boundary $\p\Om$. Further, there exist $\ell,\e_0>0$ such that, for all $\e \in (0,\e_0)$, we have $\dist (x+\e \ell \xi(x),\p\Om) \geq 2\e$, for $x\in\p\Om$. Then, for any $\psi\in \W^{2,\infty}_0(\Om)$, extended by zero on $\R^n \set \Om$, we set
\[
\mathrm K_\e \psi(x) := \int_{\R^n} \psi\big( x+\e \ell \xi(x) -\e y\big) \eta(y) \, \mathrm d y,
\]
where  $\e \in (0,\e_0)$ and $\eta \in \mathrm{C}^\infty_c(\mB)$ is a fixed mollifying kernel which satisfies $\eta \geq 0$ and $\|\eta\|_{\mathrm L^1(\R^n)}=1$. This is a variant of the standard mollifier, the difference being that it ``compresses" $\psi$ to a compactly supported function in $\Om$ before mollifying. This method gives rise to compactly supported smooth approximations, which converge globally on $\Om$. In \cite[Prop.\ 12]{K2} it is shown that $\mathrm K_\e v \in \mathrm{C}^\infty_c(\Om)$ and that $\mathrm K_\e \psi \larrow \psi$ in $\mathrm W^{1,q}_0(\Om)$ as $\e\to0$, for any $q\in (1,\infty)$. By iterating this result, we also have that $\mathrm K_\e \psi \larrow \psi$ in $\mathrm W^{2,q}_0(\Om)$ as $\e\to0$, for any $q\in (1,\infty)$. Hence, it remains to show that $ {\{\mathrm{A}\:\D^2(\mathrm K_\e\psi) \}_{\e\in(0,\e_0)}}$ is bounded in $\mathrm L^{\infty}(\Om)$. To this aim, we compute
\[
\begin{split}
\mathrm A \: \D^2(\mathrm K_\e \psi)(x) =& \ \big(\mathrm I +\e\ell \D \xi(x) \big)^{\!\!\top} \! \mathrm A\big(\mathrm I +\e\ell \D \xi(x) \big)
\\
&\ : \int_{\R^n} \D^2\psi\big( x+\e \ell \xi(x) -\e y\big) \eta(y) \, \mathrm d y,
\\
&\ +\e\ell\big( \mathrm{A}\: \D^2\xi(x)\big) \cdot \int_{\R^n} \D \psi\big( x+\e \ell \xi(x) -\e y\big) \eta(y) \, \mathrm d y,
\end{split}
\]
for any $x\in\Om$. Hence,
\[
\begin{split}
\mathrm A \: \D^2(\mathrm K_\e \psi)(x) = \,& \int_{\R^n} \mathrm{A}\: \D^2\psi\big( x+\e \ell \xi(x) -\e y\big) \eta(y) \, \mathrm d y
\\
& +\, \e\ell \Big( \D \xi(x)^\top \! \mathrm A + \mathrm A \D \xi(x)+\e\ell \D \xi(x)^\top \! \mathrm A\D\xi \Big) 
\\
&\ : \int_{\R^n} \D^2\psi\big( x+\e \ell \xi(x) -\e y\big) \eta(y) \, \mathrm d y
\\
& +\e\ell\big( \mathrm{A}\: \D^2\xi(x)\big) \cdot \int_{\R^n} \D \psi\big( x+\e \ell \xi(x) -\e y\big) \eta(y) \, \mathrm d y,
\end{split}
\]
for any fixed $x\in\Om$. Since $\eta \in \mathrm{C}^\infty_c(\mB)$, an integration by parts yields
\[
\begin{split}
\mathrm A \: \D^2(\mathrm K_\e \psi)(x) =&\, \int_{\R^n} \mathrm{A}\: \D^2\psi\big( x+\e \ell \xi(x) -\e y\big) \eta(y) \, \mathrm d y
\\
& -\, \ell \Big( \D \xi(x)^\top \! \mathrm A + \mathrm A \D \xi(x)+\e\ell \D \xi(x)^\top \! A\D\xi \Big) 
\\
&\ : \int_{\mB} \D\psi\big( x+\e \ell \xi(x) -\e y\big) \cdot \D\eta(y) \, \mathrm d y
\\
& +\e\ell\big( \mathrm{A}\: \D^2\xi(x)\big) \cdot \int_{\R^n} \D \psi\big( x+\e \ell \xi(x) -\e y\big) \eta(y) \, \mathrm d y,
\end{split}
\]
for any $x\in\Om$. Since $\eta, |\D \eta|$ are in $\mathrm L^1(\mB)$, by Young's inequality, we get the estimate
\[
\begin{split}
\big| \mathrm A \: \D^2(\mathrm K_\e \psi)(x)\big| \leq &\ \| \mathrm{A}\: \D^2\psi \|_{\mathrm L^{\infty}(\Om)} \, +\, \e_0\ell |\mathrm A| \| \D^2 \xi \|_{\mathrm L^{\infty}(\R^n)} \| \D \psi \|_{\mathrm L^{\infty}(\Om)}
\\
& +\, \ell |\mathrm A|\| \D \xi \|_{\mathrm L^{\infty}(\R^n)} \big( 2 +\e_0\ell \| \D \xi \|_{\mathrm L^{\infty}(\R^n)} \big) \| \D \psi \|_{\mathrm L^{\infty}(\Om)}  \| \D \eta \|_{\mathrm L^1(\R^n)},
\end{split}
\]
for any $x\in\Om$. This implies that $\mathrm A \: \D^2(\mathrm K_\e \psi)$ is indeed bounded in $\mathrm L^{\infty}(\Om)$ uniformly in $\e\in(0,\e_0)$. The conclusion ensues.
\qed
\ms


\subsection{Compact sets in $\W^{2,\infty}_0(\Om)$ and examples} \label{subsection2.4}
By using the analogue of the Ascoli-Arzela theorem in $\mathrm L^{\infty}(\Om)$ (see \cite{Ch}), and the results established in Subsection \ref{subsection2.3}, we can explicitly characterise compact sets $\K \sub  \W^{2,\infty}_0(\Om)$ (arising in Theorem \ref{theorem1}). To state this concisely, we need some additional notation. For any $f\in \mathrm L^{\infty}(\Om)$, the \emph{essential oscillation} of $f$ over a Lebesgue measurable set $E\sub \Om$ is defined as
\[
\underset{E}{\ess\,\osc} \,[f] \,:=\, \underset{x,y\in E}{\ess\sup} \big|f(x)-f(y)\big|.
\]
Then, a closed subset $\K \sub  \W^{2,\infty}_0(\Om)$ is (strongly) compact if and only if 
\beq
\label{1.11}
\inf_{\big\{ \{\Om_i\}_{i\in I} \text{ finite Borel partition of }\Om\big\}} \sup_{f\in\K} \max_{i \in I}  \underset{\Om_i}{\ess\,\osc} \,[\mathrm{A}\!:\! \D^2 f] =0. \ \ \ \ \
\eeq
The seemingly complicated condition \eqref{1.11} is equivalent to the following statement: for any $\e>0$, there exists a finite partition of $\Om$ to Borel sets $\{\Om_i\}_{i \in I}$, such that
\[
\big| \mathrm{A}\!:\! \D^2 f(x) -\mathrm{A}\!:\! \D^2 f(y)  \big| \, \leq\, \e,\ \ \ \text{a.e.\ }x,y \in \Om_i,
\]
for all functions $f\in\K$ and all indices $i \in I$. By utilising the above, we have the following immediate consequence of Theorem \ref{theorem1} regarding interesting compact sets of variations.

\begin{corollary}[Examples of compact sets of directions of test functions in $\W^{2,\infty}_0(\Om)$]
\label{corollary3}
In the setting of Theorem \ref{theorem1} and under the same assumptions, the following are true:

\smallskip

\noi \emph{(i)} For any $\theta>0$ and $p>n$, there exists $\rho>0$ such that
\[
\F_\infty ( u_\infty ) < \F_\infty ( u_\infty + \psi),
\]
among all $\psi \in \big(\mathrm W^{3,p} \cap \W^{2,\infty}_0\big)(\Om) \set \{0\}$ satisfying $\| \psi \|_{\W^{2,\infty}_0(\Om)}< \rho$, $\|\D^3 \psi \|_{\mathrm L^p(\Om)} < \theta\rho$.

\smallskip

\noi \emph{(ii)} Let $\{\Om_i : i \in I\}$ be a finite family of disjoint open subsets of $\Om$ with Lipschitz boundaries $\p\Om_i$, whose closures cover $\Om$. Then, for any $\theta>0$, there exists $\rho>0$ such that
\[
\F_\infty ( u_\infty ) < \F_\infty ( u_\infty + \psi),
\]
among all $\psi \in \W^{2,\infty}_0(\Om) \set\{0\}$ having piecewise H\"older continuous Hessians (namely $\psi |_{\Om_i} \in \mathrm{C}^{2,\al}(\overline{\Om}_i)$ for all $i \in I$), satisfying $\| \psi \|_{\W^{2,\infty}_0(\Om)} < \rho$, $\underset{i\in I}{\max} \,[\D^2 \psi ]_{\mathrm{C}^\al(\overline{\Om}_i)} <\theta\rho$.
\end{corollary} 

\BPCOR \ref{corollary3}. (i) Given $\theta >0$ and $p>n$, we define
\[
\begin{split}
\K' & \,:=\,  \Big\{ \psi \in \mS^{\W^{2,\infty}_0(\Om)} \cap\mathrm W^{3,p}(\Om) \, : \  \|\D^3 \psi \|_{\mathrm L^q(\Om)}\leq \theta \Big\}.
\end{split}
\]
In view of Theorem \ref{theorem1} (that will be established in Section \ref{section3}), it suffices to show that $\K'$ is a (strongly) compact set in $\W^{2,\infty}_0(\Om)$. Let $(\psi_m)_1^\infty$ be a sequence in $\mathcal K'$. Then, we have $ \|\mathrm A\: \D^2 \psi_m \|_{\mathrm L^{\infty}(\Om)} =1$ and $ \|\D^3 \psi_m \|_{\mathrm L^q(\Om)} \leq \theta$, for all $m\in\N$. By the compactness of the imbedding $\mathrm W^{3,p}(\Om) \Subset\mathrm W^{2,\infty}(\Om)$, there exists $\psi \in\mathrm W^{3,p}(\Om)$ such that $\psi_{m_k} \larrow \psi$ in $\mathrm W^{2,\infty}(\Om)$ along a subsequence as $k\to\infty$. In particular, $\psi_{m_k} \larrow \psi$  in $\W^{2,\infty}_0(\Om)$, as $k\to\infty$. Since $\D^3 \psi_{m_k} \weak \D^3 \psi$ in $\mathrm L^q(\Om;\R^{n^3})$ as well, the weak lower semicontinuity of $\| \cdot \|_{\mathrm L^q(\Om)}$ implies $\psi \in \mathcal K'$. 

\ms

\noi (ii)  Given $\theta >0$ and finitely many disjoint open sets $\{\Om_i : i\in I\}$ which have Lip- schitz boundaries $\p\Om_i \in \mathrm{C}^{0,1}$, and whose closures cover $\Om$, we define the set
\[
\begin{split}
\K'' & \,:=\, \mS^{\W^{2,\infty}_0(\Om)}  \bigcap \Big\{ \psi \in \underset{i \in I}{\oplus} \mathrm{C}^{2,\al}(\overline{\Om}_i) \, : \  \max_{i\in I} \, [\D^2 \psi ]_{\mathrm{C}^\al(\overline{\Om}_i)}  \leq \theta \Big\}.
\end{split}
\]
It is then easy to see that $\K''$ is a (strongly) compact set in $\W^{2,\infty}_0(\Om)$ (the proof follows similar lines to those of part (i), and therefore the details are omitted). \qed 
\ms


\subsection{$\mathrm L^p$-limits of $\mathrm L^{\infty}$ constrained minimisers}\label{subsection2.5}

The proofs of Theorems \ref{theorem6} and \ref{theorem7} are based on the following standalone result of independent interest, which establishes the existence and $\mathrm L^p$-approximation of narrow (and global) minimisers of \eqref{1.1} through penalised functionals, treated in a unified fashion as constrained minimisers.

\begin{theorem}[Existence and $\mathrm L^p$-approximation of $\mathrm L^{\infty}$ constrained minimisers]
\label{theorem8} Let $F$ satisfy \eqref{1.6}-\eqref{1.7} and \eqref{1.16}, and let $\bar u \in \W^{2,\infty}(\Om)$ be fixed. Suppose that $\{\E^p_{\bar u}(\Om)\}_{n\leq p \leq \infty}$ are classes of functions, which satisfy:
\beq
\label{1.25}
\left\{ \ \ 
\begin{split}
& \E^p_{\bar u}(\Om) \sub\mathrm W^{2,p}_{\bar u}(\Om), \text{ when } \, n\leq p<\infty,
\\
&\E^p_{\bar u}(\Om) \text{ is weakly closed, when } \, p<\infty,
\\
& \E^p_{\bar u}(\Om)  \sub  \E^q_{\bar u}(\Om), \text{ for } n\leq q \leq p \leq \infty,
\\
& \E^\infty_{\bar u}(\Om) \sub  \W^{2,\infty}_{\bar u}(\Om), \text{ when } \, p=\infty,
\\
& \E^\infty_{\bar u}(\Om) \, = \bigcap_{n<p<\infty} \E^p_{\bar u}(\Om) ,
\\
& \bar u  \in \E^\infty_{\bar u}(\Om).
\end{split}
\right.
\eeq
Then, the following are true:

\smallskip

\noi \emph{(i)} For any $\epsilon\geq 0$ and $p\geq n+ 1/c^3$, the functional $\F_{p,\epsilon, \bar u}$ (defined in \eqref{1.19}) has a minimiser $u_p$ over the admissible class $\E^p_{\bar u}(\Om) $. 

\smallskip

\noi \emph{(ii)} For any $\epsilon\geq 0$, the functional $ \F_\infty +  \frac{\epsilon}{2} \| \cdot - \bar u\|^2_{\mathrm L^2(\Om)}$ has a minimiser $u_\infty$ over the admissible class $\E^\infty_{\bar u}(\Om)$, where $\F_\infty$ is as in \eqref{1.1}. Additionally, there exists a sequence of indices $(p_j)_1^\infty$ satisfying $p_j\to\infty$ as $j\to\infty$, along which we have
\beq
\label{1.26}
\left\{ \ \
\begin{split}
& \F_{p,\epsilon,\bar u}(u_p) \larrow \F_\infty(u_\infty)  +  \frac{\epsilon}{2} \| u_\infty - \bar u\|^2_{\mathrm L^2(\Om)},
\\
&u_p \weak  u_\infty, \hspace{20pt} \text{ in }W^{2,q}(\Om), \ \forall\ q\in (1,\infty),
\\
&u_p \larrow   u_\infty, \hspace{21pt}  \text{ in }\mathrm{C}^{1,\ga}(\overline{\Om}), \ \forall\ \ga \in (0,1).
\end{split}
\right.
\eeq
\emph{(iii)} If finally $\bar u$ minimises $\F_\infty$ over the admissible class $\E^\infty_{\bar u}(\Om)$ and also $\epsilon>0$, then we have that $\bar u = u_\infty$. In this case, we also have that the modes of convergence in \eqref{1.26} are in fact full as $p\to\infty$.
\end{theorem}

The main technical step for the proof of Theorem \ref{theorem8} is the following coercivity estimate for \eqref{1.19}.

\begin{lemma}[Coercivity] \label{lemma13} Let $\de>0$ and let $w\in \E_{\bar u}^p(\Om)$ be a  $\de$-approximate minimiser of $\F_{p,\epsilon,\bar u}$ in $\E_{\bar u}^p(\Om)$, namely
\[
\F_{p,\epsilon,\bar u}(w) \, \leq \, \de + \inf_{\E_{\bar u}^p(\Om)} \F_{p,\epsilon,\bar u}.
\]
Then,  for any $q \in [n,p]$, we have the estimates
\beq
\label{2.20A}
\left\{\ \ 
\begin{split}
 \|w\|_{\mathrm W^{2,q}(\Om)}    \, &\leq \,  \frac{K(\Om,q)}{c(1-\al)}\Big( M+\de + \F_\infty(\bar u)\Big), \ \ \ &\mathrm{(i)}
\\
 \ \ \ \|\mathrm{A}\:\D^2 w\|_{\mathrm L^{q}(\Om)}    \, &\leq \, \frac{1}{c}\Big(1+\de + \F_\infty(\bar u)+ \|w\|_{\mathrm W^{1,q}(\Om)}^\al \Big), \ \ \ & \mathrm{(ii)}
\end{split}
\right.
\eeq
where $K(\Om,q)>0$ is the constant appearing in Lemma \ref{lemma1}, and $M>0$ is a constant depending only on $\bar u,\al,c,K(q,\Om),\mathrm{A}$. 
\end{lemma}

\BPL \ref{lemma13}. Note first that by \eqref{1.1}, \eqref{1.19}, minimality and Young's inequality, we have
\[
\inf_{\E_{\bar u}^p(\Om)} \F_{p,\epsilon,\bar u} \leq  \F_{p,\epsilon,\bar u}(\bar u) \leq  \F_{\infty}(\bar u).
\]
Further, by \eqref{1.7} we have that $F(x,\eta,\mathrm p,\cdot)^{-1}(\{0\})$ is a singleton set. Hence, if we have that if we symbolise 
\beq
\label{xibar}
F(x,\eta,\mathrm p, \cdot)^{-1}(\{0\}) \, =\, \big\{ \bar \xi (x,\eta,\mathrm p)\big\},
\eeq
again by \eqref{1.7}, we can easily deduce the inequality
\beq
\label{2.27A}
|F(x,\eta,\mathrm p,\xi)| \geq c \big| \xi -\bar \xi(x,\eta,\mathrm p)\big|.
\eeq
By Lemma \ref{lemma1} and \eqref{1.16}, since $w -\bar u \in\mathrm W^{2,q}_0(\Om)$, we can  estimate
\[
\begin{split}
\F_{p,\epsilon,\bar u}(w) \, &\geq \, \| F(\mathrm J^2w)\|_{\mathrm L^q(\Om)} 
\\
& \geq \, c \big\| \mathrm{A}\:\D^2 w -\bar \xi(\cdot,w,\D w)\big\|_{\mathrm L^q(\Om)} 
\\
& \geq \, c \big\| \mathrm{A}\:\D^2 (w -\bar u) \big\|_{\mathrm L^q(\Om)}  - c \| \mathrm{A}\: \D^2 \bar u \|_{\mathrm L^q(\Om)} - c\big\|\bar \xi(\cdot,w,\D w)\big\|_{\mathrm L^q(\Om)} 
\\
& \geq \, c\frac{\| w-\bar u\|_{\mathrm W^{2,q}(\Om)}}{K(q,\Om)}  - c  \| \mathrm{A}\: \D^2 \bar u \|_{\mathrm L^{\infty}(\Om)}  - \big\| \big(|w|+|\D w|\big)^\al \big\|_{\mathrm L^{q}(\Om)} -1 
\\
& \geq \, c\frac{\| w-\bar u\|_{\mathrm W^{2,q}(\Om)}}{K(q,\Om)}  - c  \| \mathrm{A}\: \D^2 \bar u \|_{\mathrm L^{\infty}(\Om)}  - \| w \|^\al_{\mathrm W^{1,\al q}(\Om)} -1 
\\
& \geq \, c\frac{\| w-\bar u\|_{\mathrm W^{2,q}(\Om)}}{K(q,\Om)}  - c  \| \mathrm{A}\: \D^2 \bar u \|_{\mathrm L^{\infty}(\Om)}  - \| w \|^\al_{\mathrm W^{1,q}(\Om)} -1.
\end{split}
\]
Therefore, 
\[
\begin{split}
\F_{p,\epsilon,\bar u}(w) \, & \geq \, c\frac{\| w-\bar u\|_{\mathrm W^{2,q}(\Om)}}{K(q,\Om)}  - c  \| \mathrm{A}\: \D^2 \bar u \|_{\mathrm L^{\infty}(\Om)}  - \| w -\bar u \|^\al_{\mathrm W^{1,q}(\Om)} 
\\
&\ \ \ \ - \| \bar u \|^\al_{\mathrm W^{1,\infty}(\Om)} -1 
\\
& \geq \, c\frac{\| w-\bar u\|_{\mathrm W^{2,q}(\Om)}}{K(q,\Om)}  - c  \| \mathrm{A}\: \D^2 \bar u \|_{\mathrm L^{\infty}(\Om)}  - \| w -\bar u \|^\al_{\mathrm W^{1,q}(\Om)} 
\\
&\ \ \ \ - \| \bar u \|^\al_{\mathrm W^{1,\infty}(\Om)} -1 
\\
& \geq \, \frac{c}{K(q,\Om)}\| w-\bar u\|_{\mathrm W^{2,q}(\Om)}  - \big(\| w -\bar u \|_{\mathrm W^{2,q}(\Om)} \big)^\al  
\\
&\ \ \ \ - \Big(c  \| \mathrm{A}\: \D^2 \bar u \|_{\mathrm L^{\infty}(\Om)} + \| \bar u \|^\al_{\mathrm W^{1,\infty}(\Om)} +1 \Big).
\end{split}
\]
By using the elementary inequality
\[
N t -t^\al \geq N(1-\al)\big[t-N^{1/(\al-1)}\big], \ \ \ t\geq 0,\  N>0,
\]
the previous estimate yields
\[
\begin{split}
\F_{p,\epsilon,\bar u}(w) \, & \geq \, \frac{c(1-\al)}{K(q,\Om)} \bigg[ \| w-\bar u\|_{\mathrm W^{2,q}(\Om)} - \Big(\frac{c}{K(q,\Om)}\Big)^{1/(\al-1)}\bigg]
\\
&\ \ \ \ - \Big(c  \| \mathrm{A}\: \D^2 \bar u \|_{\mathrm L^{\infty}(\Om)} + \| \bar u \|^\al_{\mathrm W^{1,\infty}(\Om)} +1 \Big).
\end{split}
\]
This implies that
\[
\begin{split}
\F_{p,\epsilon,\bar u}(w) \, &  \geq \, \frac{c(1-\al)}{K(q,\Om)} \| w\|_{\mathrm W^{2,q}(\Om)} 
\\
&\ \ \ \ - \bigg\{\frac{c(1-\al)}{K(q,\Om)}\bigg[\| \bar u\|_{\mathrm W^{2,q}(\Om)} + \Big(\frac{c}{K(q,\Om)}\Big)^{1/(\al-1)}\bigg]
\\
&\ \ \ \ \ \ \ \ + c  \| \mathrm{A}\: \D^2 \bar u \|_{\mathrm L^{\infty}(\Om)} + \| \bar u \|^\al_{\mathrm W^{1,\infty}(\Om)} +1 \bigg\}
\\
& =: \, \frac{c(1-\al)}{K(q,\Om)} \| w\|_{\mathrm W^{2,q}(\Om)} - M,
\end{split}
\]
from which \eqref{2.20A}(i) follows. Similarly, to obtain \eqref{2.20A}(ii), we have
\[
\begin{split}
\F_{p,\epsilon,\bar u}(w) \, &\geq \, \| F(\mathrm J^2w)\|_{\mathrm L^q(\Om)} 
\\
& \geq \, c \big\| \mathrm{A}\:\D^2 w -\bar \xi(\cdot,w,\D w)\big\|_{\mathrm L^q(\Om)} 
\\
& \geq \, c \| \mathrm{A}\:\D^2 w \|_{\mathrm L^q(\Om)}  - \big\| \big(|w|+|\D w|\big)^\al \big\|_{\mathrm W^{1,q}(\Om)} -1
\\
& \geq \, c \| \mathrm{A}\:\D^2 w \|_{\mathrm L^q(\Om)} - \|w\|^\al_{\mathrm W^{1,q}(\Om)}- 1.
\end{split}
\]
The lemma has been established.
\qed
\ms

Now we may continue with the proof the theorem.

\BPT \ref{theorem8}. (i) This is a simple application of the direct method. Fix an exponent $p>n+1/c^3$ and let $(w_j)_1^\infty \sub \E_{\bar u}^p(\Om) \sub W_{\bar u}^{2,p}(\Om)$ be a minimising sequence of $\F_{p,\epsilon,\bar u}$ in $\E_{\bar u}^p(\Om)$. By Lemma \ref{lemma13}, $(w_j)_1^\infty$ is bounded in $\mathrm W_{\bar u}^{2,p}(\Om)$, hence by the sequential weak compactness of bounded sets in $\mathrm W^{2,p}(\Om)$, there exists $u_p \in \E_{\bar u}^p(\Om)$ and a subsequence of indices $({j_k})_1^\infty$ such that $w_{j_k} \weak u_p$ in $\mathrm W_{\bar u}^{2,p}(\Om)$,  as $k\to\infty$. By assumption \eqref{1.25}, $\E_{\bar u}^p(\Om)$ is weakly closed, therefore $u_p \in \E_{\bar u}^p(\Om)$. Since we have assumed $p>n+1/c^3$, by \eqref{1.16} and \eqref{2.9A} tested against rank-one matrices of the form $(0,0,0,\xi) \ot (0,0,0,\xi)$, it follows that $\p^2_{\xi \xi}\big(|F(x,\eta,\mathrm p,\cdot)|^p \big)\geq 0$, which implies that $|F(x,\eta,\mathrm p,\cdot)|^p$ is convex. By \eqref{1.7} and \eqref{1.16}, we have the lower bound
\[
F(x,\eta,\mathrm p,\xi) \geq c\big(\xi-\bar \xi(x,\eta,\mathrm p)\big)\geq \, c\hspace{1pt} \xi -  \big(|\eta|+|\mathrm p|\big)^\al-1,
\]
and hence the functional \eqref{1.19} is sequentially weakly lower semi-continuous (see e.g.\ \cite[Th.\ 3.23, p.\ 96]{D}) on the weakly closed set $\E_{\bar u}^p(\Om)$. Hence, $u_p$ minimises the functional $\F_{p,\epsilon,\bar u}$ over $\E_{\bar u}^p(\Om)$.

\ms 

\noi (ii) Let $\{u_p\}_{n< p <\infty} \sub\mathrm W^{2,n}_{\bar u}(\Om)$ be the family of minimisers from part (i). By applying Lemma \ref{lemma13}, for each fixed $q\in(1,\infty)$, $\{u_p\}_{n< p <\infty}$ is bounded in $\mathrm W^{2,q}(\Om)$, as a result of \eqref{2.20A}(i). By sequential weak compactness and a standard diagonal argument, there exists a sequence of indices $(p_j)_1^\infty$ and $u_\infty \in \mathbb\mathrm W^{2,\infty}(\Om)$ (recall \eqref{2.2}) such that $u_{p_j} \weak u_\infty$ in $\mathrm W^{2,q}(\Om)$, as $j\to\infty$. By Morrey's theorem and the Lipschitz regularity of $\p\Om$, we also have that $u_{p_j} \larrow u_\infty$ in $\mathrm{C}^{1,\ga}(\overline{\Om})$, as $j\to\infty$.  Again by Lemma \ref{lemma13}, for any $q \in(1,\infty)$ fixed and $p\geq q$, by  \eqref{2.20A}(ii) we have that
\[
\|\mathrm{A}\:\D^2 u_p \|_{\mathrm L^{q}(\Om)}    \, \leq \, \frac{1+\de + \F_\infty(\bar u)+\|u_p\|^\al_{\mathrm W^{1,\infty}(\Om)}}{c}.
\]
By the sequential weak lower semi-continuity of the $\mathrm L^q$ norm, and the strong convergence of $u_{p_j} \larrow u_\infty$ in $\mathrm C^{1,\ga}(\overline{\Om})$ as $j\to\infty$, we have
\[
\|\mathrm{A}\:\D^2 u_\infty \|_{\mathrm L^{q}(\Om)} \leq \liminf_{j\to\infty} \|\mathrm{A}\:\D^2 u_{p_j} \|_{\mathrm L^{q}(\Om)}    \, \leq \, \frac{1+\de + \F_\infty(\bar u)+\|u_\infty\|^\al_{\mathrm W^{1,\infty}(\Om)}}{c}.
\]
Then, by letting $q\to\infty$ we infer that $\mathrm{A}\:\D^2 u_\infty \in \mathrm L^{\infty}(\Om)$, therefore $u_\infty \in \W^{2,\infty}_{\bar u}(\Om)$. Additionally, by \eqref{1.25}, for any fixed $q \in(1,\infty)$, the family $(u_p)_{p\geq q}$ is contained in $\E_{\bar u}^q(\Om)$. Since $\E_{\bar u}^q(\Om)$ is weakly closed, by letting $p_j\to\infty$ we infer that $u_\infty \in \E_{\bar u}^q(\Om)$ for all $q \in(1,\infty)$. By invoking \eqref{1.25} once again, we deduce that $u_\infty \in \E_{\bar u}^\infty(\Om)$. Now we show that $u_\infty$ minimises $\F_\infty +\frac{\epsilon}{2}\| \cdot - \bar u\|^2_{\mathrm L^2(\Om)}$, and also show the convergence of the energies. For any $v\in \E_{\bar u}^\infty(\Om)$ fixed and $p\geq q$, we have
\[
\begin{split}
\big\| F(\mathrm J^2u_p) \big\|_{\mathrm L^q(\Om)} +\frac{\epsilon}{2} \| u_p - \bar u\|^2_{\mathrm L^2(\Om)} & \leq \big\| F(\mathrm J^2u_p) \big\|_{\mathrm L^p(\Om)} +\frac{\epsilon}{2}\| u_p - \bar u\|^2_{\mathrm L^2(\Om)}
\\
& = \F_{p,\epsilon,\bar u}(u_p)
\\
& \leq \F_{p,\epsilon,\bar u}(v)
\\
& \leq \F_\infty(v) +\frac{\epsilon}{2}\| v - \bar u\|^2_{\mathrm L^2(\Om)}.
\end{split}
\]
By the sequential weak lower semi-continuity of the functional, we deduce
\beq
\label{2.27}
\begin{split}
\big\| F(\mathrm J^2u_\infty) \big\|_{\mathrm L^q(\Om)} +\frac{\epsilon}{2} \| u_\infty - \bar u\|^2_{\mathrm L^2(\Om)} & \leq \liminf_{p_j\to\infty} \F_{p,\epsilon,\bar u}(u_p)
\\
& \leq \limsup_{p_j\to\infty} \F_{p,\epsilon,\bar u}(u_p)
\\
& \leq \F_\infty(v) +\frac{\epsilon}{2}\| v - \bar u\|^2_{\mathrm L^2(\Om)},
\end{split}
\eeq
for any $v\in \E_{\bar u}^\infty(\Om)$. By letting $q\to\infty$  in  inequality \eqref{2.27}, the choice $v:=u_\infty$ completes the proof of \eqref{1.26}.

\ms

\noi (iii) If we have $\epsilon>0$ and $\bar u$ minimises $\F_\infty$ over $\E_{\bar u}^\infty(\Om)$, choosing $v:=\bar u$ in inequality \eqref{2.27} and letting $q\to\infty$ yields $\bar u = u_\infty$.
\qed
\ms


\subsection{A lemma of strong compactness in $\mathrm L^1$}\label{subsection2.6}

The following simple result, which is an extension of a lemma proved in the appendix of \cite{KP}, will be used in an essential fashion in Section \ref{section3}. (Let us recall here our notation $\tau_z$ for the translation operator.)

\begin{lemma} \label{lemma15} Let $\mathrm A \in \R^{n\by n}$ be a (strictly) positive matrix, $\Om \Subset \R^n$, $\sigma \in \mathrm L^1(\Om)$ and $\Sigma\in \mathrm L^1(\Om;\R^n)$. Suppose that $f \in \mathrm L^1(\Om)$ solves the PDE
\[
\ \ \ \ \ \ \ \div\big( \mathrm A \D f - \Sigma \big) + \sigma=0 \ \ \text{ in }\big(\mathrm C^\infty_c(\Om)\big)^{\!*},
\]
namely in the sense of distributions. Then, for any $\Om'\Subset \Om$, there exists a modulus of continuity $\om \in \mathrm{C}_{0^+}^\nnearrow[0,\infty)$ (depending only on $n,\Om',\Om,\mathrm A$), such that, 
\[
\|\tau_z f - f \|_{\mathrm L^1(\Om')}\, \leq\, \om(|z|)\Big(\|f\|_{\mathrm L^1(\Om)} + \|\sigma\|_{\mathrm L^1(\Om)}+ \|\Sigma\|_{\mathrm L^1(\Om)}\Big),
\]
for all $z\in \R^n$ with $|z|<\dist(\Om', \p\Om)$. Further, for any $q\in(1,n')$, where $n'=\frac{n}{n-1}$, there exists a constant $K=K(q,n,\Om',\Om,A)>0$ such that
\[
\| f\|_{\mathrm L^q(\Om')} \leq K\Big(\|f\|_{\mathrm L^1(\Om)} + \|\sigma\|_{\mathrm L^1(\Om)}+ \|\Sigma\|_{\mathrm L^1(\Om)}\Big).
\]
As a consequence, any set of distributional solutions which is bounded in $\mathrm  L^1(\Om)$, is also strongly precompact in $\mathrm L^q_{\mathrm {loc}}(\Om)$, for all $q\in(1,n')$.
\end{lemma}

\BPL \ref{lemma15}. Without loss of generality, by a change of variables argument as in Lemma \ref{lemma1}, we may assume that $\mathrm A$ is the identity matrix, and therefore we have the Poisson equation $\De f =\div(\Sigma)-\sigma$ in $\Om$. Fix $\Om'\Subset \Om$, and let fix also $\Om'',\Om'''$ such that $\Om'\Subset \Om''\Subset \Om'''\Subset \Om$. By the interior nature of the desired estimate and the linearity of the PDE, without loss of generality, a mollification argument allows us to assume that $u, \sigma \in \mathrm C^\infty(\overline{\Om'''})$, $\Sigma\in \mathrm C^\infty(\overline{\Om'''};\R^n)$ $\p\Om'''$ is piecewise $\mathrm C^\infty$ and the PDE holds classically in $\Om'''$. Let $\zeta \in \mathrm C^\infty_c(\R^n)$ be a cut off function which satisfies that $\chi_{\Om''} \leq \zeta \leq \chi_{\Om'''}$. By setting $\bar \sigma:= \zeta \sigma$ and $\bar \Sigma:= \zeta \Sigma$ the PDE implies that $\De f =\div(\bar \Sigma)-\bar \sigma$ in $\Om''$, and we have that $ \bar \sigma \in \mathrm C^\infty_c(\Om)$, $ \bar \Sigma \in \mathrm C^\infty_c(\Om;\R^n)$, and also $\bar \sigma |_{\Om''} \equiv \sigma$ and $\bar \Sigma |_{\Om''} \equiv \Sigma$. Let now $R=R_\Om>0$ be the radius of the smallest open ball such that $\Om - \Om +\mB_1 \sub \mB_R$. This means that $x-y+z \in \mB_R$, for all $x,y \in \Om$ and all $z\in \mB_1$. Let $\eta \in \mathrm C^\infty_c(\R^n)$ be a cut off function such that $\eta \equiv 1$ on $\mB_R$. By Green's formula (see e.g.\ \cite[Ch.\ 2]{GT}), we can represent the solution $f$ on $\Om''$ as
\[
f = h - \Phi *\bar \sigma + \Phi* \div (\bar \Sigma),
\]
for some harmonic function $h$ on $\Om''$, where $\Phi $ is the fundamental solution of the Laplace equation. Then, for any fixed $z\in\R^n$ with $|z|<\dist(\Om',\p\Om'')$, we have the estimate
\[
\begin{split}
\big\|\tau_z f -f \big\|_{\mathrm L^1(\Om')} & \leq \big\| \tau_z (\Phi *\div(\bar \Sigma)) -\Phi *\div(\bar \Sigma) \big\|_{\mathrm L^1(\Om')} 
\\
&\ \ \ \, + \big\| \tau_z (\Phi *\bar \sigma) -\Phi *\bar \sigma \big\|_{\mathrm L^1(\Om')} +  \|\tau_z h -h  \|_{\mathrm L^1(\Om')}.
\end{split}
\]
We will obtain the desired estimate by estimating each of the three terms in the above inequality separately. By the identities
\[
\tau_z(\Phi *(\div (\bar \Sigma)) = \tau_z(\D \Phi *\bar \Sigma) = (\tau_z\D \Phi) *\bar \Sigma,
\]
(where we interpret the convoluted inner product of vector functions in the above identities as $\D \Phi *\bar \Sigma = \sum_{i=1}^n \D_i \Phi *\bar \Sigma_i$), we estimate by using Young's inequality
\[
\begin{split}
\Big\|\tau_z(\Phi *(\div (\bar \Sigma)) -\Phi *(\div (\bar \Sigma)) \Big\|_{\mathrm L^1(\Om')}\, & =\, 
\Big\|\big(\tau_z(\D\Phi) -\D \Phi\big) *\bar \Sigma \Big\|_{\mathrm L^1(\Om')}
\\ 
& =\, 
\Big\|\big[\big(\tau_z(\eta \D\Phi) - \eta \D \Phi\big) \chi_{\mB_R}\big]* \bar \Sigma \Big\|_{\mathrm L^1(\R^n)}
\\
& \leq \, 
\big\|\tau_z(\eta\D\Phi) -\eta\D \Phi\big\|_{\mathrm L^1(\mB_R)} \| \bar \Sigma \|_{\mathrm L^1(\R^n)}
\\
& \leq \, 
\big\|\tau_z(\eta\D\Phi) -\eta\D \Phi\big\|_{\mathrm L^1(\R^n)} \| \Sigma \|_{\mathrm L^1(\Om)}
\end{split}
\]
Similarly, we have the estimate
\[
\begin{split}
\big\|\tau_z(\Phi * \bar \sigma) -\Phi *\bar \sigma \big\|_{\mathrm L^1(\Om')}\, \leq \, 
\big\|\tau_z(\eta\Phi) -\eta\Phi\big\|_{\mathrm L^1(\R^n)} \| \sigma \|_{\mathrm L^1(\Om)}.
\end{split}
\]
Further, by the mean value theorem for harmonic functions and interior derivative estimates, we have
\[
\|\tau_z h -h  \|_{\mathrm L^1(\Om')} \leq |z|\|\D h\|_{\mathrm L^\infty(\Om'')} \leq |z| \|h\|_{\mathrm L^1(\Om)}.
\]
Further, by using that $h=f+\Phi * \bar \sigma-\D \Phi * \bar \Sigma$, Young's inequality implies
\[
\begin{split}
\|h\|_{\mathrm L^1(\Om)}\, & \leq  \|f\|_{\mathrm L^1(\Om)} +  \big\|\Phi * \bar \sigma  \big\|_{\mathrm L^1(\Om)}  +  \big\|\D \Phi * \bar \Sigma  \big\|_{\mathrm L^1(\Om)}
\\
& = \|f\|_{\mathrm L^1(\Om)} +  \big\|(\Phi \chi_{\mB_R})* \bar \sigma  \big\|_{\mathrm L^1(\Om)}  +  \big\| (\D \Phi\chi_{\mB_R}) * \bar \Sigma  \big\|_{\mathrm L^1(\Om)} 
\\
&\leq \|f\|_{\mathrm L^1(\Om)} +  \|\Phi\|_{\mathrm L^1(\mB_R)} \| \sigma \|_{\mathrm L^1(\Om)}  +   \|\D \Phi\|_{\mathrm L^1(\mB_R)} \| \Sigma \|_{\mathrm L^1(\Om)} .
\end{split}
\]
Therefore, we have
\[
\|\tau_z h -h  \|_{\mathrm L^1(\Om')} \leq |z| \big(1 +  \|\Phi\|_{\mathrm W^{1,1}(\mB_R)} \big)\Big(\|f\|_{\mathrm L^1(\Om)} + \|\sigma\|_{\mathrm L^1(\Om)}+ \|\Sigma\|_{\mathrm L^1(\Om)}\Big).
\]
Finally, we conclude by putting the pieces together. Since $\eta\Phi \in \mathrm L^1(\R^n)$ and also $\eta\D\Phi \in \mathrm L^1(\R^n;\R^n)$, by the continuity of the translation operation, there exists $\om \in C^\nnearrow_0[0,\infty)$ such that
\[
\big\|\tau_z(\eta\Phi) -\eta\Phi\big\|_{\mathrm L^1(\R^n)}+ \big\|\tau_z(\eta\D\Phi) -\eta\D \Phi\big\|_{\mathrm L^1(\R^n)} \leq \om(|z|).
\]
By the above estimates, the desired $\mathrm L^1$ equi-continuity estimate ensues. Now we establish the $\mathrm L^q_{\mathrm{loc}}$ bound. Fix $q\in(1,n')$, and fixed also $\zeta \in \mathrm C^\infty_c(\Om)$ with $ \chi_{\Om'} \leq \zeta \leq \chi_{\Om'''}$. Let $\phi \in \mathrm C^\infty(\Om)$ be arbitrary, and test the PDE against $\ze \phi  \in \mathrm C^\infty_c(\Om)$, to obtain
\[
\begin{split}
\int_\Om \ze f (\mathrm A \: \D^2\phi) \,\mathrm d \mL^n &= -\int_{\Om'''} f \mathrm A \: \Big(\D  \ze \ot \D \phi + \D \phi  \ot \D \ze +\phi \D^2\ze \Big)\,\mathrm d \mL^n 
\\
&\ \ \ - \int_{\Om'''} \Sigma \cdot \Big(\phi\D \ze +\ze\D\phi \Big)\,\mathrm d \mL^n - \int_{\Om'''}  \sigma \ze \phi \,\mathrm d \mL^n.
\end{split}
\]
This identity implies the estimate
\[
\begin{split}
\left| \int_\Om \ze f (\mathrm A \: \D^2\phi) \,\mathrm d \mL^n \right| &\leq C(n)\big(|\mathrm A|+1\big)\|\ze\|_{\mathrm W^{2,\infty}(\Om)}
\\
&\ \ \ \cdot \Big(\|f\|_{\mathrm L^1(\Om)} + \|\sigma\|_{\mathrm L^1(\Om)}+ \|\Sigma\|_{\mathrm L^1(\Om)}\Big)\|\phi\|_{\mathrm W^{1,\infty}({\Om'''})},
\end{split}
\]
for some universal $C=C(n)>0$. Fix now any $\psi \in \mathrm C^\infty_c(\Om)$. By standard $\mathrm L^p$ theory for elliptic equations and bootstrap regularity (see e.g.\ \cite[Ch. 6 \& 9]{GT}, there exists a unique $\phi \in  \mathrm C^\infty(\Om)\cap  \mathrm C^1_0(\Om)$ such that $\mathrm A \: \D^2 \phi =\psi$ on $\Om$, and also
\[
\|\phi\|_{\mathrm W^{1,\infty}({\Om'''})} \leq C \|\mathrm A \: \D^2 \phi\|_{\mathrm L^{q'}({\Om'''})},
\]
for some $C=C(q,n,\Om)>0$, which is possible by Morrey's embedding theorem because $q'>n$ (since by assumption $q\in (1,n')$), and also $\Om''' \Subset \Om$ has piecewise smooth boundary. Therefore, there exists a $K=K(q,n,\Om',\Om,\mathrm A)>0$ such that
\[
\begin{split}
\left| \int_\Om (\ze f )\psi \,\mathrm d \mL^n \right| &\leq K \Big(\|f\|_{\mathrm L^1(\Om)} + \|\sigma\|_{\mathrm L^1(\Om)}+ \|\Sigma\|_{\mathrm L^1(\Om)}\Big) \|\psi\|_{\mathrm L^{q'}(\Om)},
\end{split}
\]
for any $\psi \in \mathrm C^\infty_c(\Om)$. By standard $\mathrm L^p$ theory (see e.g.\ \cite[Ch.\ 14]{KV}), this implies
\[
\begin{split}
\|f\|_{\mathrm L^{q}(\Om')} \leq \|\zeta f\|_{\mathrm L^{q}(\Om)}\leq K\Big(\|f\|_{\mathrm L^1(\Om)} + \|\sigma\|_{\mathrm L^1(\Om)}+ \|\Sigma\|_{\mathrm L^1(\Om)}\Big),
\end{split}
\]
as claimed. The last claim of the statement regarding strong compactness follows by the F\'echet-Kolmogorov and the Vitali convergence theorems. \qed
\ms


\subsection{Nodal sets and non-trivial solutions}\label{subsection2.7}

The following result lists some properties of solutions to equations of the form \eqref{1.9}, but with general $L^\infty$ coefficients. It holds true much more generally, but we specialise the statement only to the case of interest in this paper.

\begin{proposition} \label{proposition19} Let $\mathrm A \in \R^{n\by n}$ be a (strictly) positive matrix, $\Om \Subset \R^n$, $\mathrm K \in \mathrm L^\infty(\Om)$ and $\mathrm L \in \mathrm L^\infty(\Om;\R^n)$. Suppose that $f \in \mathrm L^1(\Om)$ solves the PDE
\beq
\label{PDE}
\ \ \ \ \ \ \div\big( \mathrm A \D f - \mathrm L f \big) + \mathrm K f =0 \ \, \text{ in }\Om,
\eeq
in the sense of distributions, namely in $\big(\mathrm C^\infty_c(\Om)\big)^{\!*}$. Then: 

\smallskip

\noi \emph{(i)} Any solution $f \in \mathrm L^1(\Om)$ to \eqref{PDE} satisfies $f \in {\mathrm W}^{1,p}_{\mathrm{loc}}(\Om)$, for all $p\in(1,\infty)$, and thus (the precise representative of) $f$ is in $\mathrm C^{0,\al}(\Om)$, for all $\al\in (0,1)$. 

\smallskip

\noi \emph{(ii)} There exist a non-trivial solution $f^*\not \equiv 0$ to \eqref{PDE}.

\smallskip

\noi \emph{(iii)} For any non-constant solution $f$ to \eqref{PDE}, the nodal set $\{f = 0\}$ is closed and countably $(n-1)$-rectifiable. Thus, the $\mH^{n-1}$-Hausdorff measure of $\{f = 0\}$ is locally finite, and in particular it is an $\mL^n$-nullset in $\Om$.
\end{proposition}

\BPP \ref{proposition19}. (i) Since the  coefficients $ \mathrm K, \mathrm L $ of \eqref{PDE} are essentially bounded and $\mathrm A$ is a constant positive matrix, by the result \cite[Corollary 1.2.8, p.\ 10]{BKR} of Bogachev-Krylov-R\"ockner, it follows that $f\in \mathrm W^{1,p}_{\mathrm{loc}}(\Om)$ for all $p\in(1,\infty)$. Therefore, upon identifying $f$ with its precise representative, we have $f\in \mathrm C^{0,\ga}(\Om)$ for all $\ga\in (0,1)$, as a consequence of Morrey's imbedding theorem. In particular, $\{f =0\}$ is a relatively closed set in $\Om$.

\smallskip

\noi (ii) For any $\mu \in \R$, we define
\[
\mB_\mu[f,g]:=\int_\Om\Big(\big( \mathrm A \D f - \mathrm L f \big)\cdot \D g +(\mu - \mathrm K) fg \Big)\, \mathrm d \mL^n.
\]
Then, $\mB_\mu : \mathrm W^{1,2}_0(\Om) \by \mathrm W^{1,2}_0(\Om) \larrow \R$ is a continuous bilinear form, which is also coercive for all $\mu\geq \mu^*$, where the value $\mu^*>0$ depends only on the coefficients. By the Lax-Milgram theorem, the operator $ \mathfrak{L}_\mu : \mathrm W^{1,2}_0(\Om) \larrow  \mathrm W^{-1,2}(\Om)$ given by
\[
\mathfrak{L}_\mu f := \div\big( \mathrm A \D f - \mathrm L f \big) + \mathrm K f - \mu f,
\]
is a linear isomorphism. Therefore, the inverse $ (\mathfrak{L}_\mu)^{-1} :  \mathrm W^{-1,2}(\Om) \larrow \mathrm  W^{-1,2}(\Om)$ is a compact operator, as a consequence of the compactness of the canonical inclusions $\mathrm W^{1,2}_0(\Om) \Subset  \mathrm L^{2}(\Om)\Subset  \mathrm W^{-1,2}(\Om)$. By the Fredholm alternative, the equation 
\beq
\label{inteq}
f+(\mu (\mathfrak{L}_\mu)^{-1})f = (\mathfrak{L}_\mu)^{-1}\phi
\eeq
either it has a solution $f\not\equiv0$ when $\phi \equiv 0$, or else it is solvable for all $\phi \in \mathrm  W^{-1,2}(\Om)$. By noting that \eqref{inteq} can be rewritten as $\mathfrak{L}_\mu f = \phi -\mu f$, which is equivalent to $\mathfrak{L}_0 f=\phi$, it follows that either $\mathfrak{L}_0 f=0$ has a solution $f^*\in \mathrm W^{1,2}_0(\Om)\set \{0\}$, or else the problem $\mathfrak{L}_0 h=\mathrm K - \div \,\mathrm L$ has a solution $h\in \mathrm W^{1,2}_0(\Om)$. The latter is equivalent to the solvability of the boundary value problem $\mathfrak{L}_0 f=0$ in $\Om$ with $f=1$ on $\p\Om$ by some $f^*\in \mathrm W^{1,2}(\Om)\set\{0\}$. In either case, it follows that \eqref{PDE} has a non-trivial solution $f^*\not\equiv 0$.

\smallskip

\noi (iii) We begin by considering the auxiliary boundary value problem
\[
\div(\mathrm A \D \psi) = \div \, \mathrm L \, \text{ on }\Om,\ \ \ \psi =0 \, \text{ on }\p\Om.
\]
By the Lax-Milgram theorem, there exists a unique solution $\psi \in \mathrm W^{1,2}_0(\Om)$, which by a result of Di Fazio \cite[Th.\ 3.1]{DiF}, it satisfies $\psi \in \mathrm W^{1,p}_{\mathrm{loc}}(\Om)$ for all $p\in(1,\infty)$. By Morrey's theorem, (the precise representative of) $\psi$ is in $\mathrm C^{0,\al}(\Om)$, for all $\al\in(0,1)$. We now set $\phi:= e^{-\psi}f$, where $f \in \mathrm L^1(\Om)$ is any given non-constant solution to \eqref{PDE}. By part (i), all three function $\phi,\psi, f$ lie in $\mathrm W^{1,p}_{\mathrm{loc}}(\Om)$ for all $p\in(1,\infty)$, whilst $\mathrm L \in \mathrm L^\infty(\Om;\R^n)$. Therefore, by using that $\div  ( \mathrm A\D \psi - \mathrm L )=0$ in $\Om$, we may rigorously compute
\[
\begin{split}
0 &= \div\big( \mathrm A \D f - \mathrm L f \big) + \mathrm K f
\\
&= \div\Big( e^{\psi}\mathrm A\D \phi + e^{\psi}\phi \big(\mathrm A\D \psi -\mathrm L\big)\Big) + \mathrm K e^{\psi}\phi
\\
&= \div\big( (e^{\psi}\mathrm A)\D \phi\big) + \big(\mathrm A\D \psi -\mathrm L\big) \cdot [e^{\psi}\D\phi+e^{\psi}\phi\D\psi] + (e^{\psi}\mathrm K)\phi ,
\end{split}
\]
in $\Om$. Hence, the function $\phi$ is a weak solution to the divergence PDE
\[
\div\big( (e^{\psi}\mathrm A)\D \phi\big) + \big[e^{\psi}\big(\mathrm A\D \psi -\mathrm L\big)\big] \cdot \D\phi + \big[e^{\psi}\big(\!\big(\mathrm A\D \psi -\mathrm L\big) \cdot \D\psi + \mathrm K\big) \big]\phi= 0 \ \text{ in }\Om.
\]
By the previous, the coefficients satisfy
\[
e^{\psi}\big(\mathrm A\D \psi -\mathrm L\big) , \ e^{\psi}\big(\!\big(\mathrm A\D \psi -\mathrm L\big) \cdot \D\psi + \mathrm K\big)\, \in \mathrm L^{p}_{\mathrm{loc}}(\Om), \ \ e^{\psi}\mathrm A \in \mathrm C^{0,\al}(\Om;\R^{n\by n}),
\]
for all $p\in(1,\infty)$, and all $\al\in(0,1)$. By the result \cite[Th.\ 1.7 \& Rem.\ 1.8(3)]{HS} of Hardt and Simon, it follows that for any point $\bar x \in \Om$ at which $\phi$ has finite order of vanishing, the set $\{\phi = 0\}\cap \mB_\rho(\bar x)$ is countably $(n-1)$-rectifiable, for all sufficiently small $\rho>0$. Since $f= e^{\psi}\phi$, we have $\{f = 0\} = \{\phi = 0\}$, and since $e^{\pm\psi} \in \mathrm L^\infty_{\mathrm{loc}}(\Om)$, it follows that $f$ and $\phi$ have the same sets of points with finite order of vanishing. On the other hand, by the strong unique continuation principle applied to \eqref{PDE}, it follows that $f$ has finite order of vanishing everywhere on $\Om$ (see e.g.\ Koch-Tataru \cite{KT2, KT}). Hence,  $\{f = 0\}$ is countably $(n-1)$-rectifiable, and in particular it is a Lebesgue nullset. The proposition ensues.  \qed
\ms


\section{Proofs of the main results}
\label{section3}

We may now establish our first main result.

\BPT \ref{theorem1}. \noi {\bf Step 1.} Let $u_\infty \in \W^{2,\infty}(\Om)$ and $f_\infty \in \mathrm L^1(\Om)$ be given, and suppose they satisfy the system of equations \eqref{1.8}-\eqref{1.9}. We assume first that 
\[
f_\infty \not\equiv 0.
\]
\noi {\bf Step 2.} We define the functions
\beq
\label{2.11}
\mathrm K_\infty :=\frac{\p_\eta F(\mathrm J^2u_\infty)}{\p_\xi F(\mathrm J^2u_\infty)},\ \ \ \ \mathrm L_\infty :=\frac{\p_{\mathrm p} F(\mathrm J^2u_\infty)}{\p_\xi F(\mathrm J^2u_\infty)},
\eeq
and note that $\mathrm K_\infty,|\mathrm L_\infty|$ are in $\mathrm L^{\infty}(\Om)$: indeed, by \eqref{1.7} we have that $\p_\xi F(\mathrm J^2u_\infty) \geq c$ a.e.\ on $\Om$, and also
\[
\big| \p_\eta F(\mathrm J^2u_\infty) \big|+ \big| \p_{\mathrm p} F(\mathrm J^2u_\infty) \big| \leq C\Big( \|u_\infty\|_{\mathrm W^{1,\infty}(\Om)} + \|\mathrm{A}\:\D^2 u_\infty\|_{\mathrm L^{\infty}(\Om)} \Big),
\]
a.e.\ on $\Om$. 

\smallskip

\noi {\bf Step 3.} By \eqref{1.9} and \eqref{2.11}, $f_\infty \in \mathrm L^1(\Om)$ is a distributional solution to
\beq
\label{1.19A}
\ \ \ \ \div\big( \mathrm{A} \D f_\infty -\mathrm L_\infty f_\infty\big) + \mathrm K_\infty f_\infty \, =\, 0 \  \ \ \ \text{in }\Om. 
\eeq
By Proposition \ref{proposition19}, we have $\mL^n\big(\{f_\infty = 0\}\big)=0$ and $f_\infty  \in  {\mathrm W}^{1,p}_{\mathrm{loc}}(\Om)$ for all $p\in(1,\infty)$. Additionally, upon identifying $f_\infty$ with its precise representative, the nodal set $\{f_\infty = 0\}$ is countably $(n-1)$-rectifiable and closed.

\smallskip

\noi {\bf Step 4.} We will now establish that there exists some $C \in \mathrm{C}^\nnearrow[0,\infty)$ such that
\[
 \big|F(x,\eta,\mathrm p, \xi) \big| \,\leq\, C\big(|\eta|+|\mathrm p|+|\xi| \big),
\]
for a.e.\ $x\in\Om$ and all $(\eta,\mathrm p,\xi) \in \R \by\R^n \by \R$. To this aim, recalling \eqref{xibar} and assumption \eqref{1.7}, it follows that
\[
\begin{split}
 \big|F(x,\eta,\mathrm p,\xi) \big| \, & =\,  \Big|F(x,\eta,\mathrm p,\xi) -  F(x,\eta,\mathrm p,\bar \xi(x,\eta,\mathrm p)) \Big|
\\
& =\,  \left| \int_{\bar \xi(x,\eta,\mathrm p)}^\xi   
\p_\xi F(x,\eta,\mathrm p,\cdot) \, \mathrm d \mL^1 \right|,
\end{split}
\]
which yields
\[
\begin{split}
 \big|F(x,\eta,\mathrm p,\xi) \big| \, &  \leq \, \big| \bar \xi(x,\eta,\mathrm p) - \xi \big|  C\Big(| \bar \xi(x,\eta,\mathrm p)| +|\xi|+ |\eta| +|\mathrm p| \Big)
\\
& \leq \, \big( C\big(|\eta|+|\mathrm p|\big)  +|\xi |\big)  C\Big(C\big(|\eta|+|\mathrm p|\big) +|\xi|+ |\eta|+|\mathrm p| \Big).
\end{split}
\]
Hence, the desired estimate ensues for the new function in $\mathrm{C}^\nnearrow[0,\infty)$ given by $t\mapsto \big[C(t)+t\big]C\big(C(t)+t\big)$.

\smallskip

\noi {\bf Step 5.} By assumption \eqref{1.7}, we may fix $x\in \Om$ (this is possible on a subset of full measure) and apply Taylor's theorem to $|F(x,\cdot,\cdot,\cdot)|^2 \in \mathrm{C}^2 (\R\by \R^n\by \R)$, to deduce that
\beq
\label{2.14}
\begin{split}
\big|F(\mathrm J^2u_\infty &+\mathrm J^2\psi) \big|^2 = \big|F(\mathrm J^2u_\infty) \big|^2 \, +\, 2 F(\mathrm J^2u_\infty) \Big\{ \p_* F(\mathrm J^2u_\infty) \cdot \mathrm J^2_*\psi\Big\} 
\\
& +\bigg( \int_0^1 (1-\la ) \p_{**}^2 (|F|^2) \big(\mathrm J^2u_\infty+ \la \mathrm J^2\psi \big) \, \mathrm d \la \bigg) : \mathrm J^2_*\psi \ot \mathrm J^2_*\psi,
\end{split}
\eeq
a.e.\ on $\Om$. By the identity \eqref{2.9A} for $p=2$, assumption \eqref{1.7} and Step 4, we have the bound
\[
\big| \p_{**}^2 (|F|^2)(x,\eta,\xi) \big| \,\leq\, C\big(|\eta|+|\mathrm p|+|\xi| \big)
\]
for some function $C \in \mathrm{C}^\nnearrow[0,\infty)$. By using this bound and equation \eqref{1.8} (which implies $\big|F(\mathrm J^2u_\infty) \big| = \F_\infty(u_\infty)$ a.e.\ on $\Om$), \eqref{2.14} yields
\beq
\label{2.15}
\begin{split}
\big|F(\mathrm J^2u_\infty +\mathrm J^2\psi) \big|^2  & \geq \,  |\F_\infty(u_\infty)|  ^2 
\\
&+\, 2 \F_\infty(u_\infty) \sgn(f_\infty)\Big\{ \p_* F(\mathrm J^2u_\infty) \cdot \mathrm J^2_*\psi\Big\} 
\\
&  - \,C\Big(\| \mathrm J^2_* u_\infty\|_{\mathrm L^{\infty}(\Om)}+  \| \mathrm J^2_*\psi \|_{\mathrm L^{\infty}(\Om)} \Big)   \| \mathrm J^2_*\psi \|^2_{\mathrm L^{\infty}(\Om)},
\end{split}
\eeq
a.e.\ on $\Om$. 

\smallskip

\noi {\bf Step 6.} Let us now define the functional $\Theta_\infty : \W^{2,\infty}_0(\Om) \larrow \R$ by setting
\beq
\label{2.16}
\Theta_\infty(\psi) \, :=\, \underset{\Om}{\ess\sup} \Big[  \sgn(f_\infty)\Big\{ \p_* F(\mathrm J^2u_\infty) \cdot \mathrm J^2_*\psi\Big\}\Big],
\eeq
and set 
\beq
\label{2.17}
N_\infty \,:=\, C\big( 1 +  \| \mathrm J^2_* u_\infty \|_{\mathrm L^{\infty}(\Om)} \big). 
\eeq
In view of \eqref{2.16}-\eqref{2.17}, we may deduce from \eqref{2.15} that
\beq
\label{2.18}
\begin{split}
(\F_\infty(u_\infty+\psi))^2  & \geq \,  (\F_\infty(u_\infty))^2  +  \Big\{2 \F_\infty(u_\infty) \Theta_\infty(\psi) - N_\infty  \| \mathrm J^2_*\psi \|^2_{\mathrm L^{\infty}(\Om)}\Big\},
\end{split}
\eeq
for all $\psi \in \W^{2,\infty}_0(\Om)$ satisfying $ \| \mathrm J^2_*\psi \|_{\mathrm L^{\infty}(\Om)} \leq 1$ (namely in the closed unit ball of $\W^{2,\infty}_0(\Om)$ with respect to the $\| \cdot \|'_{\W^{2,\infty}_0(\Om)}$-norm in \eqref{norms}).

\ms

\noi {\bf Step 7.}  Let $\K \sub \mS^{\W^{2,\infty}_0(\Om)}$ be a compact set. We will regard the unit sphere of $\W^{2,\infty}_0(\Om)$ to be with respect to the $\smash{\| \cdot \|'_{\W^{2,\infty}_0(\Om)}}$-norm in  \eqref{norms}. We aim at showing that there exists $\theta_\infty>0$ such that
\beq
\label{2.19}
\ \ \ \Theta_\infty(\psi)  \geq \, \theta_\infty, \ \ \ \forall\, \psi \in \K.
\eeq
Temporarily assuming that \eqref{2.19} holds true, we will now show how we can complete the proof of the theorem. From \eqref{2.18}-\eqref{2.19}, for any $\psi \in \K$ and $0<t <1$, by noting that $\Theta_\infty$ is $1$-homogeneous, we have 
\[
(\F_\infty(u_\infty+ t\psi))^2  \geq \,  (\F_\infty(u_\infty))^2  + \, t\Big\{2 \F_\infty(u_\infty) \theta_\infty - t N_\infty \Big\} ,
\]
for all $\psi \in \K$. By restricting $t\in (0,\rho)$, where
\[
\rho := \frac{2 \F_\infty(u_\infty) \theta_\infty}{N_\infty},
\]
the desired conclusion \eqref{1.10} ensues (having assumed that \eqref{2.19} holds true).

\ms

\noi {\bf Step 8.} The remaining of the proof is devoted to establishing the lower bound \eqref{2.19} for the functional $\Theta_\infty : \W^{2,\infty}_0(\Om) \larrow \R$ defined in \eqref{2.16}. We will first show
\beq
\label{2.20}
\Theta_\infty(\psi) > 0, \ \ \  \psi \in \W^{2,\infty}_0(\Om)  \set \{0\}.
\eeq
To this aim, fix $\psi \in \W^{2,\infty}_0(\Om)  \set \{0\}$. By Lemma \ref{lemma2}(ii), the divergence equation \eqref{1.9} is satisfies also weakly in $(\W^{2,\infty}_0(\Om))^*$, namely when testing against functions in the space $\W^{2,\infty}_0(\Om)$. Thus, \eqref{1.9} implies
\beq
\label{2.21}
\int_{\Om} f_\infty \Big( \mathrm{A}\:\D^2\psi + \mathrm L_\infty \! \cdot \D \psi+ \mathrm K_\infty \psi \Big) \, \d\mL^n \, =\, 0.
\eeq
We now claim that 
\beq
\label{2.22}
\mL^n\Big( \Big\{ f_\infty \Big( \mathrm{A}\:\D^2\psi + \mathrm L_\infty \! \cdot \D \psi+ \mathrm K_\infty \psi \Big) \neq 0 \Big\} \Big)>0. 
\eeq
If for the sake of contradiction we suppose that 
\[
\mathrm{A}\:\D^2\psi + \mathrm L_\infty \!\cdot \D \psi+ \mathrm K_\infty \psi =0, \ \text{ a.e.\ on }\Om.
\]
Since $|\psi|= |\D \psi| = 0$ on $\p\Om$ and $\psi \in W_0^{2,q}(\Om)$ for all $q\in(1,\infty)$, the extension $\bar \psi $ of $\psi$ by zero to $\R^n\set\Om$ is compactly supported in $\R^n$ and satisfies that $\bar \psi \in\mathrm W^{2,q}_c(\R^n)$, solving strongly the PDE
\[
\mathrm{A}\:\D^2 \bar \psi + \bar {\mathrm L}_\infty \! \cdot \D \bar \psi + \bar {\mathrm K}_\infty \bar \psi =0 \ \ \ \text{ a.e.\ on }\R^n,
\] 
where $\bar {\mathrm K}_\infty , |\bar {\mathrm L}_\infty | \in \mathrm L^{\infty}_c(\R^n)$ are the extensions of ${\mathrm K}_\infty, {\mathrm L}_\infty $ by zero on $\R^n \set\Om$. In that case, by \cite[Th.\;(1.7) and Rem.\;(1.8)]{HS} we infer that $\mL^n(\{\bar \psi=0\})=0$, contradicting that we have selected $\bar \psi \equiv 0$ on $\R^n\set\{0\}$. As a conclusion, \eqref{2.22} is true. Since we already know from Step 3 that $f_\infty=0$ at most on an $\mL^n$-nullset, \eqref{2.22}  follows. Further, we already know from \eqref{2.21} that the average must vanish over $\Om$. As a result, there must exist Borel sets $\Om^\pm \sub \Om$ with $\mL^n(\Om^\pm)>0$ such that 
\[
\frac{f_\infty}{|f_\infty|}\Big(\mathrm{A}\:\D^2\psi + \mathrm L_\infty \! \cdot \D \psi + \mathrm K_\infty \psi\Big) \gtrless \, 0 \ \ \ \text{ a.e.\ on }\Om^\pm.
\]
In view of \eqref{2.11}, the above yields
\[
\ \ \ \ \ \frac{f_\infty}{|f_\infty|}\Big( \p_\xi F(\mathrm J^2 u_\infty) \mathrm{A}\:\D^2\psi + \p_{\mathrm p} F(\mathrm J^2 u_\infty) \cdot \D \psi+  \p_\eta F(\mathrm J^2 u_\infty) \psi\Big) > \, 0, \ \ \ \text{ a.e.\ on }\Om^+.
\]
Using our notation for jets, we can rewrite the above as
\[
\sgn(f_\infty) \Big\{ \p_* F(\mathrm J^2 u_\infty) \cdot \mathrm J^2_*\psi\Big\} > \, 0 \ \ \ \text{ a.e.\ on }\Om^+,
\]
and $\mL^n(\Om^+)>0$. In conclusion, we obtain (in view of \eqref{2.16}) that
\[
\Theta_\infty(\psi) \, =\, \underset{\Om}{\ess\sup} \Big[ \sgn(f_\infty) \Big\{ \p_* F(\mathrm J^2 u_\infty) \cdot \mathrm J^2_*\psi\Big\}  \Big] >0,
\]
for any $\psi \in \W^{2,\infty}_0(\Om)  \set \{0\}$.

\ms

\noi {\bf Step 9.} Now we establish the bound \eqref{2.19} for the functional $\Theta_\infty$. In view of Step 8, it suffices to show that $\Theta_\infty$ is continuous. Indeed, on the one hand, for any $\phi, \psi \in \W^{2,\infty}_0(\Om)$ we have the estimate
\[
\Theta_\infty(\phi) =  \, \Theta_\infty(\phi +\psi -\psi) \, \leq \, \Theta_\infty(\phi+\psi) \, + \, \| \mathrm J^2_* u_\infty \|_{\mathrm L^{\infty}(\Om)} \| \mathrm J^2_*\psi \|_{\mathrm L^{\infty}(\Om)} .
\]
On the other hand, by the convexity of $\Theta_\infty$, we have that
\[
\begin{split}
\Theta_\infty(\phi+\psi) - \Theta_\infty(\phi) \, &\leq  \bigg( \Theta_\infty\bigg(\phi + \frac{\psi}{\| \mathrm J^2_*\psi \|_{\mathrm L^{\infty}(\Om)} }\bigg) - \Theta_\infty(\phi)\bigg) \| \mathrm J^2_*\psi \|_{\mathrm L^{\infty}(\Om)}
\\
&\leq \, 2  \| \mathrm J^2_* u_\infty \|_{\mathrm L^{\infty}(\Om)} \Big( 1+ \| \mathrm J^2_*\phi \|_{\mathrm L^{\infty}(\Om)} \Big) \| \mathrm J^2_*\psi \|_{\mathrm L^{\infty}(\Om)} .
\end{split}
\]
\noi {\bf Step 10.} In conclusion, $\Theta_\infty$ is continuous and positive on the compact set $\K$, yielding that
\[
\theta_\infty\,:=\, \inf_{\K} \Theta_\infty  \, > \, 0.
\]
The proof of Theorem \ref{theorem1} is now complete.
\qed
\ms

We may now establish our second main result.

\BPT \ref{theorem2}. Let $p\geq 2$ and fix $\psi \in \W^{2,\infty}_0(\Om)$. We continue our arguments from the proof of Step 5 of Theorem \ref{theorem1}. By Taylor's theorem for $|F(x,\cdot,\cdot)|^p$, we have
\beq
\label{2.23}
\begin{split}
\big|F(\mathrm J^2u_\infty + & \mathrm J^2\psi)  \big|^p  = \big|F(\mathrm J^2u_\infty) \big|^p
\\
 &+\, p \big| F(\mathrm J^2u_\infty) \big|^{p-2} F(\mathrm J^2u_\infty) \Big\{ \p_* F(\mathrm J^2u_\infty) \cdot \mathrm J^2_*\psi\Big\} 
\\
& +\bigg( \int_0^1 (1-\la ) \p_{**}^2 (|F|^p) \big(\mathrm J^2u_\infty+ \la \mathrm J^2\psi \big) \, \mathrm d \la \bigg) : \mathrm J^2_*\psi \ot \mathrm J^2_*\psi,
\end{split}
\eeq
a.e.\ on $\Om$. We will now show that, under \eqref{1.13}-\eqref{1.15}, the last term of \eqref{2.23} above is non-negative for all sufficiently large $p\geq 2$, when $\psi$ is restricted appropriately. More precisely: 

\ms

\noi \underline{Case 1}: Suppose that \eqref{1.13} is satisfied. Then, the essential range of $\mathrm J^2 u_\infty$ satisfies $\ess ((\mathrm J^2u_\infty)(\Om)) \sub \Om \by L \by R$. Since $L,R$ are open intervals and the essential range is a compact set, for any $\psi \in \W^{2,\infty}_0(\Om)$ and any $\la \in(0,1)$, we may select $\rho >0$ small enough such that $\|\mathrm J^2_* \psi\|_{\mathrm L^{\infty}(\Om)}<\rho$, in order to obtain
\[
\ess \big((\mathrm J^2u_\infty+ \la \mathrm J^2\psi)(\Om)\big) \sub \Om \by L \by R.
\]
By choosing the corresponding $\bar p \geq 2$, it follows that $\p_{**}^2 (|F|^p) \big(\mathrm J^2u_\infty+ \la \mathrm J^2\psi \big) \geq 0$ for all $p\geq \bar p$, when $\psi \in \mB_\rho^{\W^{2,\infty}_0(\Om)}$.

\smallskip

\noi \underline{Case 2}: Suppose \eqref{1.14} is satisfied. Then, the essential range satisfies $\ess (u_\infty(\Om)) \sub L $. Since $L$ is an open interval and the essential range is a compact set, for any $\psi \in \W^{2,\infty}_0(\Om)$ and any $\la \in(0,1)$, we may select $\rho >0$ small enough such that $\|\psi\|_{\mathrm L^{\infty}(\Om)}<\rho$, and any sufficiently large interval $R\sub \R$, in order to obtain
\[
\ess \big((\mathrm J^2u_\infty+ \la \mathrm J^2\psi)(\Om)\big) \sub \Om \by L \by R.
\]
Then, for the corresponding $\bar p \geq 2$, it follows that $\p_{**}^2 (|F|^p) \big(\mathrm J^2u_\infty+ \la \mathrm J^2\psi \big) \geq 0$ for all $p\geq \bar p$, when $\psi \in \mB_\rho^{\mathrm L^{\infty}(\Om)} \cap \W^{2,\infty}_0(\Om)$.

\smallskip

\noi \underline{Case 3}: Suppose \eqref{1.15} is satisfied. Fix any $\psi \in \W^{2,\infty}_0(\Om)$ and choose $R,L\sub \R$ large enough and the corresponding $\bar p \geq 2$, so that the essential range satisfies
\[
\ess \big(\big(\mathrm J^2u_\infty+ \la \mathrm J^2\psi\big)(\Om)\big) \sub \Om \by L \by R.
\]
It follows that $\p_{**}^2 (|F|^p) \big(\mathrm J^2u_\infty+ \la \mathrm J^2\psi \big) \geq 0$ for all $p\geq \bar p$.

\ms

In conclusion, for any fixed $\psi \in \W^{2,\infty}_0(\Om)$  appropriately restricted depending on the assumptions, we have
\[
\begin{split}
\big|F(\mathrm J^2u_\infty + & \mathrm J^2\psi)  \big|^p  - \big|F(\mathrm J^2u_\infty) \big|^p 
\\
 &\geq \, p \big| F(\mathrm J^2u_\infty) \big|^{p-2} F(\mathrm J^2u_\infty) \Big\{ \p_* F(\mathrm J^2u_\infty) \cdot \mathrm J^2_*\psi\Big\} ,
 \end{split}
\]
a.e.\ on $\Om$. In view of \eqref{1.8} and \eqref{1.1}, the above yields
\[
\begin{split}
\underset{\Om}{\ess\sup}\big|F(\mathrm J^2u_\infty + & \mathrm J^2\psi)  \big|^p  - (\F_\infty(u_\infty))^p 
\\
 &\geq \, p (\F_\infty(u_\infty))^{p-1}  \underset{\Om}{\ess\sup} \left[\sgn(f_\infty) \Big\{ \p_* F(\mathrm J^2u_\infty) \cdot \mathrm J^2_*\psi\Big\} \right].
\end{split}
\]
By \eqref{2.16} and \eqref{1.1}, since powers commute with the essential supremum of non-negative functions, the above inequality can be rewritten as
\[
\begin{split}
\big(F_\infty(u_\infty + \psi)  \big)^p  - (\F_\infty(u_\infty))^p 
\, \geq \, p (\F_\infty(u_\infty))^{p-1} \Theta_\infty(\psi).
\end{split}
\]
By Step 8 in the proof of Theorem \ref{theorem1}, we know that $\Theta_\infty(\psi)>0$ for all functions $\psi \in \W^{2,\infty}_0(\Om)\set\{0\}$. Therefore, if $\F_\infty(u_\infty)>0$, we have obtained the inequality
\[
 \F_\infty(u_\infty)  <  \F_\infty(u_\infty+\psi),
\]
when $\psi \neq 0$. If on the other hand $\F_\infty(u_\infty)=0$, the previous arguments do not allow to conclude that $\F_\infty(u_\infty+\psi)>0=\F_\infty(u_\infty)$, when $\psi \in \W^{2,\infty}_0(\Om)\set\{0\}$. To show that the last inequality is indeed strict, suppose for the sake of contradiction that instead $F(\mathrm J^2u_\infty+\mathrm J^2\psi)=0$ a.e.\ on $\Om$, for some $\psi \in \W^{2,\infty}_0(\Om)\set\{0\}$. In view of \eqref{xibar}, this can be restated as
\beq 
\label{E1}
\mathrm A\:\D^2 \psi  = \bar \xi \big(\cdot, u_\infty+\psi, \D u_\infty+ \D \psi \big) - \mathrm A\:\D^2 u_\infty, \ \ \ \text{a.e.\ in }\Om.
\eeq
Further, since we also have $F(\mathrm J^2u_\infty)=0$ a.e.\ on $\Om$, \eqref{E1} is true when $\psi \equiv 0$ as well, namely
\beq
\label{E2}
\mathrm A\:\D^2 u_\infty  = \bar \xi \big(\cdot, u_\infty,   \D u_\infty\big), \ \ \ \text{a.e.\ in }\Om.
\eeq
By \eqref{E1}-\eqref{E2}, we deduce
\beq 
\label{E3}
\mathrm A\:\D^2 \psi  =  \bar \xi \big(\cdot, u_\infty+\psi, \D u_\infty+ \D \psi \big) -  \bar \xi \big(\cdot, u_\infty,   \D u_\infty\big), \ \ \ \text{a.e.\ in }\Om.
\eeq
By Taylor's theorem, we can write
\beq 
\label{E4}
\begin{split}
  \bar \xi \big(\cdot, u_\infty+\psi, \D u_\infty+ \D \psi \big) &-  \bar \xi \big(\cdot, u_\infty,   \D u_\infty\big) 
 \\
 = &\left[\int_0^1 \p_\eta \bar \xi \big(\cdot, u_\infty+ \la \psi, \D u_\infty+ \la \D \psi \big) \, \mathrm d \la \right] \psi
 \\
   &\   + \left[\int_0^1 \p_{\mathrm p} \bar \xi \big(\cdot, u_\infty+ \la \psi, \D u_\infty+ \la \D \psi \big) \, \mathrm d \la\right] \cdot \D\psi,
  \end{split}
\eeq
a.e.\ on $\Om$. Note also that $|\psi|= |\D \psi| = 0$ $\mH^{n-1}$-a.e.\ on $\p\Om$, and $\psi \in \mathrm W_0^{2,q}(\Om)$ for all $q\in(1,\infty)$. Let $\psi$ be extended by zero on $\R^n\set\Om$ to a compactly supported function $\psi \in\mathrm W^{2,q}_c(\R^n)$. Then, by \eqref{E3}-\eqref{E4}, there exist functions $\mathfrak B \in\mathrm L^{\infty}_c(\R^n;\R^n)$ and $\mathfrak C \in\mathrm L^{\infty}_c(\R^n)$, such that $\psi$ satisfies the next PDE in the strong sense:
\[
\mathrm A\:\D^2 \psi  + \mathfrak B \cdot \D \psi + \mathfrak C \psi = 0, \ \ \ \text{ a.e.\ in }\R^n.
\] 
By \cite[Th.\ (1.7) and Rem.\ (1.8)]{HS}, we infer that $\mL^n(\{\psi=0\})=0$, which contradicts that $ \psi \equiv 0$ on $\R^n \set \Om$. Therefore, the proof of Theorem \ref{theorem2} is complete.
\qed
\ms

Now we will establish the proofs of Theorems \ref{theorem6} and  \ref{theorem7} together, by utilising  Theorem \ref{theorem8} (established in Subsection \ref{subsection2.5}). 

\ms

\noi \textbf{Proofs of Theorems} \ref{theorem6} \textbf{and} \ref{theorem7}. (i) Let $u_\infty \in \W^{2,\infty}_{u_0}(\Om)$ be a narrow minimiser of \eqref{1.1}. By Definition \ref{definition1}, this means there exists $\rho>0$ such that
\[
\ \ \ \F_\infty (u_\infty) \leq \F_\infty(u_\infty+\psi), \ \ \  \psi \in \mB_\rho^{\mathrm W^{1,\infty}(\Om)} \cap \W_0^{2,\infty}(\Om). 
\]
Let us consider now the functional $\F_{p,1,u_\infty}$ given by \eqref{1.19} for $\epsilon=1$ (and $\bar u = u_\infty$). For $p\geq p_0 := \bar p +n+1/c^3$, we set
\[
\ \ \ \ \ \ \left\{ \ \
\begin{split}
\E^\infty_{u_\infty}(\Om) \, &:=\, \overline \mB_{\rho/2}^{\mathrm W^{1,\infty}(\Om)}(u_\infty)   \cap \W^{2,\infty}_{u_0} (\Om), 
\\
\E^p_{u_\infty}(\Om) \, &:=\, \overline \mB_{\rho/2}^{\mathrm W^{1,\infty}(\Om)}(u_\infty) \cap\mathrm W^{2,p}_{u_0} (\Om),\ \ \ \ \ p_0 \leq p<\infty.
\end{split}
\right.
\]
Then, the class of sets $\{ \E^p_{u_\infty}(\Om)\}_{p_0 \leq p \leq \infty}$ satisfies the properties set out in assumption \eqref{1.25}. In suffices to confirm only that $\E^p_{u_\infty}(\Om)$ is weakly closed in $\mathrm W^{2,p}_{u_0} (\Om)$ when $p\geq p_0$: indeed, since $\E^p_{u_\infty}(\Om)$ is convex and norm-closed in $\mathrm W^{2,p}(\Om)$, its weak closedness is a consequence of Mazur's theorem. Therefore, by Theorem \ref{theorem8}, $\F_{p,1,u_\infty}$ has a minimiser $u_p \in \E^p_{u_\infty}(\Om)$, namely 
\[
\ \ \ \ \F_{p,1,u_\infty}(u_p) \leq \, \F_{p,1,u_\infty}(v), \ \ \  \ v \in \E^p_{u_\infty}(\Om).
\]
By setting $\phi := v - u_\infty \in\mathrm W^{2,p}_0(\Om)$, we can rewrite the above condition as
\[
\ \ \ \F_{p,1,u_\infty}(u_p) \leq \, \F_{p,1,u_\infty}(u_\infty + \phi), \ \ \ \ \phi \in  \overline \mB_{\rho/2}^{\mathrm W^{1,\infty}(\Om)}  \cap\mathrm W^{2,p}_0 (\Om).
\]
Note further that, again by Theorem \ref{theorem8}, there exists a sequence $(p_j)_1^\infty$ along which $\F_{p_j,1,u_\infty}(u_p) \larrow \F_\infty(u_\infty)$, and also $u_{p_j} \weak u_\infty$ in $\mathrm W^{2,q}_{u_0}(\Om)$ for any fixed $q\in (n,\infty)$, both as $j\to \infty$. By Morrey's theorem, we have $u_{p_j}  \larrow u_\infty$ in $\mathrm W^{1,\infty}(\Om)$, as $j\to\infty$. This implies for $j\geq j_0$ sufficiently large, that
\[
u_{p_j} \in \mB_{\rho/4}^{\mathrm W^{1,\infty}(\Om)}(u_\infty) \cap\mathrm W^{2,p}_{u_0} (\Om),
\]
and therefore $u_{p_j}$ is a point in the strong interior of $\E^p_{u_\infty}(\Om)$. It follows that 
\[
\ \ \ \ \F_{p,1,u_\infty}(u_p) \leq \, \F_{p,1,u_\infty}(u_p + \psi), \ \ \ \psi \in \,  \mB_{\rho/4}^{\mathrm W^{1,\infty}(\Om)} \cap\mathrm W^{2,p}_{0} (\Om),
\]
which establishes that $u_p$ is indeed a narrow minimiser of $\F_{p,1,u_\infty}$ in $\mathrm W^{2,p}_{u_0} (\Om)$.

\ms

\noi (ii) Let $u_p$ be the narrow minimiser of $\F_{p,1,u_\infty}$ in $\mathrm W^{2,p}_{u_0} (\Om)$ from part (i). Since the corresponding set of variations
\[
\mB_{\rho/4}^{\mathrm W^{1,\infty}(\Om)}(u_\infty) \cap\mathrm W^{2,p}_{u_0} (\Om)
\]
is strongly (norm) open in $\mathrm W^{2,p}_{u_0} (\Om)$, and $F$ satisfies the growth and regularity conditions given in \eqref{1.6}-\eqref{1.7} and \eqref{1.17}, it follows that $\F_{p,1,u_\infty}$ is Gateaux differentiable in $\mathrm W^{2,p}_{u_0} (\Om)$ at $u_p$, therefore $u_p$ is a weak solution to the Euler-Lagrange equations, whose weak formulation is given by \eqref{1.21}-\eqref{1.22}. Let $f_p$ be defined by \eqref{1.23} and let us also define the functions
\beq
\label{2.28}
\left\{\ \ 
\begin{split}
\mathrm K_p & : = \frac{\p_\eta F(\mathrm J^2u_p)}{\p_\xi F(\mathrm J^2u_p)}\ : \ \ \Om \larrow\R,
\\
\mathrm L_p & : = \frac{\p_{\mathrm p} F(\mathrm J^2u_p)}{\p_\xi F(\mathrm J^2u_p)} \ : \ \ \Om \larrow\R^n.
\end{split}
\right.
\eeq
By \eqref{1.7} and \eqref{1.17}, we have the bounds
\[
\left\{ \ \ 
\begin{split}
|\mathrm K_p|  = \frac{\big|\p_\eta F(\mathrm J^2u_p)\big|}{\p_\xi F(\mathrm J^2u_p)} \leq  \frac{C\big(\|u_p\|_{\mathrm W^{1,\infty}(\Om)}\big)}{c} \big(1 + |\mathrm{A}\:\D^2 u_p|\big),
\\
|\mathrm L_p|  = \frac{\big|\p_{\mathrm p} F(\mathrm J^2u_p)\big|}{\p_\xi F(\mathrm J^2u_p)} \leq  \frac{C\big(\|u_p\|_{\mathrm W^{1,\infty}(\Om)}\big)}{c} \big(1 + |\mathrm{A}\:\D^2 u_p|\big).
\end{split}
\right.
\]
In view of Lemma \ref{lemma13}, the above estimate implies that, for any fixed $q\in(1,\infty)$, $\{\mathrm K_p\}_{p\geq p_0}$ and and $\{|\mathrm L_p|\}_{p\geq p_0}$ are bounded in $\mathrm L^q(\Om)$. By the identities
\beq
\label{KpLp}
\left\{\ \ \
\begin{split}
\mathrm K_p f_p \, &=\, e_p^{1-p} \big(|F|^{p-2}F \p_\eta F\big)(\mathrm J^2u_p),
\\
\mathrm L_p f_p \, &=\, e_p^{1-p} \big(|F|^{p-2}F \p_{\mathrm p} F\big)(\mathrm J^2u_p),
\end{split}
\right.
\eeq
which follow by \eqref{1.23} and \eqref{2.28}, and in view of our notation for reduced jets and reduced derivatives of $F$, we can rewrite \eqref{1.21}-\eqref{1.23} as
\beq
\label{2.29}
\ \ \div\big(\mathrm{A}\D f_p -\mathrm L_p f_p\big)+ \mathrm K_p f_p +u_p-u_\infty =0,\ \  \text{ weakly in $(\mathrm W^{2,p}_0(\Om))^*$}.
\eeq
We now show that $\{f_p\}_{p\geq p_0}$ and $\{\mathrm K_p f_p\}_{p\geq p_0}$ are bounded in $\mathrm L^1(\Om)$, and $\{\mathrm L_p f_p\}_{p\geq p_0}$ is bounded in $\mathrm L^1(\Om;\R^n)$, at least along a subsequence of indices $(p_j)_1^\infty$. By \eqref{1.22}, \eqref{1.23} and \eqref{1.17}, if $p'=p/(p-1)$ is the conjugate exponent, we estimate
\[
\begin{split}
\| f_ p \|_{\mathrm L^{p'}(\Om)} &= \, e_p^{1-p} \bigg( \, {\av_\Om}  \Big( \big(|F|^{p-1} \p_\xi F\big) (\mathrm J^2 u_p)\Big)^{p'} \mathrm d \mL^n \! \bigg)^{\!\!\frac{1}{p'}}
\\
&\leq \, e_p^{1-p} C\Big(\|u_\infty\|_{\mathrm W^{1,\infty}(\Om)}+1 \Big)\bigg( \, {\av_\Om}\big|F (\mathrm J^2 u_p)\big|^p\,\mathrm d \mL^n \! \bigg)^{\!\!\frac{p-1}{p}}
\\
&=\, C\Big(\|u_\infty\|_{\mathrm W^{1,\infty}(\Om)}+1 \Big),
\end{split}
\]
for all sufficiently large $p\geq p_0$ (made possible by \eqref{1.20}). Further, by \eqref{2.27A} and \eqref{1.16}, we have the estimate
\beq
\label{2.31}
|\mathrm{A}\:\D^2 w| \leq \frac{|F(\mathrm J^2w)| +(|w|+|\D w|)^\al + 1}{c} \ \ \text{ a.e.\ on }\Om,
\eeq
for any $w \in\mathrm W^{2,1}_{\mathrm{loc}}(\Om)$. Therefore, by \eqref{1.17} and \eqref{KpLp}, we may estimate
\[
\begin{split}
\| \mathrm K_p f_ p \|_{\mathrm L^1(\Om)} &\leq \, e_p^{1-p}   {\av_\Om}  \big| \big(|F|^{p-1} \p_\eta F \big) (\mathrm J^2 u_p)\big|  \,\mathrm d \mL^n  
\\
&\leq e_p^{1-p} C\Big(\|u_\infty\|_{\mathrm W^{1,\infty}(\Om)}+1 \Big)   {\av_\Om}\big|F (\mathrm J^2 u_p)\big|^{p-1}  \Big( 1+ |\mathrm{A}\:\D^2 u_p| \Big)\,\mathrm d \mL^n .
\end{split}
\]
This implies
\[
\begin{split}
\| \mathrm K_p f_ p \|_{\mathrm L^1(\Om)} & \leq e_p^{1-p} C\Big(\|u_\infty\|_{\mathrm W^{1,\infty}(\Om)}+1 \Big)   \bigg( \, {\av_\Om}\big|F (\mathrm J^2 u_p)\big|^p\,\mathrm d \mL^n \! \bigg)^{\!\!\!\frac{p-1}{p}} 
\\
&\ \ \ \ \cdot \bigg( \, {\av_\Om} \big( 1  +|\mathrm{A}\:\D^2 u_p| \big)^p\,\mathrm d \mL^n \! \bigg)^{\!\!\frac{1}{p}}
\\
&= C\Big(\|u_\infty\|_{\mathrm W^{1,\infty}(\Om)}+1 \Big) \bigg( \, {\av_\Om} \big( 1  +|\mathrm{A}\:\D^2 u_p| \big)^p\,\mathrm d \mL^n \! \bigg)^{\!\!\frac{1}{p}}.
\end{split}
\]
Hence, in view of \eqref{2.31}, the above estimate yields (for large enough $p$),
\[
\begin{split}
\| \mathrm K_p f_ p \|_{\mathrm L^1(\Om)} &\leq C\Big(\|u_\infty\|_{\mathrm W^{1,\infty}(\Om)}+1 \Big)  
\\
&\ \ \ \cdot  \bigg[ \,\, {\av_\Om} \bigg(1+\frac{1}{c} + \frac{\big(|u_p|+|\D u_p|\big)^\al }{c} + \frac{|F(\mathrm{J}^2 u_p)|}{c} \bigg)^{\!p}\,\mathrm d \mL^n  \bigg]^{\!\!\frac{1}{p}}
\\
& \leq C\Big(\|u_\infty\|_{\mathrm W^{1,\infty}(\Om)}+1 \Big) \bigg[ 1+\frac{1}{c} +\frac{\|u_\infty\|_{\mathrm W^{1,\infty}(\Om)}^\al +1}{c} +\frac{e_p}{c}  \bigg] 
\\
& = C\Big(\|u_\infty\|_{\mathrm W^{1,\infty}(\Om)}+1 \Big)\frac{1}{c} \bigg( 2 +c +\|u_\infty\|_{\mathrm W^{1,\infty}(\Om)}^\al  + e_p   \bigg) .
\end{split}
\]
Again by \eqref{1.17}, since $\p_\eta F$ and  $\p_{\mathrm p} F$ satisfy the same bounds, by  \eqref{KpLp} we have that the same estimate is true also for $\mathrm L_p f_ p$:
\[
\begin{split}
\| \mathrm L_p f_ p \|_{\mathrm L^1(\Om)} \leq C\Big(\|u_\infty\|_{\mathrm W^{1,\infty}(\Om)}+1 \Big)\frac{1}{c} \bigg( 2 +c +\|u_\infty\|_{\mathrm W^{1,\infty}(\Om)}^\al  + e_p   \bigg) .
\end{split}
\]
The above inequalities implies the desired bound for sufficiently large $p\geq p_0$, upon nothing that, by \eqref{1.22}, \eqref{1.19} and \eqref{1.20}, we have
\beq
\label{2.32}
e_p \, =\, \F_{p,1,u_\infty}(u_p) - \frac{1}{2} \| u_p -u_\infty \|^2_{\mathrm L^2(\Om)} \larrow \F_\infty (u_\infty),
\eeq
as $p_j\to\infty$. In view of the previous uniform $\mathrm L^1$-bounds, there exist Radon measures
\[
\varphi_\infty \in \mM(\overline{\Om}), \ \  \kappa_\infty \in \mM(\overline{\Om}), \ \ \la_\infty \in \mM(\overline{\Om};\R^n),
\]
and by applying Lemma \ref{lemma15} (Subsection \ref{subsection2.6}), there exists 
\[
f_\infty \in \, \mathrm L^1(\Om)\cap \bigcap_{1<q<n'} \mathrm L^q_{\mathrm{loc}}(\Om),
\]
where $n'=n/(n-1)$, and a further subsequence, denoted again by $(p_j)_1^\infty$, along which
\[
\left\{ \ \
\begin{split}
& f_p \larrow f_\infty, \hspace{98pt} \text{ in }\mathrm L^q_{\mathrm{loc}}(\Om),\ q\in(1,n'),
\\
& f_p \mL^n \LL_\Om \weakstar \, \varphi_\infty  ,  \hspace{73pt} \text{ in }\mM(\overline{\Om}),
\\
& \big(\mathrm K_p f_p +u_p-u_\infty\big)  \mL^n \LL_\Om \weakstar \, \kappa_\infty,  \ \text{ in }\mM(\overline{\Om}),
\\
&(\mathrm L_p f_p)  \mL^n \LL_\Om \weakstar \, \la_\infty,   \hspace{56pt}  \text{ in }\mM(\overline{\Om};\R^n),
\end{split}
\right.
\]
as $p_j \to\infty$. Further, since $(\mathrm K_{p_j})_1^\infty$ and $(|\mathrm L_{p_j}|)_1^\infty$ are both bounded in $\mathrm L^q(\Om)$ for all $q\in(1,\infty)$, there exist
\beq 
\label{2.34A}
\ \ \ \ \mathrm K_\infty \in \bigcap_{1<q<\infty} \mathrm L^q(\Om),\ \ \ \mathrm L_\infty \in \bigcap_{1<q<\infty} \mathrm L^q(\Om;\R^n),
\eeq
such that $\mathrm K_p \weak \mathrm K_\infty$ in $\mathrm L^q(\Om)$, and also $\mathrm L_p \weak \mathrm L_\infty$ in $\mathrm L^q(\Om;\R^n)$, as $p_j\to\infty$ for all $q\in(1,\infty)$. Since $f_p \larrow f_\infty$ in $\mathrm L^q_{\mathrm{loc}}(\Om)$ for all $q\in(1,n')$, we have $\mathrm K_p f_p \weak \mathrm K_\infty f_\infty$ in $\mathrm L^r_{\mathrm{loc}}(\Om)$ for all $r\in (1,n')$ as $p_j\to\infty$, and similarly $\mathrm L_p f_p \weak \mathrm L_\infty f_\infty$ in $\mathrm L^r_{\mathrm{loc}}(\Om;\R^n)$ for all $r\in (1,n')$, as $p_j\to\infty$. By \eqref{1.20} we also have that $u_p \larrow u_\infty$ in $\mathrm C(\overline{\Om})$, as $p_j\to\infty$. Hence, by the uniqueness of limits,  the measures $\varphi_\infty,\kappa_\infty,\la_\infty$ are absolutely continuous with respect to the Lebesgue measure restricted to $\Om$, and in addition
\[
\left\{ \ \
\begin{split}
\varphi_\infty  \LL_\Om &= f_\infty \mL^n \LL_\Om,
\\
\kappa_\infty  \LL_\Om &= (\mathrm K_\infty f_\infty) \mL^n \LL_\Om,
\\
\la_\infty  \LL_\Om &= (\mathrm L_\infty f_\infty) \mL^n \LL_\Om.
\end{split}
\right.
\]
We therefore have for any fixed $q\in (1,n')$ the following convergences
\beq
\label{2.33}
\left\{ \ \
\begin{split}
& f_p \larrow f_\infty, \, \hspace{190pt} \text{ in }\mathrm L^q_{\mathrm{loc}}(\Om),
\\
& f_p \mL^n \LL_\Om \weakstar \, f_\infty  \mL^n\LL_\Om + \varphi_\infty \LL_{\p\Om},  \ \hspace{94pt} \text{ in }\mM(\overline{\Om}),
\\
& \big(\mathrm K_p f_p +u_p-u_\infty\big)  \mL^n \LL_\Om \weakstar \, (\mathrm K_\infty f_\infty) \mL^n \LL_\Om + \kappa_\infty \LL_{\p\Om},  \ \text{ in }\mM(\overline{\Om}),
\\
& (\mathrm L_p f_p)  \mL^n \LL_\Om \weakstar \, (\mathrm L_\infty f_\infty) \mL^n \LL_\Om + \la_\infty \LL_{\p\Om},  \hspace{57pt} \text{ in }\mM(\overline{\Om};\R^n),
\end{split}
\right.
\eeq
as $p_j\to \infty$. By passing to the limit as $j\to\infty$ in the PDE \eqref{2.29} along the sequence $(p_j)_1^\infty$, the previous considerations allow us to deduce 
\beq
\label{2.34}
\ \ \ \div\big(\mathrm A \D f_\infty -\mathrm L_\infty f_\infty\big) +\mathrm K_\infty f_\infty =0 \ \ \ \text{ in }(\mathrm{C}^\infty_c(\Om))^*,
\eeq
namely in the sense of distributions. By applying Lemma \ref{lemma2}(ii), it follows that \eqref{2.34} is also satisfied weakly in $(\W^{2,\infty}_0(\Om))^*$. 

\ms


We now show that, if \eqref{1.17A} is satisfied, then $\varphi_\infty \LL{\p\Om}=0$. (This is actually the only place in the paper that the $\mathrm{C}^2$ regularity assumption on $\p\Om$ is utilised).

\begin{claim} \label{claim14} If $\p\Om$ is $\mathrm{C}^2$, we have that $\varphi_\infty \LL_{\p\Om}=0$.
\end{claim}

\BPC \ref{claim14}. By equations \eqref{2.34} and \eqref{2.29}, we have
\beq
\label{2.36A}
\begin{split}
\div \big(\mathrm{A}\D f_p -\mathrm{A}\D  f_\infty\big) &= \mathrm K_\infty f_\infty - \mathrm K_p f_p  +u_\infty-u_p
\\
&\ \ \ -\div\big( \mathrm L_\infty f_\infty - \mathrm L_p f_p \big)  ,
\end{split}
\eeq
weakly in $(\W^{2,\infty}_0(\Om))^*$. Since by assumption $\p\Om$ is $\mathrm{C}^2$, there exists $r_0>0$ such that the distance function $\dist(\cdot,\p\Om)$ is $\mathrm{C}^2$ when restricted to the inner tubular neighbourhood $\overline{\Om}\cap \{\dist(\cdot,\p\Om) < r_0\}$. Let $d \in \mathrm{C}^2(\overline{\Om})$ be any extension of the distance function $\dist(\cdot,\p\Om)$ from $\overline{\Om} \cap \{\dist(\cdot,\p\Om) < r_0\}$ to $\overline{\Om}$. Let us fix also any $g \in \mathrm{C}^2(\overline{\Om})$. We set 
\[
\psi  :=  d^2 g \, \in  \mathrm{C}^2(\overline{\Om}) \cap \W^{2,\infty}_0(\Om).
\]
By testing $\psi$ in the weak formulation of \eqref{2.36A}, we have
\beq
\label{2.36B}
\begin{split}
\int_\Om (f_p  - f_\infty)\mathrm{A}\: \D^2\psi \, \mathrm d\mL^n &= \int_\Om \Big(\mathrm K_\infty f_\infty - \mathrm  K_p f_p +u_\infty -u_p \Big) \psi \, \mathrm d\mL^n
\\
&\ \ \ + \int_\Om \big(\mathrm L_\infty f_\infty - \mathrm  L_p f_p  \big)\cdot \D \psi \, \mathrm d\mL^n .
\end{split}
\eeq
By a computation, we have
\[
\D^2 \psi = 2d\Big( \D g \ot \D d + \D d \ot \D g + g\D^2d + \frac{d}{2} \D^2g\Big) + 2 g \D d \ot \D d, \ \text{ on }\overline{\Om},
\]
noting also that
\[
\D^2\psi \big|_{\p\Om} =\, 2g \, \mathrm{n} \ot  \mathrm{n},
\]
where $ \mathrm{n} =\D d$ is the outwards unit normal vector field on $\p\Om$. By \eqref{2.33}, since $\psi \in \mathrm C_0(\overline{\Om})$ and $\D\psi \in \mathrm C_0(\overline{\Om};\R^n)$, we may pass to the limit as $p_j\to\infty$ in \eqref{2.36B}, to obtain
\[
 \int_{\p\Om} (g \, \mathrm{A}\: \mathrm{n} \ot  \mathrm{n}) \, \mathrm d\varphi_\infty  =0,
\]
for any $g \in \mathrm{C}^2(\overline{\Om})$. Fix now $\e>0$ and an arbitrary $h \in \mathrm{C}(\p \Om)$. We may extend $h$ to a function in  $\mathrm{C}(\overline{\Om})$, denoted again by $h$ (e.g.\ by using its McShane-Whitney extension $x \mapsto \max_{y\in\p\Om}\{h(y)-|x-y|\}$). It follows that there exists $g \in \mathrm{C}^2(\overline{\Om})$ such that 
\[
\left\| \frac{h}{\mathrm{A}\: \mathrm{n} \ot  \mathrm{n}} -g \,\right\|_{\mathrm C(\p\Om)} < \frac{\e}{\|A\|}.
\]
In conclusion,
\[
\left| \int_{\p\Om} h \, \mathrm d \varphi_\infty \right| \leq  \int_{\p\Om} \big| h-g \, \mathrm{A}\: \mathrm{n} \ot  \mathrm{n} \big| \, \mathrm d \|\varphi_\infty\| <\e \|\varphi_\infty\|(\p\Om),
\]
which yields that $\varphi_\infty \LL_{\p\Om}=0$, as claimed.       \qed
\ms


Note now that, by the definition \eqref{1.23}, and the fact that the inverse function of $t\mapsto |t|^{p-2}t$ on $\R$ is given by $s\mapsto \sgn(s) |s|^{1/(p-1)}$, we have the identity
\beq
\label{2.35}
\ \ F(\mathrm J^2u_p) \, = \, e_p \sgn(f_p) |f_p|^{\frac{1}{p-1}} \big(\p_\xi F(\mathrm J^2u_p)\big)^{-\frac{1}{p-1}}  \ \ \ \text{ a.e.\ on }\Om.
\eeq 
Let us fix $x\in\Om\set \{f_\infty = 0\}$. Since by \eqref{2.33} we have that $f_p \larrow f_\infty$ a.e.\ on $\Om$ as $p_j\to\infty$, it follows that $f_p(x) \in \big[ |f_\infty(x)|/2, 2|f_\infty(x)| \big]$ for large enough $j\in\N$. By assumptions \eqref{1.7} and \eqref{1.17}, as well as \eqref{2.32}, we may pass to the a.e.\ limit in \eqref{2.35} to obtain
\beq
\label{2.36}
\ \ F(\mathrm J^2u_p) \larrow \sgn(f_\infty)\F(u_\infty)  \ \ \ \text{ a.e.\ on }\Om\set\{f_\infty = 0\},
\eeq 
as $p_j\to\infty$. Again by \eqref{1.7}, the function $F$ is invertible with respect to last variable, therefore \eqref{2.36} implies
\[
\ \ \mathrm{A}\:\D^2u_p \larrow \big(F(\cdot,u_\infty,\D u_\infty,\_)\big)^{-1}\big(\sgn(f_\infty)\F(u_\infty)\big)  \ \ \ \text{ a.e.\ on }\Om\set \{f_\infty = 0\}, 
\]
as $p_j\to\infty$. However, by \eqref{1.20} we know that $\mathrm{A}\:\D^2u_p \weak \mathrm{A}\:\D^2u_\infty$ in $\mathrm L^q(\Om)$ as $p_j\to\infty$, for all $q\in(1,\infty)$. By the Vitali convergence theorem and the uniqueness of limits, it follows that 
\beq
\label{2.37}
\mathrm{A}\:\D^2u_p \larrow \mathrm{A}\:\D^2u_\infty = \big(F(\cdot,u_\infty,\D u_\infty,\_)\big)^{-1}\big(\sgn(f_\infty)\F(u_\infty)\big), 
\eeq 
in $\mathrm L^q\big( \Om\set\{f_\infty = 0\} \big)$ for all $q\in(1,\infty)$, and additionally a.e.\ on $\Om\set\{f_\infty = 0\}$, as $p_j\to\infty$. The convergence in \eqref{2.37} has two consequences. On the one hand, by applying the function $F(\cdot,u_\infty,\D u_\infty,\_)$ to \eqref{2.37}, we have
\beq
\label{2.38}
\ \ F(\mathrm J^2 u_\infty) = \sgn(f_\infty)\F(u_\infty) \ \ \ \text{ a.e.\ on }\Om\set\{f_\infty = 0\}.
\eeq
The equality \eqref{2.38} implies that the pair $(u_\infty,f_\infty)$ satisfies the equation \eqref{1.8} a.e.\ on $\Om\set\{f_\infty = 0\}$. Further, by \eqref{2.28} and \eqref{1.20}, \eqref{2.37} implies
\beq
\label{2.39}
\left.\ \ 
\begin{split}
\mathrm K_p &= \frac{\p_\eta F(\mathrm J^2u_p)}{\p_\xi F(\mathrm J^2u_p)} \larrow \frac{\p_\eta F(\mathrm J^2u_\infty)}{\p_\xi F(\mathrm J^2u_\infty)}
\\
\mathrm L_p &= \frac{\p_{\mathrm p} F(\mathrm J^2u_p)}{\p_\xi F(\mathrm J^2u_p)} \larrow \frac{\p_{\mathrm p} F(\mathrm J^2u_\infty)}{\p_\xi F(\mathrm J^2u_\infty)}\ 
\end{split}
\right\} \text{ a.e.\ on }\Om\set\{f_\infty = 0\}, \text{ as $p_j\to\infty$}.
\eeq
Recalling that $\mathrm K_p \weak \mathrm K_\infty$ in $\mathrm L^q(\Om)$ and also $\mathrm L_p \weak \mathrm L_\infty$ in $\mathrm L^q(\Om;\R^n)$, as $p_j\to\infty$, for all $q\in (1,\infty)$, the Vitali convergence theorem yields that the modes of convergence in \eqref{2.39} are valid in $\mathrm L^q\big( \Om\set\{f_\infty = 0\}\big)$, for all $q\in (1,\infty)$. Further, by the uniqueness of limits, we infer that $\mathrm K_\infty,\mathrm L_\infty$ are characterised as
\beq
\label{2.40}
\left.
\begin{split}
\mathrm K_\infty &= \frac{\p_\eta F(\mathrm J^2u_\infty)}{\p_\xi F(\mathrm J^2u_\infty)}\ 
\\
\mathrm L_\infty &= \frac{\p_{\mathrm p} F(\mathrm J^2u_\infty)}{\p_\xi F(\mathrm J^2u_\infty)}
\end{split}
\right\} \ \text{ a.e.\ on }\Om\set\{f_\infty = 0\}.
\eeq
In view \eqref{2.40} and \eqref{2.34}, it follows that the PDE \eqref{1.9} is indeed satisfied by the pair $(u_\infty,f_\infty)$. By Proposition \eqref{proposition19} applied to \eqref{1.9}, we have that $f_\infty \in \mathrm W^{1,q}_{\mathrm{loc}}(\Om)$ for all $q\in(1,\infty)$, and $\{f_\infty = 0\}$ is a countably $(n-1)$-rectifiable subset of $\Om$. In particular, $\{f_\infty = 0\}$ is a Lebesgue nullset. We may further assume that $\F_\infty(u_\infty)>0$. If instead $\F_\infty(u_\infty)=0$, then \eqref{1.8} is trivially satisfied a.e.\ on $\Om$ for any $f_\infty$, whilst \eqref{1.9} is also satisfied by choosing any non-trivial solution $f_\infty \neq 0$, guaranteed to exist by Proposition \eqref{proposition19}. The proofs of Theorems \ref{theorem6} and \ref{theorem7} will therefore be completed once we prove that $f_\infty \not\equiv 0$. The remainder of the proof is therefore devoted to showing $f_\infty \not\equiv 0$, by using assumption \eqref{1.17A}. (This is the only point where the $\mathrm{C}^2$ regularity assumption for $u_0|_{\p\Om}$ is utilised.) As a first step, we  establish the estimate
\beq
\label{2.44}
F(\mathrm J^2u_p) \p_\xi F(\mathrm J^2u_p) \Big( \mathrm{A}\:\D^2u_p - \bar \xi(\cdot,u_p,\D u_p) \Big) \geq \Big(\frac{c}{C_\infty}\Big)^{\!2}\big(F (\mathrm J^2u_p)\big)^2,
\eeq
a.e.\ on $\Om$, for all sufficiently large $p\geq p_0$, where 
\[
C_\infty := C\big(1+ \|u_\infty\|_{\mathrm W^{1,\infty}(\Om)} \big).
\]
To establish \eqref{2.44}, we use \eqref{1.7} and \eqref{1.17} to estimate as follows:
\[
\begin{split}
F(x,\eta,\xi,\mathrm p) \big( \xi - \bar \xi(x,\eta,\mathrm p)\big) &= \Big[F(x,\eta,\mathrm p,\xi) -F\big(x,\eta,\mathrm p,\bar \xi(x,\eta, \mathrm p)\big) \Big] \big( \xi - \bar \xi(x,\eta, \mathrm p)\big) 
\\
&= \bigg( \int_{\bar \xi(x,\eta,\mathrm p)}^\xi \p_\xi F(x,\eta,\mathrm p, \cdot)\, \mathrm d \mL^1 \!\bigg) \big( \xi - \bar \xi(x,\eta,\mathrm p)\big)
\\
&\geq c \big| \xi - \bar \xi(x,\eta,\mathrm p)\big|^2
\\
& \geq \frac{c}{C(|\eta|+|\mathrm p|)^2}\big( F(x,\eta,\mathrm p,\xi) \big)^2,
\end{split}
\]
which, again by  \eqref{1.7}, yields
\beq
\label{2.45}
F(x,\eta,\mathrm p,\xi) \p_\xi F(x,\eta,\mathrm p,\xi) \big( \xi - \bar \xi(x,\eta,\mathrm p)\big)  \geq \frac{c^2\big( F(x,\eta,\mathrm p,\xi) \big)^2}{C(|\eta|+|\mathrm p|)^2},
\eeq
for a.e.\ $x\in \Om$  and all $(\eta,\mathrm p,\xi) \in \R\by \R^n\by \R$. Inequality \eqref{2.45} readily implies \eqref{2.44}. Next, we establish the following inequality
\beq
\label{2.46}
\left\{ \ \ \
\begin{split}
e_p \Big(\frac{c}{C_\infty}\Big)^{\!2}  & \leq \, {\av_\Om}\,  f_p \Big[ \mathrm{A}\: \D^2 w - \bar \xi (\cdot, u_p,\D u_p)\Big] \,\mathrm d \mL^n
\\
& \ \ \  + \, {\av_\Om}\,  f_p \Big[ \mathrm K_p(w-u_p)-\mathrm L_p\cdot \D(w-u_p)\Big] \,\mathrm d \mL^n
\\
& \ \ \ + \| w-u_p\|_{\mathrm L^1(\Om)} \| u_\infty -u_p\|_{\mathrm L^{\infty}(\Om)}, \phantom{\Big[}
\end{split}
\right.
\eeq
for any fixed $w\in \W^{2,\infty}_{u_0}(\Om)$, where $e_p$ is given by \eqref{1.22}, and satisfies \eqref{2.32} as $j\to\infty$ along a subsequence $(p_j)_1^\infty$. To this aim, let $w\in \W^{2,\infty}_{u_0}(\Om)$ be fixed. Since $w-u_p \in\mathrm W^{2,p}_0(\Om)$, by testing this in the weak formulation of \eqref{2.29} and employing  \eqref{2.44} and \eqref{1.23}, we have
\[
\begin{split}
\ \ \ &{\av_\Om}\,  \Big[\Big(\mathrm K_p f_p +(u_p-u_\infty)\Big)(w-u_p) -   (\mathrm L_p f_p) \cdot \D(w-u_p) \Big]\,\mathrm d \mL^n
\\
&\  + {\av_\Om}\,  f_p \Big[ \mathrm{A}\: \D^2 w- \bar \xi (\cdot, u_p,\D u_p) \Big] \,\mathrm d   \mL^n = \, {\av_\Om}\,  f_p \Big[ \mathrm{A}\: \D^2 u_p - \bar \xi (\cdot, u_p,\D u_p)\Big] \,\mathrm d \mL^n  
\\
&\hspace{165pt} \geq \,e_p^{1-p}\Big(\frac{c}{C_\infty}\Big)^{\!2}  {\av_\Om}\, \big(F (\mathrm J^2u_p)\big)^p\,\mathrm d \mL^n .  
\end{split}
\]
The above inequality readily implies \eqref{2.46}. 

By noting that assumption \eqref{1.17A} requires that $u_0|_{\p\Om}$ is $\mathrm{C}^2$ on the closed set $\p\Om$, Whitney's extension theorem \cite[Th.\ I, p.\ 65]{W}, yields the existence of some $u_* \in \mathrm{C}^2(\R^n)$ such that $u_0 =u_*$ on $\p\Om$, and also $\D u_0 = \D u_*$ on $\p\Om$. For any $\e>0$, we set $\Om^\e := \Om \cap \{\dist(\cdot, \p\Om)<\e\}$. Since $\| u_* -u_0 \|_{\mathrm C^1(\Om^\e)}\larrow 0$, as $\e\to0$, by using a standard regularisation argument involving partitions of unity, we can construct a function $u_\e \in \mathrm{C}^2(\overline{\Om})$ satisfying that $u_\e \equiv u_*$ in a inner neighbourhood of the boundary, and $u_\e$ equals the standard mollification of $u_\infty$ (for small enough regularisation parameter), on the complement of such an inner neighbourhood. Then, by construction $u_\e$ satisfies $u_\infty =u_\e$ on $\p\Om$ and also $\D u_\infty = \D u_\e$ on $\p\Om$, with the additional property $\| u_\infty - u_\e \|_{\mathrm W^{1,\infty}(\Om)}\larrow 0$, as $\e\to0$. By recalling that $\{\mathrm K_p f_p\}_{p> p_0}$ and $\{|\mathrm L_p f_p|\}_{p> p_0}$ are bounded in $\mathrm L^1(\Om)$, the choice $w:=u_\e$ in \eqref{2.46} yields
\[
\begin{split}
e_p \Big(\frac{c}{C_\infty}\Big)^{\!2} & \leq \, {\av_\Om}\,  f_p \Big[ \mathrm{A}\: \D^2 u_\e - \bar \xi (\cdot, u_p,\D u_p)\Big] \,\mathrm d \mL^n  \ \ \ \ \ \ \ \
\\
& \ \ \ + \| u_\e-u_p\|_{\mathrm W^{1,\infty}(\Om)}\Big[   \| u_\infty -u_p\|_{\mathrm L^{\infty}(\Om)} + \| \mathrm K_p f_p\|_{\mathrm L^1(\Om)} + \| \mathrm L_p f_p\|_{\mathrm L^1(\Om)} \Big] .
\end{split}
\]
By letting $p\to\infty$ along a subsequence in the above inequality, we may use that $f_p \mL^n\LL_\Om \weakstar f_\infty \mL^n \LL_\Om$ in $\mM(\overline{\Om})\cong (\mathrm{C}(\overline{\Om}))^*$ (recall Claim \ref{claim14}), and the fact that $\mathrm{A}\: \D^2 u_\e \in \mathrm{C}(\overline{\Om})$ and that $\bar \xi (\cdot, u_p,\D u_p) \larrow \bar \xi (\cdot, u_\infty,\D u_\infty)$ uniformly on $\overline{\Om}$ as $p\to\infty$ along a subsequence, to deduce the inequality
\[
\begin{split}
\ \ \ \F_\infty(u_\infty) \Big(\frac{c}{C_\infty}\Big)^{\!2} - o(1)  \leq \, {\av_\Om}\,  f_\infty \Big[ \mathrm{A}\: \D^2 u_\e - \bar \xi (\cdot, u_\infty,\D u_\infty)\Big] \,\mathrm d \mL^n, \ \ \text{ as }\e\to0.
\end{split}
\]
Since by assumption $\F_\infty(u_\infty)>0$, if $\e>0$ is chosen small enough, the left hand side in the above inequality becomes positive, and this implies $f_\infty \not\equiv 0$. In conclusion, the proof of Theorems \ref{theorem6} and \ref{theorem7} is complete.
\qed
\ms

We now conclude by establishing Theorem \ref{UniquenessTheorem}.

\BPT \ref{UniquenessTheorem}. (I) Let $u_\infty \in \W^{2,\infty}_{u_0}(\Om)$ be a narrow minimiser of the functional $\F_\infty$. By Definition \ref{definition1}, there exists $\rho>0$, such that
\[
\ \ \ \ \F_\infty(u_\infty) \leq \F_\infty(u_\infty+\psi), \ \ \ \ \psi \in \smash{\mB_{\rho}^{\mathrm W^{1,\infty}(\Om)} \cap \W^{2,\infty}_0(\Om)}.
\]
By Theorem \ref{theorem6}, it follows that there exists $f_\infty \in \mathrm L^1(\Om)\set\{0\}$, such that the pair $(u_\infty, f_\infty)$ satisfies the system of equations \eqref{1.8}-\eqref{1.9}. By Theorem \ref{theorem2}, it follows that $u_\infty$ is actually a \emph{strict} narrow minimiser, therefore there exists $\rho'>0$ such that
\[
\ \ \ \ \ \ \ \F_\infty(u_\infty) < \F_\infty(u_\infty+\psi), \ \ \ \ \psi \in \smash{\big(\mB_{\rho'}^{\mathrm W^{1,\infty}(\Om)} \cap \W^{2,\infty}_0(\Om)\big)} \set\{0\}.
\]
Hence, $u_\infty$ is the unique narrow minimiser of the functional  $\F_\infty$ (at the same energy level) within in the admissible class
\[
\mB_{\min\{\rho,\rho'\}}^{\mathrm W^{1,\infty}(\Om)}(u_\infty) \cap \W^{2,\infty}_{u_0}(\Om).
\]

\noi (II) Let $u_\infty \in \W^{2,\infty}_{u_0}(\Om)$ be a global minimiser of the functional $\F_\infty$ (Definition \ref{definition1}). By Theorem \ref{theorem6}, it follows that there exists $f_\infty \in \mathrm L^1(\Om)\set\{0\}$, such that the pair $(u_\infty, f_\infty)$ satisfies the system of equations \eqref{1.8}-\eqref{1.9}. By Theorem \ref{theorem2}, it follows that $u_\infty$ is a strict global minimiser (Definition \ref{definition1}). In conclusion, $u_\infty$ is the unique global minimiser of the functional  $\F_\infty$ in the space $\W^{2,\infty}_{u_0}(\Om)$.  \qed
\ms

\BPCOR \ref{corollary7}. Let $\mathrm J^1u(x):=\big(x,u(x),\D u(x)\big)$ denote the first order jet operator. Let us also set $\Om^*:=\Om\set\{f_\infty =0\}$. Since $\{f_\infty=0\}$ is nowhere dense, $\Om^*$ is an open dense set in $\Om$. By equation \eqref{1.8} and assumption \eqref{1.7}, we have
\[
\ \ \ \ \ \mathrm A \! : \! \D^2u_\infty = F(\mathrm J^1 u_\infty, \cdot)^{-1} \big(\F_\infty(u_\infty)\sgn(f_\infty)\big), \ \ \text{a.e.\ on }\Om^*.
\]
Since $u_\infty \in \mathrm C^{1,\ga}(\Om)$ for all $\ga\in(0,1)$, the right hand side is in $ \mathrm C^{0,\ga}(\Om^*)$. By the Schauder interior estimates (see e.g.\ \cite[Ch.\ 6]{GT}), we have that $u_\infty \in \mathrm C^{2,\ga}(\Om^*)$. Hence, since $F$ is by assumption $\mathrm C^1$, the right hand side is in $ \mathrm C^{1,\ga}(\Om^*)$, which implies $u_\infty \in \mathrm C^{3,\ga}(\Om^*)$. The conclusion for general $k\geq 2$ is a consequence of a standard bootstrap argument and finite induction. Finally, by Proposition \ref{proposition19}, the nodal set $\{f_\infty=0\}$ is countably $(n-1)$-rectifiable.   \qed
\ms

\BPCOR \ref{corollary10}. Since any increasing $\Phi : \R \larrow \R$ which satisfies $\Phi(-\,\cdot )=-\Phi$ commutes with the absolute value and the (essential) supremum, we have
\[
\bar \F_\infty (u) =  \underset{\Om}{\ess \sup}\hspace{1pt} \big|\Phi\big(F(\mathrm J^2 u)\big)\big| =\Phi\Big(\underset{\Om}{\ess \sup}\hspace{1pt} |F(\mathrm J^2 u)|\Big) = \Phi( \F_\infty (u)),
\]
yielding that $\F_\infty$ and $\bar \F_\infty$ have the same sets of global/narrow/local minimisers. 
\qed
\ms

\subsection*{Acknowledgement} N.K.\ would like to thank Daniel Tataru and Alexander Logunov for their insights on the satisfaction of the unique continuation principle and the size of nodal sets of solutions to the equations arising in this work.

\bibliographystyle{amsplain}

\begin{thebibliography}{30}
\bibitem{Ad} R. A. Adams, \emph{Sobolev spaces}, Volume 140 of Pure and applied mathematics, Academic Press, 2003.

\bibitem{AP} N. Ansini, F. Prinari, \emph{On the lower semicontinuity of supremal functional under differential constraints}, ESAIM - Control, Opt. and Calc. Var. 21(4), 1053-1075 (2015).


\bibitem{A1} G. Aronsson, \emph{Minimization problems for the functional $sup_x F(x,
f(x), f'(x))$ I, II, III}, Arkiv f\"ur Mat. (1965), 33 - 53; (1966), 409 - 431; (1969), 509 - 512.

\bibitem{A3} G. Aronsson, \emph{Extension of functions satisfying Lipschitz
 conditions}, Arkiv f\"ur Mat. 6 (1967), 551 - 561.

\bibitem{A4} G. Aronsson, \emph{On the partial differential equation $u_x^2 u_{xx} + 2u_x u_y u_{xy} + u_y^2 u_{yy} = 0$}, Arkiv f\"ur Mat. 7
(1968), 395 - 425.

\bibitem{AB} G. Aronsson, E.N. Barron, \emph{$\mathrm L^{\infty}$ Variational Problems with Running Costs and Constraints}, Appl. Math. Optimization 65, 53-90 (2012).

\bibitem{BJ} E.N. Barron, R. Jensen, \emph{Minimizing the $\mathrm L^{\infty}$ norm of the gradient with an energy constraint}, Comm. Partial Differential Equations 30, 10-12, 1741-1772 (2005).

\bibitem{BJW2} E. N. Barron, R. Jensen, C. Wang, \emph{Lower Semicontinuity of $\mathrm L^{\infty}$ Functionals} Ann. I. H. Poincar\'e AN 18, 4 (2001) 495-517.

\bibitem{B} C. B\"ar, \emph{Zero sets of solutions to semilinear elliptic systems of first order}, Invent. Math. 138(1), 183-202 (1999)

\bibitem{BBH} T. Beck, S. Becker-Kahn, B. Hanin, \emph{Nodal sets of smooth functions with finite vanishing order and p-sweepouts}, Calc. Var. 57, 140 (2018). 

\bibitem{BKR} V. I. Bogachev, N. V. Krylov, and M. R\"ockner, \emph{Elliptic and parabolic equations for measures}, Uspekhi Mat. Nauk 64 (2009), 5–116, translation in Russian Math. Surveys 64 (2009), 973-1078.

\bibitem{BDP} C. Brizzi, L. De Pascale, \emph{A property of Absolute Minimizers in $\mathrm L^{\infty}$  Calculus of Variations and of solutions of the Aronsson-Euler equation}, Advances in Differential Equations 28 (3/4), 287-310 (2023).

\bibitem{BK} L. Bungert, Y. Korolev, \emph{Eigenvalue problems in $\mathrm L^{\infty}$: optimality conditions, duality, and relations with optimal transport}, Communications of the American Mathematical Society 2:345  (2022).

\bibitem{CDP} T. Champion, L. De Pascale, F. Prinari, \emph{$\Ga$-convergence and absolute minimizers for supremal functionals}, COCV ESAIM: Control, Optimisation and Calculus of Variations (2004), Vol. 10, 14-27.

\bibitem{Ch} B. M. Cherkas, \emph{Compactness in $L_\infty$ spaces}, Proceedings of the AMS 25:2, 347-350 (1970).

\bibitem{CKM} E. Clark, N. Katzourakis, B. Muha, \emph{Data assimilation for the Navier-Stokes equations through PDE-constrained optimisation in $\mathrm L^{\infty}$}, Nonlinearity 35:1 470 (2021).

\bibitem{CK2} E. Clark, N. Katzourakis, \emph{Second order generalised vectorial $\infty$-eigenvalue problems},  Proceedings of the Royal Society of Edinburgh, in Press (2024).

\bibitem{C} J. B. Conway, \emph{A course in Functional Analysis}, Graduate texts in mathematics, Springer, second edition, 2007.

\bibitem{D} B. Dacorogna, \emph{Direct Methods in the Calculus of Variations}, $2$nd Edition, Volume 78, Applied Mathematical Sciences, Springer, 2008.

\bibitem{Da} J. M. Danskin, \emph{The theory of max-min with applications}, J. SIAM Appl. Math.  Vol. 14, No. 4, 641-665 (1966). 

\bibitem{DiF} G. Di Fazio, \emph{$L^p$ estimates for divergence form elliptic equations with discontinuous coefficients}, Boll. Un. Mat. Ital. A (7) 10(2), 409-420 (1996).

\bibitem{EG} L.C. Evans, R. Gariepy, \emph{Measure theory and fine properties of functions}, Studies in advanced mathematics, CRC press, 1992.

\bibitem{F} H. Federer, \emph{Geometric measure theory}, Classics in Mathematics, Springer Berlin, Heidelberg, 1996.
 
\bibitem{FL} I. Fonseca, G. Leoni, \emph{Modern methods in the Calculus of Variations: ${\mathrm{L}}^p$ spaces}, Springer Monographs in Mathematics, 2007.

\bibitem{GT} D. Gilbarg, N. S. Trudinger, \emph{Elliptic Partial Differential Equations of Second Order}, 2nd edition, Springer-Verlag GmbH Germany, part of Springer Nature 2001.

\bibitem{HS} R. Hardt, L. Simon, \emph{Nodal sets for solutions of elliptic equations}, J. Differ. Geom., 30, 505–522 (1989).

\bibitem{JWY} R. Jensen, C. Wang, Y. Yu, \emph{Uniqueness and Nonuniqueness of Viscosity Solutions to Aronsson’s Equation}, Arch. Rational Mech. Anal. 190, 347-370 (2008).

\bibitem{K0} N. Katzourakis, An Introduction to Viscosity Solutions for Fully Nonlinear PDE with Applications to Calculus of Variations in $\mathrm L^{\infty}$, Springer Briefs in Mathematics, 2015, DOI 10.1007/978-3-319-12829-0.

\bibitem{K1} N. Katzourakis, \emph{Generalised solutions for fully nonlinear PDE systems and existence-uniqueness theorems}, Journal of Differential Equations 23, 641-686 (2017).

\bibitem{K2} N. Katzourakis, \emph{Generalised vectorial $\infty$-eigenvalue nonlinear problems for $L^\infty$ functionals}, Nonlinear Analysis 219, 112806 (2022).

\bibitem{KM} N. Katzourakis, R. Moser, \emph{Existence, Uniqueness and Structure of Second Order Absolute Minimisers}, Archives for Rational Mechanics and Analysis 231, Issue 3, 1615-1634 (2019).

\bibitem{KM2} N. Katzourakis, R. Moser, \emph{Variational problems in $\mathrm L^{\infty}$ involving semilinear second order differential operators}, ESAIM: COCV 29 (2023) 76.

\bibitem{KP} N. Katzourakis, E. Parini, \emph{The eigenvalue problem for the $\infty$-Bilaplacian}, Nonlinear Differential Equations and Applications NoDEA 24:68, (2017).

\bibitem{KPr} N. Katzourakis, T. Pryer, \emph{2nd order $\mathrm L^{\infty}$ Variational Problems and the $\infty$-Polylaplacian}, Advances in Calculus of Variations 13:2, 115-140 (2020).

\bibitem{KV} N. Katzourakis, E. Varvaruca,  \emph{An Illustrative Introduction to Modern Analysis}, CRC Press / Taylor \& Francis, 560 pages, 2018.

\bibitem{KT} H. Koch, D.  Tataru, \emph{Carleman estimates and unique continuation for second-order elliptic equations with nonsmooth coefficients}, Comm. Pure Appl. Math. 54 : 339-360 (2001).

\bibitem{KT2} H. Koch, D.  Tataru, \emph{Recent Results on Unique Continuation for Second Order Elliptic Equations}, In: Colombini, F., Zuily, C. (eds) Carleman Estimates and Applications to Uniqueness and Control Theory, Progress in Nonlinear Differential Equations and Their Applications vol. 46, Birkh\"auser, Boston, MA (2001). 

\bibitem{KZ} C. Kreisbeck, E. Zappale, \emph{Lower semicontinuity and relaxation of nonlocal $\mathrm L^{\infty}$-functionals}, Calculus of Variations and PDE 59 (4), 1-36 (2020).

\bibitem{KT} J. Kristensen, A. Taheri, \emph{Partial Regularity of Strong Local Minimizers in the Multi-Dimensional Calculus of Variations}, Arch. Rational Mech. Anal. 170 63-89 (2003).

\bibitem{MWZ} Q. Miao, C. Wang, Y. Zhou, \emph{Uniqueness of Absolute Minimizers for $\mathrm L^{\infty}$-Functionals Involving Hamiltonians $H(x,p)$}, Archive for Rational Mechanics and Analysis 223 (1), 141-198 (2017).

\bibitem{M1} R. Moser, \emph{Structure and classification results for the $\infty$-elastica problem}, American Journal of Mathematics 144, no. 5 (2022): 1299-1329.

\bibitem{MS} R. Moser, H. Schwetlick, \emph{Minimizers of a weighted maximum of the Gauss curvature}, Annals of Global Analysis and Geometry 41 (2012), 199-207.

\bibitem{O} M. I. Ostrovskii, \emph{Weak* sequential closures in Banach space theory and their applications, in: ``General Topology in Banach Spaces"}, ed. by T. Banakh and A. Plichko, Nova Science, New York, 2001, pp.\ 21-34.

\bibitem{PP} G. Papamikos, T. Pryer, \emph{A Lie symmetry analysis and explicit solutions of the two-dimensional $\infty$-Polylaplacian}, Studies in applied mathematics, Vol.\ 142, Issue1, 48-64 (2019).

\bibitem{PWZ} F. Peng, C. Wang, Y. Zhou, \emph{Regularity of absolute minimizers for continuous convex Hamiltonians}, Journal of Differential Equations, Volume 274 (2021), 1115-1164.

\bibitem{PZ} F. Prinari, E. Zappale, \emph{A Relaxation Result in the Vectorial Setting and Power Law Approximation for Supremal Functionals}, J Optim. Theory Appl. 186, 412-452 (2020).

\bibitem{RZ} A.N. Ribeiro, E. Zappale, \emph{Existence of minimisers for nonlevel convex functionals}, SIAM J. Control Opt., Vol. 52, No. 5,  (2014) 3341-3370.

\bibitem{R} W. Rudin, \emph{Functional Analysis}, second edition, McGraw-Hill, 1991.

\bibitem{S} Z. Sakellaris, \emph{Minimization of scalar curvature in conformal geometry}, Annals of Global Analysis and Geometry volume 51, pages 73-89 (2017).

\bibitem{T} J. Toland, \emph{The Dual of $\mathrm L^{\infty}(X,\mL,\la)$, Finitely Additive Measures and Weak Convergence: A Primer},  SpringerBriefs in Mathematics, 2020. 

\bibitem{W} H. Whitney, \emph{Analytic Extensions of Differentiable Functions Defined in Closed Sets}, Transactions of the American Mathematical Society Vol.\ 36, No.\ 1 (1934), 63-89.

\bibitem{Y} Y. Yu, \emph{$\mathrm L^{\infty}$ Variational Problems and Aronsson Equations}, Arch. Rational Mech. Anal. 182  153-180 (2006).




\end{thebibliography}

\end{document}